\documentclass[final,3p,times]{elsarticle}
\usepackage{hyperref}

\journal{}

\usepackage{amssymb}
\usepackage{amsthm}
\usepackage{amsmath}
\usepackage{graphicx}
\usepackage{float}
\usepackage{subfigure}
\usepackage{caption}

\begin{document}

\begin{frontmatter}

\title{Explicit and Implicit Kinetic Streamlined-Upwind Petrov Galerkin Method for Hyperbolic Partial Differential Equations}




\author{Ameya Dilip Jagtap}
\ead{ameya.aero@gmail.com, ameya@aero.iisc.ernet.in}

\author{S.V. Raghurama Rao}
\ead{raghu@aero.iisc.ernet.in}

\cortext[cor1]{Corresponding author}
\address{Department of Aerospace Engineering, Indian Institute of Science, Bengaluru, Karnataka, India-560012}

\begin{abstract}
A novel explicit and implicit Kinetic Streamlined-Upwind Petrov Galerkin (KSUPG) scheme is presented for hyperbolic equations such as Burgers equation and compressible Euler equations. The proposed scheme performs better than the original SUPG stabilized method in multi-dimensions. To demonstrate the numerical accuracy of the scheme, various numerical experiments have been carried out for 1D and 2D Burgers equation as well as for 1D and 2D Euler equations using Q4 and T3 elements. Furthermore, spectral stability analysis is done for the explicit 2D formulation. Finally, a comparison is made between explicit and implicit versions of the KSUPG scheme.  
\end{abstract}


\end{frontmatter}


\section{Introduction}
\label{Introduction} 
     Finite element method is one of the popular numerical methods for solving the partial differential equations numerically on digital computers.  The standard finite element method is well-suited to solve elliptic partial differential equations efficiently \cite{OCZT, SF, JTO}.  But this method produces oscillations for hyperbolic partial differential equations; for example, governing equations of convection dominated flows require additional stabilization for the standard finite element based discretization.  For such flows, many stabilized finite element methods are available in the literature like  Streamline-Upwind Petrov Galerkin (SUPG), Discontinuous-Galerkin method, Taylor Galerkin method, Galerkin Least-Squares method etc. The detailed discussion about these methods are given in \cite{OCZTF, PBGU, JSHTW, JJD, JJDB, TJC}.  Among them, SUPG method is one of the popular stabilized finite element methods used to solve high speed compressible flows governed by Euler equations \cite{TH, HT}.  This method introduces diffusion along the streamline direction of the flow which makes it stable.  Apart from diffusion requirement along the streamline, SUPG method needs additional diffusion across the high gradient regions especially for multidimensional case.  This diffusion can be controlled by using a shock capturing parameter which senses the shock region and adds the diffusion appropriately.  Many nonlinear discontinuity capturing terms are available in the literature \cite{TO, TV, TS}).  

In the finite volume methods, kinetic schemes (also known as Boltzmann scheme) are interesting alternatives for the popular Riemann solvers.  Development of these schemes are based on the fact that one can recover the Euler equations by applying a suitable moment method strategy to the Boltzmann equation. The Boltzmann equation is given by 
\begin{equation}  
\frac{\partial f}{\partial t} + v. \nabla f = \left( \frac{d f}{d t}\right)_{\mathcal{C}}
\end{equation} 
where $f$ and $v$ are velocity distribution function and molecular velocity respectively. The right hand side is the collision term and left hand side consists of an unsteady term and a convection term.  The well-known BGK model simplifies the collision term and converts the otherwise integro-differential equation to a partial differential equation with a relaxation source term 
\cite{BGK}.  Using an operator splitting strategy by which the solution of the Boltzmann equation is split into a convection step and a collision step and further employing an instantaneous relaxation to equilibrium in the collision step leads to a simplification which is often used in Boltzmann schemes.  The moments of the resulting Boltzmann equation then leads to Euler equations of gas dynamics, with the equilibrium distribution $f^{M}$ being a Maxwellian.  The Euler equations can then be written in the following moment form.  
\begin{equation} 
\left<\Psi \left(\frac{\partial f}{\partial t} + v. \nabla f = 0 \right), \ 
f = f^{M} \right> 
\end{equation} 
Here $<.>$ is an appropriate moment and $\Psi$ is the moment function vector. One can also obtain the Burgers equation using same procedure stated above by defining an appropriate equilibrium distribution function. The advantage of this procedure is, instead of dealing with a nonlinear hyperbolic system of equations (Euler equations) we are dealing with a linear scalar equation (Boltzmann equation without collision term). There are many kinetic schemes available in the literature like Beam scheme of Sanders \& Prendergast \cite{SP9}, the method of Rietz \cite{RDR}, the Equilibrium Flux Method of Pullin \cite{DIP}, Kinetic Flux Vector Splitting (KFVS) method of Deshpande \cite{SMD,SMJC}, the compactly supported distribution based methods of Kaniel \cite{Kaniel} and Perthame 
\cite{Perthame}, the Peculiar Velocity based Upwind (PVU) method of Raghurama Rao \& Deshpande \cite{SVRRao_SMD} and the BGK scheme of Prendergast \& Xu 
\cite{PK8}.  These methods were developed in the framework of finite difference or finite volume methods.  The application of finite element methods in the framework of kinetic schemes is of currently ongoing interest, with some of the works in this category being due to Deshpande \& Pironneau \cite{SMD_Pironneau} and Deshpande \cite{SSD}, Yu \& Dai \cite{Yu_Dai}, Khobalatte \& Leyland 
\cite{Khobalatte_Leyland}, Tang \& Warnecke \cite{Tang_Warnecke}, Liu ans Xu 
\cite{Liu_Xu}, Ren {\em et al.} \cite{Ren_Hu_Xiong_Qiu}, Gassner \cite{GJG}.  

In this paper an attempt has been made to take advantage of the strategy of kinetic schemes for developing an efficient SUPG scheme in the framework of Boltzmann schemes. Along with this novel scheme (KSUPG), we have also developed a simple shock capturing parameter which senses the jump inside the element for 2D Euler equations.  However, unlike the traditional SUPG method, the shock capturing parameter is needed only for 2D Euler equations (not in one dimension) and not even for 2D Burgers equation. Constructing the stabilization parameter $\tau$ (which is the intrinsic time scale) in mutidimensional SUPG framework is not a trivial task. Many methods are available in the literature \cite{BT,TO}. But, in the proposed KSUPG scheme $\tau$ is defined for both scalar and vector equations simply from the linear scalar formulation. The efficiency of the new scheme is demonstrated by solving various test cases. This paper is arranged as follows.  In section 2 governing equations for high speed flows (Euler equations) and scalar Burgers equation are given.  Section 3 and 4 give the 1D and 2D explicit KSUPG weak formulation for both Burgers equation and Euler equations. In section 5, a simple shock capturing parameter is introduced. Section 6 explains the spectral stability analysis for explicit KSUPG scheme. An implicit KSUPG formulation for 1D and 2D Euler equations is given in section 7 followed by section 8 where various numerical test cases for both explicit and Implicit formulation are solved. Before ending section 8, comparison for explicit and implicit KSUPG schemes is made based on the number of iterations required to bring down the residue below the specified tolerance limit, the  computational cost and the sparsity pattern of the coefficient matrix of global system of equations. 
			  
\section{Governing Equations}

The governing equations are for the inviscid compressible flows, given by Euler equations as
\begin{equation}\label{qeul}
\frac{\partial \mathbf{U}}{\partial t}+ \frac{\partial \mathbf{G}_i}{\partial x_i} = 0  \ \  \in \,\, \Omega \times [0, \,\,T]
\end{equation}
where 
$\mathbf{U}~=~[\rho,\,\rho u_1,\,\rho u_2,\,\rho u_3,\, \rho E]^T$ is the conserved variable vector and $\mathbf{G}_i = [\rho u_i , \, \delta_{i1}p+ \rho u_1u_i ,\, \delta_{i2}p+ \rho u_2u_i , \delta_{i3}p+ \rho u_3u_i ,\,\, p u_i+ \rho u_i E]^T$ is the inviscid  flux vector. $\rho, u_1, u_2,u_3, E, p$ are density, velocity components in $x, y$ and $z$ directions, total energy and pressure respectively and $\delta_{ij}$ is a kronecker delta.  Total energy is given by
$$E = \frac{p}{\rho (\gamma-1)} + \frac{1}{2} (u_1^2+ u_2^2+ u_3^2)$$
Note that $\mathbf{A}_i = \frac{\partial \mathbf{G}_i}{\partial \mathbf{U}}$ is an inviscid flux jacobian matrix for the domain $\Omega \in \mathbb{R}^\mathcal{D}$ (where $\mathcal{D}$ is the spatial dimension) with boundary $\Gamma = \Gamma_D \cup \Gamma_N$ and final time is given by $T \in \mathbb{R}^{+}$.  As the eigenvalues of $\mathbf{A}_i$ are real and eigenvectors are linearly independent, the system of equations is hyperbolic.  Beyond being hyperbolic, these equations are nonlinear and produce shock waves, expansion waves and contact discontinuities which need to be resolved in numerical simulations.   

We also consider a scalar hyperbolic conservation law as 
\begin{equation}
\frac{\partial u}{\partial t}+ \frac{\partial g_{i}(u)}{\partial x_i} = 0  \ \  \in \,\, \Omega \times [0, \,\,T]
\end{equation} 
where $u$ is the conserved variable.  The fluxes $g_{i}(u)$ can be linear or nonlinear.  One example is the inviscid Burgers equation in the the fluxes are nonlinear and produce shock waves and expansion waves.  

\vspace{-6pt}

\section{One Dimensional Explicit KSUPG Weak Formulation}

\vspace{-2pt}
The standard Galerkin finite element approximation for molecular velocity distribution function is
\begin{equation}
f \approx f^h = \sum_{\forall i} N_i^h f_i
\end{equation}
where the domain is divided into $n_{el}$ elements.
\begin{equation}
\Omega = \bigcup_{\forall i}^{n_{el}} \Omega_i \ \ \textrm{and} \ \ \Omega_i  \cap \Omega_j = \O , \,\,\, \forall i \neq j
\end{equation}
Defining the appropriate test and trial functions spaces as $\mathcal{V}^h~=~\{N^h~\in~\mathcal{H}^{1} \linebreak \, \textrm{and} \, N^h  = 0 \,\,\text{on}\,\, \Gamma_D \}$ and $\mathcal{S}^h  =  \{  f^h \in \mathcal{H}^{1} \, \textrm{and} \, f^h = f^h_D \,\,\text{on}\,\, \Gamma_D \}  $ where $\Gamma_D$ is the Dirichlet boundary, the weak formulation is written as, find $f^h \in \mathcal{S}^h $ such that $\forall \, N^h \in \mathcal{V}^h$
\begin{align}
 \sum_{i =1}^{n_{el}} \int_{\Omega_i} \left( \frac{\partial f^h}{\partial t} + v \frac{\partial f^h}{\partial x} \right)  \left( N_i^h +  \tau v \frac{dN_i^h}{dx}\right) \, d\Omega_i  = 0
\end{align}
where $\tau = h/(2|v|)$.
The global system of equations are obtained as
\begin{align}
& \int_{\Omega} \left( \frac{\partial f^h}{\partial t} + v \frac{\partial f^h}{\partial x} \right)  \left( N^h +  \frac{h}{2|v|} v \frac{dN^h}{dx}\right) \, d\Omega   = 0 \nonumber 
\\ 
\textrm{or} \nonumber 
\\  & \int_{\Omega} \left(N^h\right)^T (N^h)  \, d\Omega \frac{d f}{d t}  + \int_{\Omega} \left(N^h\right)^T \left(\frac{\partial N^h}{\partial x} \right)  \, d\Omega\,\, vf \nonumber \\ & + \frac{h}{2}  \int_{\Omega} \left(\frac{\partial N^h}{\partial x}\right)^T  \left(\frac{\partial N^h}{\partial x}  \right) \, d\Omega \,\, \text{sign}(v) \, v f +\int_{\Gamma_N} \frac{\partial f^h}{\partial x}  d\Gamma_N= 0 
\end{align}
where basis functions $N^h \in C^0(\Omega)$. It is important to note that, the test function is enriched with additional term which is multiplied only with the convection term. That gives a required diffusion term.

In matrix form,
\begin{equation}
M \frac{d f}{d t}  +C  vf + \frac{h}{2} D \ \text{sign}(v) \, v f  + f_N = 0 
\end{equation}
where Mass matrix $M$, Convection matrix $C$, Diffuion matrix $D$ and Neumann boundary condition $f_N$ are given by
\begin{align}
M &=  \int_{\Omega} \left(N^h\right)^T (N^h)   \, d\Omega \nonumber 
\\ C &= \int_{\Omega} \left(N^h\right)^T \left(\frac{\partial N^h}{\partial x} \right)   \, d\Omega \nonumber 
\\  D  & = \int_{\Omega} \left(\frac{\partial N^h}{\partial x}\right)^T  \left(\frac{\partial N^h}{\partial x}  \right)   \, d\Omega \nonumber 
\\ f_N & = \int_{\Gamma_N} \frac{\partial f^h}{\partial x}  d\Gamma_N  \nonumber 
\end{align}
All the integrals are evaluated with full Gauss-Quadrature integration.  Taking moments with the suitable moment function vector $\Psi$
\begin{align}\label{be4}
\left<\Psi, M \frac{d f}{d t}  + C vf + \frac{h}{2} D \text{sign}(v) \, v f   + f_N  \right>& = 0 \nonumber  
\\ \textrm{or} \ M \frac{d <\Psi, f>}{d t}  +  C <\Psi, vf>+ \frac{h}{2} D <\Psi,  \text{sign}(v) \, v f >  + < f_N > & = 0
\end{align}
Equation \eqref{be4} is the semi-discrete weak formulation.

\subsection{One Dimensional Burgers Equation}
The 1D Burgers equation is given by

\begin{equation}
\frac{\partial u}{\partial t}+ \frac{\partial g(u)}{\partial x}= 0  \ \  \in \,\, \Omega \times [0, \,\,T]
\end{equation} 
For the sake of convenience, we write the flux $g(u) = \frac{1}{2} u^{2}$ as 
$g(u) = c u$ with $c= \frac{1}{2} u$.  In case of one dimensional Burgers equation $\Psi = 1$ and Maxwellian distribution function to recover the Burgers equation as a moment from the Boltzmann equation is given by $$f^M = u \left(\frac{\beta}{\pi}\right)^{\frac{\mathcal{D}}{2}} e^{-\beta (v -c)^2}$$ where $c$ is constant to recover the linear convection equation and is a function of $u$, {\em i.e.}, $c = u/2$ to recover the inviscid  Burgers equation.  For one dimensional problem $\mathcal{D} = 1$. Let us now evaluate the terms for the case of one dimensional Burgers equation. 

\begin{align}
<\Psi, f^M > & = \int_{-\infty}^{\infty} f^M \, dv = u
\\ <\Psi,  v f^M > & =  \int_{-\infty}^{\infty} v f^M \, dv =cu
\end{align}
\begin{align}\label{kd3}
< \Psi,\text{sign}(v) \, v  f^M > & = \int_{-\infty}^{\infty} \text{sign}(v) \, 
v f^M \, dv \nonumber
\\ & =  -\int_{-\infty}^{0}  v f^M \, dv  +  \int_{0}^{\infty}  v f^M \, dv \nonumber
\\ & =  cu\,\text{erf}(s) + \frac{u}{\sqrt{\pi \beta}} e^{-s^2} 
\end{align}
\begin{align}
<\Psi,f_N > & = \int_{\Gamma_N} \frac{\partial <\Psi, f^h>}{\partial x}   d\Gamma_N  = \int_{\Gamma_N} \frac{\partial u}{\partial x}   d\Gamma_N 
\end{align}
where $s = u\sqrt{\beta}/2$ and $\beta = 1$. Note that, since no energy equation is involved (no pressure and temperature terms) so, $\beta$ is just a constant value calculated from the standard Maxwellian distribution function. Moments of  last expression lead to the Neumann boundary condition in macroscopic variable $u$.

Substituting these values in equation \eqref{be4}, we get
\begin{equation}\label{ksupg1}
M \frac{d u}{d t}  + C cu  + \frac{h}{2} D \left( cu\,\text{erf}(s) + \frac{u}{\sqrt{\pi \beta}} e^{-s^2}  \right) + u_N   = 0
\end{equation}
where $u_N$ is the Neumann boundary condition in macroscopic variable.

\subsubsection{Temporal Discretization:}   
In this work finite difference approach is adopted for temporal discretization using $\theta$ method as 
\begin{align}
M \frac{u^{n+1}-u^{n}}{\Delta t}  & + (1-\theta)\left( Cc^nu^n  + \frac{h}{2}  D\left[ c^nu^n \,\text{erf}(s)  + \frac{u^n}{\sqrt{\pi \beta}} e^{-s^2}  \right]   \right) \nonumber
\\ & + \theta \left( Cc^nu^n   + \frac{h}{2} D \left[ c^{n+1}u^{n+1}\,\text{erf}(s) + \frac{u^{n+1}}{\sqrt{\pi \beta}} e^{-s^2}  \right]  \right) +u_N = 0
\end{align}
Thus, $\theta = 0$ gives an explicit method and $\theta = 1$ gives an implicit method. Semi-Implicit methods can be obtained with $0 < \theta < 1$.  For example, $\theta = 1/2$ gives Crank-Nicolson method. 

\subsection{One Dimensional Euler Equations}
The 1D Euler equations are given by 
\begin{equation}
\frac{\partial \mathbf{U}}{\partial t} + \frac{\partial \mathbf{G}}{\partial x} =0
\end{equation}
where $$\mathbf{U} = \left \{ 
\begin{array}{c}
\rho \\
\rho u \\
\rho E \\
\end{array} \right  \}, \ \  \mathbf{G} = \left \{ 
\begin{array}{c}
\rho u \\
p+\rho u^2 \\
pu+\rho u E \\
\end{array} \right  \}$$ are the solution vector and the flux vector respectively.  
For recovering the 1D Euler equations as moments of the Boltzmann equation, the Maxwellian distribution function is given by
\begin{equation}
f^M = \frac{\rho}{I_o} \left( \frac{\beta}{\pi} \right)^{\frac{\mathcal{D}}{2}}  
e^{\left(-\beta(v^2-u^2)- \frac{I}{I_o}\right)}
\end{equation}
where $v$ is the molecular velocity, $I$ is the internal energy variable corresponding to the non-translational degrees of freedom  and $I_o$ is the average internal energy corresponding to the non-translational degrees of freedom which is given by 
\begin{equation}
I_o = \frac{3-\gamma}{2(\gamma-1)}RT
\end{equation}
and $\gamma$ is the ratio of specific heats.

For 1D Euler equations, moment function vector $\Psi$ is defined as
\begin{equation}
\Psi = \left \{ 
\begin{array}{c}
1 \\
v \\
I+\frac{v^2}{2}\\
\end{array} \right  \} 
\end{equation}
With the Maxwellian $f = f^M$ the other terms in weak formulation given by equation \eqref{be4} can be obtained as
\begin{equation}
<\Psi, f^M>  =\int_0^{\infty} dI  \int_{-\infty}^{\infty} \Psi f^M dv = \left \{ 
\begin{array}{c}
\rho \\
\rho u \\
\rho E \\
\end{array} \right  \}  = \mathbf{U}
\end{equation}
\begin{equation}
 <\Psi, vf^M>  = \int_0^{\infty} dI  \int_{-\infty}^{\infty} \Psi v f^M dv = \left \{ 
\begin{array}{c}
\rho u \\
p+\rho u^2 \\
pu+\rho u E \\
\end{array} \right  \} =  \mathbf{G}
\end{equation}

\begin{align}\label{ks1e}
 < \Psi,\text{sign}(v) \, v f^M > & = \int_0^{\infty} dI   \int_{-\infty}^{\infty} \text{sign}(v) \, v f^M \, dv \nonumber
\\ & =  -\int_0^{\infty} dI   \int_{-\infty}^{0}  v f^M\, dv  + \int_0^{\infty} dI   \int_{0}^{\infty}  v f^M \, dv 
\end{align}
Taking the first integral term
\begin{align}
\int_0^{\infty} dI   \int_{-\infty}^{0}  v f^M\, dv   & =\left \{ 
\begin{array}{c}
\rho u A^{-} -\rho B\\
(p+\rho u^2) A^{-} -\rho u B\\
\left( p+\rho E \right)u A^{-} - \rho\left( \frac{p}{2\rho} +E \right)B \\
\end{array} \right  \}  
\end{align}
Similarly, second integral term can be evaluated as
\begin{align}
\int_0^{\infty} dI   \int_{0}^{\infty}  v f^M\, dv   & =\left \{ 
\begin{array}{c}
\rho u A^{+} +\rho B\\
(p+\rho u^2) A^{+} +\rho u B\\
\left( p+\rho E \right)u A^{+} + \rho\left( \frac{p}{2\rho} +E \right)B \\
\end{array} \right  \} 
\end{align}
where 
$$A^{\pm} = \frac{1 \pm \text{erf}(s)}{2} \ \textrm{and} \ B = \frac{1}{2\sqrt{\pi \beta}} 
e^{-s^2}$$
Here, $s = u\sqrt{\beta}$, $\beta = 1/2RT$.  Substituting these values in equation \eqref{ks1e} and then simplifying we get
\begin{align}
 <\Psi, \text{sign}(v) \, v f^M > & = \left \{ 
\begin{array}{c}
\rho u \,\text{erf}(s) + \frac{\rho}{\sqrt{\pi \beta}} e^{-s^2}\\
(p+\rho u^2) \,\text{erf}(s) + \frac{\rho u}{\sqrt{\pi \beta}} e^{-s^2}\\
\left( p+\rho E \right)u \,\text{erf}(s) + \rho\left( \frac{p}{2\rho} +E \right) \frac{1}{\sqrt{\pi \beta}} e^{-s^2}  \\
\end{array} \right  \} 
\end{align}
Substituting these values in \eqref{be4} we get
\begin{align}\small 
& M \frac{d  }{d t} \left \{ 
\begin{array}{c}
\rho \\
\rho u \\
\rho E \\
\end{array} \right  \}  
+ C\left \{ 
\begin{array}{c}
\rho u \\
p+\rho u^2 \\
pu+\rho u E \\
\end{array} \right  \} \nonumber 
\\ 
& + \frac{h}{2}  D\left \{ 
\begin{array}{c}
\rho u \,\text{erf}(s) + \frac{\rho}{\sqrt{\pi \beta}} e^{-s^2}\\
(p+\rho u^2) \,\text{erf}(s) + \frac{\rho u}{\sqrt{\pi \beta}} e^{-s^2}\\
\left( p+\rho E \right)u \,\text{erf}(s) + \rho\left( \frac{p}{2\rho} +E \right) \frac{1}{\sqrt{\pi \beta}} e^{-s^2}  \\
\end{array} \right  \} +u_N  & = 0
\end{align} 

\subsubsection{Temporal Discretization:}  
Temporal discretization is done using $\theta$ Method with $\theta = 0$ as
\begin{align}
M \frac{\mathbf{U}^{n+1}-\mathbf{U}^n }{\Delta t} + \left(C \mathbf{A} \mathbf{U}^n  + \frac{h}{2} D < \Psi,\text{sign}(v) \, v f^M >^n  \right) + u_N& = 0
\end{align}
In the above discretized form the flux vector $\mathbf{G}$ is written as $\mathbf{AU}$ where $\mathbf{A}$ is the flux Jacobian matrix given by
\begin{equation}
\mathbf{A}  = \left [ 
\begin{array}{ccc}
0&1 &0 \\
\left( \frac{\gamma-3}{2}\right) u^2& (3-\gamma)u & (\gamma-1) \\
(\gamma-1)u^3 - \gamma u E  & \gamma E- \frac{3}{2}(\gamma -1) u^2 & \gamma u\\
\end{array} \right  ] 
\end{equation}

\subsection{Linearization and Iterative Solver}
The global nonlinear fully discretized equation of the form $K(u)u = \phi$ can be linearized by using Picard iteration technique as 
$$K(u^{n})u^{n+1} = \phi $$   
Then, the linearized system of equations is solved using bi-conjugate gradient stabilized method.

\section{Two Dimensional Explicit KSUPG Weak Formulation}

The standard Galerkin finite element approximation for molecular velocity distribution function is
\begin{equation}
f \approx f^h = \sum_{\forall i} N_i^h f_i
\end{equation}
where the domain is divided into $n_{el}$ elements.
\begin{equation}
\Omega = \bigcup_{\forall i}^{n_{el}} \Omega_i \ \ \textrm{and} \ \ \Omega_i  \cap \Omega_j = \O , \,\,\, \forall i \neq j
\end{equation}
Defining the test and trial functions spaces as $\mathcal{V}^h = \{ N^h \in \mathcal{H}^{1} \, \textrm{and} \linebreak N^h = 0 \,\,\text{on}\,\, \Gamma_D \}$ and $\mathcal{S}^h  =  \{  f^h \in \mathcal{H}^{1} \, \textrm{and} \, f^h = f^h_D \,\,\text{on}\,\, \Gamma_D \}  $ where $\Gamma_D$ is the Dirichlet boundary, the weak formulation is written as, find $f^h \in \mathcal{S}^h $ such that $\forall \, N^h \in \mathcal{V}^h$
\begin{align}
 \sum_{i =1}^{n_{el}} \int_{\Omega_i} \left( \frac{\partial f^h}{\partial t} + v_1 \frac{\partial f^h}{\partial x} + v_2 \frac{\partial f^h}{\partial y} \right)  \left( N_i^h+  \left[ \tau_1 v_1 \frac{dN_i^h}{dx} + \tau_1 v_2 \frac{dN_i^h}{dy}   \right]  \right)  \, d\Omega_i  = 0
\end{align}
where $\tau_1 = h/(2|v_1|)$ and $\tau_2 = h/(2|v_2|)$.
The global system of equations are obtained as
\begin{equation}  
\int_{\Omega} \left( \frac{\partial f^h}{\partial t} + v_1 \frac{\partial f^h}{\partial x} + v_2 \frac{\partial f^h}{\partial y} \right)  \left( N^h + \left[ \tau_1 v_1 \frac{dN^h}{dx} + \tau_1 v_2 \frac{dN^h}{dy}   \right]   \right)  \, d\Omega   = 0 
\end{equation}
or 
\begin{align}
& \int_{\Omega} (N^h)^T (N^h)  \, d\Omega \frac{d f}{d t}  +\int_{\Omega} (N^h)^T \left( \frac{\partial N^h}{\partial x} \right)  \, d\Omega \, v_1f +  \int_{\Omega} (N^h)^T \left(\frac{\partial N^h}{\partial y}\right)   \, d\Omega\,v_2 f \nonumber
\\ & + \frac{h}{2}  \int_{\Omega} \left(\frac{\partial N^h}{\partial x}\right)^T  \left(\frac{\partial N^h}{\partial x}\right)   \, d\Omega \,\text{sign}(v_1) \, v_1 f \nonumber \\
& + \frac{h}{2}  \int_{\Omega} \left(\frac{\partial N^h}{\partial x} \right)^T \left(\frac{\partial N^h}{\partial y}\right)   \, d\Omega\,\text{sign}(v1) \, v2 f\nonumber
\\ & + \frac{h}{2}  \int_{\Omega} \left(\frac{\partial N^h}{\partial x} \right)^T \left(\frac{\partial N^h}{\partial y} \right)  \, d\Omega \text{sign}(v2) \, v1 f \nonumber \\ 
& + \frac{h}{2}  \int_{\Omega} \left(\frac{\partial N^h}{\partial y}\right)^T  \left(\frac{\partial N^h}{\partial y}\right)   \, d\Omega \,\text{sign}(v2) \, v2 f \nonumber
\\ & +\int_{\Gamma_N} \frac{\partial f^h}{\partial n}   d\Gamma_N= 0 
\end{align}
where basis functions $N^h \in C^0(\Omega)$. Again, the enriched terms present in the test function are multiplied only with convective terms which gives diffusion terms in $x$, $y$ directions and cross-diffusion terms in $x-y$ directions. \\

\noindent In matrix form,
\begin{align}
 M \frac{d f}{d t} & + C_xv_1f  +C_y v_2f + \frac{h}{2} D_x\text{sign}(v_1) \, v_1 f   + \frac{h}{2}D_{xy}  (\text{sign}(v_1) \, v_2f+ \text{sign}(v_2) \, v_1f) \nonumber
\\ & + \frac{h}{2} D_y\text{sign}(v_2) \, v_2 f  +f_N= 0 
\end{align}
where
\begin{align}
M & =  \int_{\Omega} (N^h)^T (N^h)  \, d\Omega \nonumber 
\\ C_x & = \int_{\Omega} (N^h)^T \left(\frac{\partial N^h}{\partial x}\right) \nonumber 
\\ C_y & = \int_{\Omega} (N^h)^T \left(\frac{\partial N^h}{\partial y} \right) \nonumber 
\\ D_x & = \int_{\Omega} \left(\frac{\partial N^h}{\partial x}\right)^T  \left(\frac{\partial N^h}{\partial x}\right)   \, d\Omega \nonumber 
\end{align}
\begin{align}
 D_{xy} & = \int_{\Omega} \left(\frac{\partial N^h}{\partial x} \right)^T \left(\frac{\partial N^h}{\partial y} \right)  \, d\Omega \nonumber 
\\ D_y & = \int_{\Omega} \left(\frac{\partial N^h}{\partial y}\right)^T  \left(\frac{\partial N^h}{\partial y}\right)   \, d\Omega \nonumber 
\\ f_N & = \int_{\Gamma_N} \frac{\partial f^h}{\partial n}   d\Gamma_N \nonumber 
\end{align}
All integrals are evaluated using full Gauss-Quadrature integration.  Taking moments
\begin{align}
&\left< \Psi, M \frac{d f}{d t}  + C_xv_1f  +C_y v_2f + \frac{h}{2} D_x \text{sign}(v_1) \, v_1 f   + \frac{h}{2} D_{xy} (\text{sign}(v_1) \, v_2f+ \text{sign}(v_2) \, v_1f)   \right. \nonumber
\\ &  \left.+ \frac{h}{2} D_y \text{sign}(v_2) \, v_2 f    \right> + <f_N>= 0 \nonumber
 \\ &  M \frac{d <\Psi,f>}{d t}  + C_x<\Psi,v_1f>  + C_y<v_2f> + \frac{h}{2}D_x  <\Psi,\text{sign}(v_1) \, v_1 f >  \nonumber
\\& + \frac{h}{2} D_{xy}(<\Psi,\text{sign}(v_1) \, v_2f >+ <\Psi,\text{sign}(v_2) \, v_1f> )  + \frac{h}{2}  D_y <\Psi, \text{sign}(v_2) \, v_2 f>\nonumber
\\ & + <\Psi,f_N> = 0  \nonumber
\end{align}
Now lets evaluate these moments for 2D Burgers equation and 2D Euler equations.

\subsection{2D Burgers Equation}
The 2D Burgers equation is given by 
\begin{equation} 
\frac{\partial u}{\partial t}+ \frac{\partial g_1 (u)}{\partial x}+\frac{\partial  g_2 (u)}{\partial y}= 0 
\end{equation} 
is written as 
\begin{equation}
\frac{\partial u}{\partial t}+ \frac{\partial c_1 u}{\partial x}+\frac{\partial  c_2 u}{\partial y}= 0 
\end{equation}
where $c_1$ and $c_2$ can be functions of $u$ for obtaining nonlinearity or can be constants for keeping them as linear.  Maxwellian distribution function for recovering the 2D Burgers Equation as a moment of the Boltzmann equation is given by
\begin{equation}
f^M = u \left( \frac{\beta}{\pi} \right)^{\frac{\mathcal{D}}{2}} 
e^{(-\beta(v_1-c_1)^2-\beta(v_2-c_2)^2)}
\end{equation}
Here, $\Psi = 1$ and value of $\beta$ is fixed as unity.

For 2D Burgers equation one can evaluate the integrals as
\begin{align}
<\Psi,  f^M > &  =\int_{-\infty}^{\infty} \int_{-\infty}^{\infty} f^M \, dv_1  dv_2= u
\end{align}
\begin{align}
<\Psi,  v_1 f^M > &  =\int_{-\infty}^{\infty}\int_{-\infty}^{\infty} v_1 f^M \, dv_1  dv_2 = c_1u
\end{align}
\begin{align}
<\Psi,  v_2 f^M > &  = \int_{-\infty}^{\infty}\int_{-\infty}^{\infty} v_2 f^M \, dv_1  dv_2 = c_2u
\end{align}
\begin{align}
<\Psi,  f_N > &  = \int_{\Gamma_N} \frac{\partial <\Psi,f^h>}{\partial n}   d\Gamma_N  = \int_{\Gamma_N} \frac{\partial u}{\partial n}   d\Gamma_N 
\end{align}
\begin{align}
<\Psi,\text{sign}(v_1)\,v_1 f^M> &= \int_{-\infty}^{\infty} \int_{-\infty}^{\infty} \text{sign}(v_1)\,v_1 f^M \, dv_1  dv_2 \nonumber \\ & = \sqrt{\frac{\beta}{\pi}} u \left[ \frac{e^{-s_1^2}}{\pi} + c_1\text{erf}(s_1)  \right]   
 \\ <\Psi, \text{sign}(v_2)\ ,v_2 f^M> &=\int_{-\infty}^{\infty}\int_{-\infty}^{\infty} \text{sign}(v_2)\,v_2 f^M \, dv_1  dv_2 \nonumber \\ & =\sqrt{\frac{\beta}{\pi}} u \left[ \frac{e^{-s_2^2}}{\pi} + c_2\text{erf}(s_2)  \right] 
\\ <\Psi , \text{sign}(v_2)\,v_1 f^M>  & =\int_{-\infty}^{\infty}\int_{-\infty}^{\infty} \text{sign}(v_2)\,v_1 f^M \, dv_1  dv_2 \nonumber \\ & = \sqrt{\frac{\beta}{\pi}} uc_1 \text{erf}(s_2) 
\\ <\Psi , \text{sign}(v_1)\,v_2 f^M>  & = \int_{-\infty}^{\infty}\int_{-\infty}^{\infty} \text{sign}(v_1)\,v_2 f^M \, dv_1  dv_2 \nonumber \\ & =   \sqrt{\frac{\beta}{\pi}} uc_2 \text{erf}(s_1) 
\end{align}

\subsection{2D Euler Equations}
The 2D Euler equations are given by 
\begin{equation}
\frac{\partial \mathbf{U}}{\partial t} + \frac{\partial \mathbf{G}_1}{\partial x} + \frac{\partial \mathbf{G}_2}{\partial x} =0
\end{equation}
where $$\mathbf{U} = \left \{ 
\begin{array}{c}
\rho \\
\rho u_1 \\
\rho u_2
\rho E \\
\end{array} \right  \}, \ \  \mathbf{G}_1 = \left \{ 
\begin{array}{c}
\rho u_1 \\
\rho u_1 u_2 \\
p+\rho u_1^2 \\
pu_1+\rho u_1 E \\
\end{array} \right  \}, \ \  \mathbf{G}_2 = \left \{ 
\begin{array}{c}
\rho u_2 \\
p+\rho u_2^2 \\
\rho u_1 u_2 \\
pu_2+\rho u_2 E \\
\end{array} \right  \}$$ are the solution vector and the flux vectors in $x$ and $y$ directions resptively.

In case of 2D Euler equations, the Maxwellian distribution function is given as
\begin{equation}
 f^M = \frac{\rho}{I_o} \frac{\beta}{\pi} e^{\left(-\beta(v_1-u_1)^2-\beta(v_2-u_2)^2- \frac{I}{I_o}\right)}
\end{equation}
where $v_1$ and $v_2$ are molecular velocities in $x$ and $y$ directions and $I_o$ is defined as
\begin{equation}
I_o = \frac{2-\gamma}{(\gamma-1)}RT
\end{equation}
The vector $\Psi$ is defined as
\begin{equation}
\Psi = \left \{ 
\begin{array}{c}
1 \\
v_1 \\
v_2 \\
I+\frac{v_1^2+v_2^2}{2}\\
\end{array} \right  \} 
\end{equation}
and  $\beta = \frac{1}{2RT}$. Integrals are evaluated as
\begin{align}
<\Psi,  f^M > &  =\int_0^{\infty} dI \, \int_{-\infty}^{\infty}  dv_1 \int_{-\infty}^{\infty} dv_2 \Psi f^M \nonumber
\\ & = \left \{ 
\begin{array}{c}
\rho  \\
\rho u_1 \\
\rho  u_2\\
 \rho  E \\
\end{array} \right  \}  = \mathbf{U}
\end{align}
\begin{align}
<\Psi,  v_1 f^M > &  =\int_0^{\infty} dI \, \int_{-\infty}^{\infty}  dv_1 \int_{-\infty}^{\infty} dv_2 \Psi v_1 f^M \nonumber
\\ & = \left \{ 
\begin{array}{c}
\rho u_1 \\
p+\rho u_1 \\
\rho u_1 u_2\\
pu_1+ \rho u_1 E \\
\end{array} \right  \}  = \mathbf{G}_1
\end{align}
Flux vector $G_1$ is further decompose by using homogenity property $\mathbf{G}_1 = \mathbf{A}_1 \mathbf{U}$ where
\begin{equation}
\mathbf{A}_1  =   \left [ 
\begin{array}{cccc}
0&1&0 &0 \\
-u_1^2 + \frac{\gamma-1}{2} (u_1^2+u_2^2) & (3-\gamma) u_1&-(\gamma-1)u_2 &\gamma-1 \\
-u_1 u_2&u_2&u_1  &0 \\
-(\gamma E - (\gamma-1)(u_1^2+u_2^2)) u_1 &\gamma E - \frac{\gamma-1}{2}(2u_1^2 (u_1^2+u_2^2)) &-(\gamma-1)u_1 u_2& \gamma u_1\\
\end{array} \right ]  \nonumber
\end{equation}
\begin{align}
<\Psi,  v_2 f^M > &  = \int_0^{\infty} dI \, \int_{-\infty}^{\infty}  dv_1 \int_{-\infty}^{\infty} dv_2 \Psi v_2 f^M \nonumber
\\ & = \left \{ 
\begin{array}{c}
\rho u_2 \\
\rho u_1 u_2\\
p+\rho u_2\\
pu_2+ \rho u_2 E \\
\end{array} \right  \}  = \mathbf{G}_2
\end{align}
Similarly, flux vector $G_2$ is further decompose by usign homogenity property $\mathbf{G}_2 = \mathbf{A}_2 \mathbf{U}$ where
\begin{equation}
\mathbf{A}_2  =  \left [ 
\begin{array}{cccc}
0&0&1 &0 \\
-u_1 u_2&u_1&u_2  &0 \\
-u_2^2 + \frac{\gamma-1}{2} (u_1^2+u_2^2) &-(\gamma-1)u_1 & (3-\gamma) u_2&\gamma-1 \\
-(\gamma E - (\gamma-1)(u_1^2+u_2^2)) u_2 &-(\gamma-1)u_1 u_2&\gamma E - \frac{\gamma-1}{2}(2u_2^2 (u_1^2+u_2^2)) & \gamma u_2\\
\end{array} \right ]  \nonumber
\end{equation}

\begin{align}
<\Psi,  f_N > &  = \int_{\Gamma_N} \frac{\partial <\Psi,f^h>}{\partial n} N^h  d\Gamma_N \nonumber
\\ & = \int_{\Gamma_N} \frac{\partial u}{\partial n} N^h  d\Gamma_N 
\end{align}
\begin{align}
<\Psi,\text{sign}(v_1)\,v_1 f^M> &= \left \{ 
\begin{array}{c}
\rho u_1 \text{erf}(s_1) + \rho \frac{e^{-s_1^2}}{\sqrt{\pi \beta}}  \\
(p+\rho u_1^2) \text{erf}(s_1) + \rho u_1 \frac{e^{-s_1^2}}{\sqrt{\pi \beta}}  \\
\rho u_1 u_2 \text{erf}(s_1) + \rho u_2 \frac{e^{-s_1^2}}{\sqrt{\pi \beta}}  \\
\left( \frac{\gamma}{\gamma-1} p + \frac{1}{2} \rho (u_1^2+u_2^2)\right) u_1\text{erf}(s_1) + \left( \frac{\gamma+1}{2(\gamma-1)} p + \frac{1}{2} \rho (u_1^2+u_2^2)\right)\frac{e^{-s_1^2}}{\sqrt{\pi \beta}} \\
\end{array} \right  \}
\end{align}
\begin{align}
 <\Psi, \text{sign}(v_2)\,v_2f^M> &= \left \{ 
\begin{array}{c}
\rho u_2 \text{erf}(s_2) + \rho \frac{e^{-s_2^2}}{\sqrt{\pi \beta}}  \\
\rho u_1 u_2 \text{erf}(s_2) + \rho u_1 \frac{e^{-s_2^2}}{\sqrt{\pi \beta}}  \\
(p+\rho u_2^2) \text{erf}(s_2) + \rho u_2 \frac{e^{-s_2^2}}{\sqrt{\pi \beta}}  \\
\left( \frac{\gamma}{\gamma-1} p + \frac{1}{2} \rho (u_1^2+u_2^2)\right) u_2\text{erf}(s_2) + \left( \frac{\gamma+1}{2(\gamma-1)} p + \frac{1}{2} \rho (u_1^2+u_2^2)\right)\frac{e^{-s_2^2}}{\sqrt{\pi \beta}} \\
\end{array} \right  \}
\end{align}
\begin{align}
& <\Psi, \text{sign}(v_1) \, v_2 f^M > = \nonumber
\\ &  \left \{ 
\begin{array}{c}
\rho\, u_2 \,\text{erf}(s_1)   \\
\rho u_2   \left( \frac{e^{-s_1^2}}{\sqrt{\pi \beta}} + u_1\text{erf}(s_1)   \right)  \\
 \rho \, \text{erf}(s_1) \left(  \frac{1}{2 \beta}  + u_2^2 \right) \\
---------------------------\\
\rho\,I_0  \, u_2 \text{erf}(s_1)    +\frac{\rho}{2} \, \text{erf}(s_1)  \left(  \frac{3u_2}{2\beta} +u_2^3\right)     \\
+\frac{\rho \, u_2}{2} \sqrt{\frac{\beta}{\pi}} \left( \frac{2}{\beta \sqrt{\beta}} \left[ \frac{-s_1e^{-s_1^2}}{2} + \frac{\sqrt{\pi}}{2} \text{erf}(s_1)\right]+\frac{2u_1}{\beta} e^{-s_1^2} +\frac{u_1^2}{\sqrt{\beta}} \text{erf}(s_1) \sqrt{\pi} \right)  \\
\end{array} \right  \}
\end{align}

Similarly,
\begin{align}
& <\Psi, \text{sign}(v_2) \, v_1 f^M > = \nonumber
\\ &  \left \{ 
\begin{array}{c}
\rho\, u_1 \,\text{erf}(s_2)    \\
 \rho \, \text{erf}(s_2)   \left(  \frac{1}{2 \beta}  + u_1^2 \right) \\
\rho u_1   \left( \frac{e^{-s_2^2}}{\sqrt{\pi \beta}} + u_2 \text{erf}(s_2)   \right)  \\
---------------------------\\
\rho\,I_0  \, u_1 \text{erf}(s_2) 
+\frac{\rho}{2} \, \text{erf}(s_2) \left(  \frac{3u_1}{2\beta} +u_1^3\right)  \\
+ \frac{\rho \, u_1}{2}  \sqrt{\frac{\beta}{\pi}} \left( \frac{2}{\beta \sqrt{\beta}} \left[ \frac{-s_2 e^{-s_2^2}}{2} + \frac{\sqrt{\pi}}{2} \text{erf}(s_2)\right]+\frac{2u_2}{\beta} e^{-s_2^2} +\frac{u_2^2}{\sqrt{\beta}} \text{erf}(s_2) \sqrt{\pi} \right)
\end{array} \right  \}
\end{align}

As usual temporal discretization is done by using $\theta$ method with $\theta =0$ for explicit KSUPG scheme. The nonlinear system of equations is linearized using Picard iteration technique and is solved by using bi-conjugate gradient stabilized method.

\section{Shock Capturing Paramter}

In case of multidimensional KSUPG method, diffusion along streamline direction is not sufficient to suppress the oscillations near high gradient region. Hence additional diffusion terms with a shock capturing parameter is required which can sense these high gradient regions and adds additional diffusion. There are many shock capturing parameters available in the literature \cite{TH, TS}. In this work we present a simple gradient based shock capturing parameter as follows. 

\begin{figure} [h!] 
\centering
\includegraphics[scale=0.75]{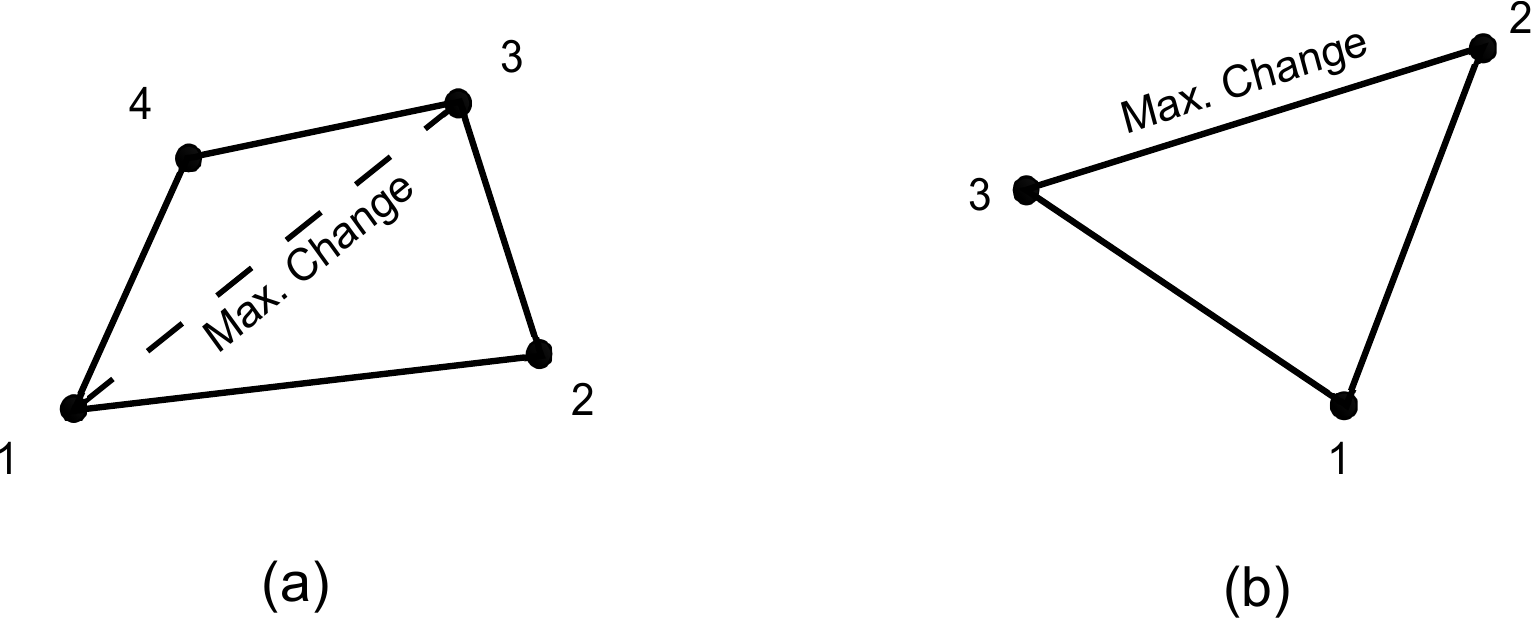}
\caption{Four node quadrilateral element in physical domain}
\label{fig:4nq}
\end{figure}

We define a simple element-wise gradient based shock capturing parameter 
$\delta^{\text{ele}}$ which introduces diffusion along high gradient directions. Figure ~\ref{fig:4nq} (a) shows a typical four node quadrilateral element. As shown in figure, the maximum change in $\Psi$ (where $\Psi$ could be density, temperature or even pressure; in present work, density is used for all numerical test cases because density is the only primitive variable which jumps across all the three waves: shocks, expansion waves and contact disctonintuies) occurs across node 1 and 3. The element based shock capturing parameter is then defined for node 1 and 3 as

\begin{equation}
\delta_{1,3}^{\text{ele}} = \frac{h}{\alpha}\left(\frac{ || \nabla \Psi ||_{\infty}^{\textrm{ele}} }{ ||  \Psi ||_{\infty}^{\textrm{ele}}  }\right)
\end{equation}
where subscripts 1 and 3 represent node numbers.  For nodes 2 and 4, it is defined as
\begin{equation}
\delta_i^{\text{ele}} = \frac{h}{\alpha}\left(\frac{ \Psi_{\text{Max}}^{\textrm{ele}} - \Psi_i^{\text{ele}} }{ ||\Psi ||_{\infty}^{\textrm{ele}}  }\right) \ \  \text{where} \,\, i =2,4.  
\end{equation}
Here $1.4< \alpha \leq 2$.  For most of the test cases $\alpha = 2$ works fine. At element level matrix form, the shock capturing parameter is given by

\begin{equation}
\delta^{\text{ele}} = \left [ 
\begin{array}{cccc}
\delta_{1}^{\text{ele}}  & 0& 0& 0 \\
 0& \delta_{2}^{\text{ele}} & 0  & 0 \\
 0& 0 & \delta_{3}^{\text{ele}} & 0 \\
 0& 0 & 0 &  \delta_{4}^{\text{ele}} \\
\end{array}\right  ] 
\end{equation}
The upper and lower bound on the value $|| \nabla\Psi ||_{\infty}^{\text{ele}}$ is given by  
\begin{equation}
0 \leq || \nabla \Psi ||_{\infty}^{\text{ele}} \leq ||\Psi ||_{\infty}^{\textrm{ele}}
\end{equation}
It is important to note that, the addition of extra shock capturing term in the weak formulation makes the formulation inconsistent with the original equation. Thus, we define $\delta^{\text{ele}}$ such that as $h \rightarrow 0$, $\delta^{\text{ele}}$ should disappear. This condition is achieved by including $h$ in the numerator, which vanishes as we refine the mesh.  Similarly, one can define such a delta parameter for triangular elements shown in figure ~\ref{fig:4nq} (b).  The additional diffusion term along with the shock capturing parameter is then given by 
\begin{equation}\label{SCPn}
\delta (D_x +D_y)
\end{equation} where $\delta$ is the global matrix obtained by assembly and $D_x,D_y$ are the diffusion matrices in $x, y$ direction respectively. These diffusion matrices are defined in 2D Euler KSUPG formulation.

\section{Spectral Stability Analysis}

Stability analysis of a numerical scheme gives the acceptable value of time step $\Delta t$ within which the scheme is stable. In other words, error does not grow with time. Unlike von Neumann stability analysis, spectral stability analysis includes the boundary points too.  In the following analysis, we consider the 2D weak formulation of a linear equation. The global system of equation can be written as
\begin{equation}\label{stab1}
U^{t+\Delta t} = \mathcal{A} U^{t}
\end{equation}
where $\mathcal{A}$ is amplification matrix and $U^{t+\Delta t}$, $U^{t}$ are the numerical solution at time level $t+\Delta t$ and $t$ respectively. Let $U$ be the exact solution, then the error $\epsilon$ is given by
\begin{equation}
\epsilon^{\tau} = U -  U^{\tau}, \ \ \forall \, \tau \in \mathbb{R}^{+}
\end{equation}
Substituting this in equation \eqref{stab1} one can obtain
\begin{equation}
\epsilon^{t+\Delta t} =\mathcal{A} \epsilon^{t}
\end{equation}
\begin{align}
\epsilon^{t+\Delta t} & =\mathcal{A} \epsilon^{t}  = \mathcal{A}^{1+\Delta t} \epsilon^{t-\Delta t} = \cdots  =\mathcal{A}^{t+\Delta t}  \epsilon^{0} 
\end{align}
Rearranging, we get
$$ \frac{\epsilon^{t+\Delta t}}{\epsilon^{0}}  = \mathcal{A}^{t+\Delta t} $$
For stable solution
$$ \left|\left| \frac{\epsilon^{t+\Delta t}}{\epsilon^{0}}  \right|\right| \leq 1 $$
which gives $||\mathcal{A}^{t+\Delta t}  || \leq 1 \rightarrow  ||\mathcal{A} ||\leq 1$.  We use the following relation 
\begin{equation}\label{stks1}
|\varrho(\mathcal{A})| \leq  ||\mathcal{A} || \leq 1
\end{equation}
where $\varrho(\mathcal{A}) $ is the spectral radius of amplification matrix.

Thus, error $\epsilon^{t+\Delta t} $ remains bounded when the maximum eigenvalue of amplification matrix $\mathcal{A}$ is less than or equal to unity.  To find 
$\mathcal{A}$ matrix, we use explicit weak formulation for 2D scalar linear problem which is given by
\begin{align}
 &  M \frac{u^{n+1} -u^{n}}{\Delta t}  +c_1 C_x  u^n  + c_2 C_y u^n + \frac{h}{2}  \sqrt{\frac{\beta}{\pi}} \left\{ D_x\left[ \frac{e^{-s_1^2}}{\pi} + c_1\text{erf}(s_1)  \right]  u^n  \right. \nonumber
\\ &  \left. + D_{xy} (c_2 \text{erf}(s_1)+ c_1 \text{erf}(s_2) )  u^n   +  D_y \left[ \frac{e^{-s_2^2}}{\pi} + c_2\text{erf}(s_2)  \right]  u^n \right\} = 0   \nonumber
\end{align}
here $c_1$ and $c_2$ are constants. Rearranging above equation, we get
\begin{align}
 &  M u^{n+1}  = \left( M  -\Delta t  c_1 C_x +  c_2 C_y+ \frac{h}{2} \sqrt{\frac{\beta}{\pi}} \left\{  D_x  \left[ \frac{e^{-s_1^2}}{\pi} + c_1\text{erf}(s_1)  \right]  \right. \right.\nonumber
\\ &  \left. \left. + D_{xy}   (c_2 \text{erf}(s_1)+ c_1 \text{erf}(s_2) ) +D_y  \left[ \frac{e^{-s_2^2}}{\pi} + c_2\text{erf}(s_2)  \right]  \right\} \right) u^n  \nonumber
\end{align}
which gives matrix $\mathcal{A}$ as
\begin{align}
 & \mathcal{A}  = \left\{ M  -\Delta t \left(  c_1 C_x +  c_2 C_y+ \frac{h}{2} \sqrt{\frac{\beta}{\pi}} \left\{D_x \left[ \frac{e^{-s_1^2}}{\pi} + c_1\text{erf}(s_1)  \right]   \right. \right. \right. \nonumber
\\ &  \left. \left.  \left.+ D_{xy}  (c_2 \text{erf}(s_1)+ c_1 \text{erf}(s_2) ) + D_y \left[ \frac{e^{-s_2^2}}{\pi} + c_2\text{erf}(s_2)  \right] \right\} \right) \right\} / M  \nonumber
\end{align}

Using the stability condition given by equation \eqref{stks1}, we get 
\begin{equation}
\Delta t \leq \frac{1}{ \varrho  \left \{\frac{I}{\Delta t}  - M^{-1}\left(  c_1 C_x +  c_2 C_y+ \frac{h}{2} \sqrt{\frac{\beta}{\pi}} \Upsilon \right) \right\} } 
\end{equation}
where $$ \Upsilon = D_x  \left[ \frac{e^{-s_1^2}}{\pi} + c_1\text{erf}(s_1)  \right]  +   D_{xy} (c_2 \text{erf}(s_1)+ c_1 \text{erf}(s_2) ) + D_y \left[ \frac{e^{-s_2^2}}{\pi} + c_2\text{erf}(s_2)  \right] $$

The maximum eigenvalue $\lambda_{max} $ (which is the absolute maximum of all eigenvalues ) is computed numerically using Rayleigh quotient over a $32 \times 32$ grid for 2D linear convection equation with unity wavespeeds in both directions. The initial condition is a cosine pulse convecting diagonally in a square domain $[0,1] \times [0,1]$.

\begin{figure} [h!] 
\centering
\includegraphics[scale=0.55]{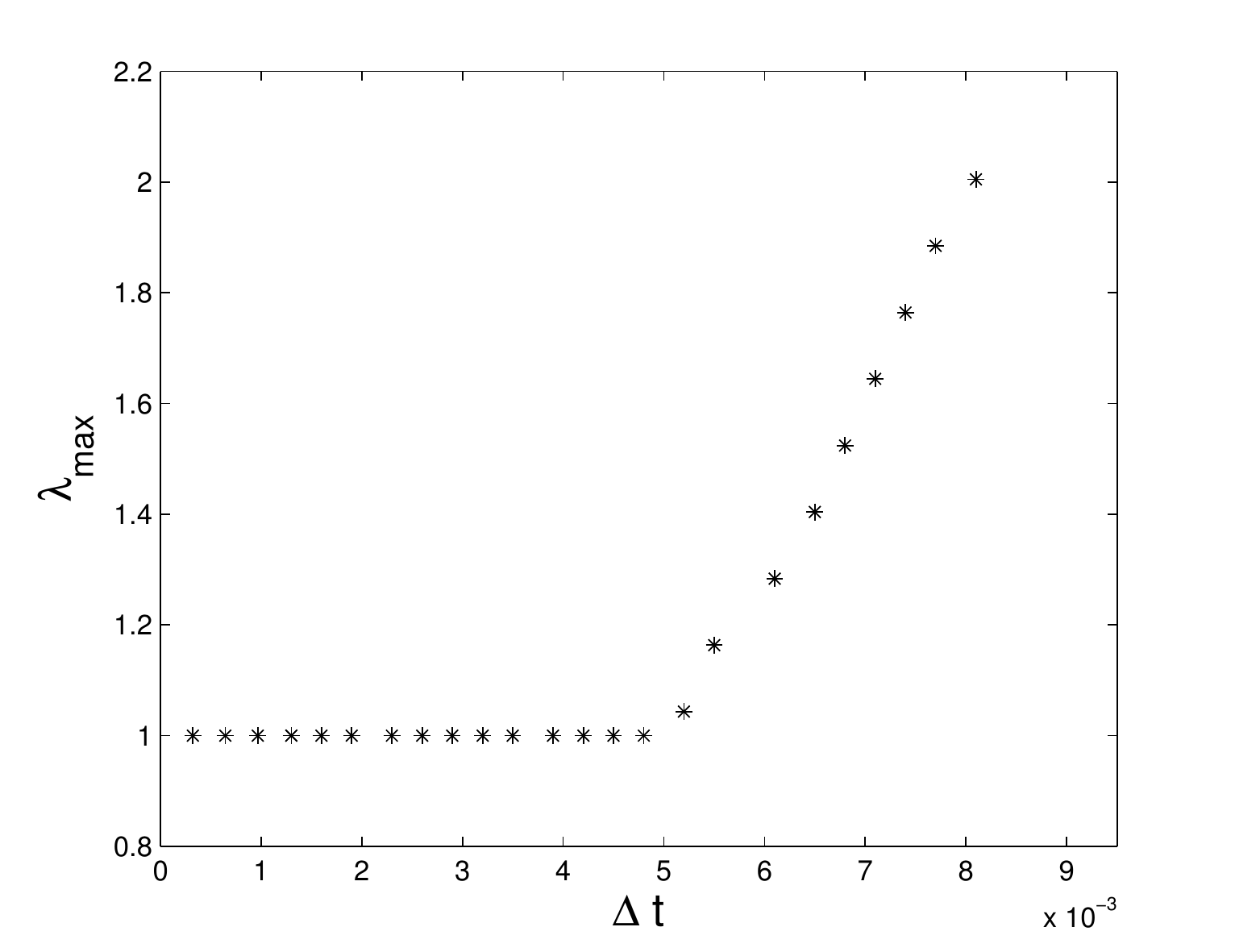}
\caption{Spectral radii of amplification matrix $\mathcal{A}$ at different time steps $\Delta t$.}
\label{fig:stM}
\end{figure} 
Figure ~\ref{fig:stM} shows the spectral radii of amplification matrix $\mathcal{A}$ at different time steps $\Delta t$. It shows the solution become unstable after $\Delta t = 0.005$.

\section{Implicit KSUPG Formulation}

Explicit schemes are slow because of the small values of time-steps used based on the conditional stability limits.  CFL condition puts time step restriction on such schemes.  Therefore, solving problems like vector conservation laws 
({\em e.g.}, Euler equations) using explicit schemes are computationally very expensive.  An alternative is to use implicit schemes which are unconditionally stable without any restriction on time step.  But this advantage comes with extra liability. In this case the coefficient matrix of the solution which needs to be inverted is quite complicated, hence require utmost care.
 
In this section, an implicit formulation of the KSUPG scheme is derived which is computationally very efficient. This formulation is derived for 1D as well as for 
2D Euler equations.

\subsection{One Dimensional Euler Equations}
Using $\theta$ method for temporial discretization of 1D Euler equations
\begin{align}
& M \frac{\mathbf{U}^{n+1}-\mathbf{U}^n }{\Delta t} +(1-\theta) \left(C \mathbf{A} \mathbf{U}^{n+1}  + \frac{h}{2} D< \Psi,\text{sign}(v) \, v f^M >^{n+1}  \right)  \nonumber
\\ & +\theta \left(C \mathbf{A} \mathbf{U}^{n} + \frac{h}{2} D< \Psi,\text{sign}(v) \, v f^M >^{n}  \right)  = 0\nonumber
\end{align}
for $\theta =1$, we get
\begin{align}
M \frac{\mathbf{U}^{n+1}-\mathbf{U}^n }{\Delta t} +\left(C \mathbf{A} \mathbf{U}^{n+1}  + \frac{h}{2} D< \Psi,\text{sign}(v) \, v f^M >^{n+1}  \right)  & = 0\nonumber
\end{align}
where $< \Psi,\text{sign}(v) \, v f^M > = \mathbf{\tilde{D}} \mathbf{U}$
\begin{equation}
\mathbf{\tilde{D}} =   \left [ 
\begin{array}{ccc}
 u \,\textrm{erf}(s) + \frac{e^{-s^2}}{\sqrt{\pi \beta}}  & 0& 0 \\
\frac{p\,\textrm{erf}(s) }{\rho} & u \,\textrm{erf}(s) + \frac{e^{-s^2}}{\sqrt{\pi \beta}}  & 0 \\
 \frac{p\,e^{-s^2} }{2\rho \sqrt{\pi\beta}}  & \frac{p\,\textrm{erf}(s) }{\rho} & u \,\textrm{erf}(s) + \frac{e^{-s^2}}{\sqrt{\pi \beta}}    \\
\end{array} \right ]  
\end{equation}

\subsection{Two Dimensional Euler Equations}
Again, using $\theta$ method for 2D Euler equations
\begin{align}
&  M \frac{\mathbf{U}^{n+1}-\mathbf{U}^n }{\Delta t}  + (1-\theta) \left( C_x \mathbf{G}_x^n  + C_y \mathbf{G}_y^n + \frac{h}{2}D_x  <\Psi,\text{sign}(v_1) \, v_1 f^M >^n \right.  \nonumber
\\& \left. + \frac{h}{2} D_{xy}(<\Psi,\text{sign}(v_1) \, v_2f^M >^n+ <\Psi,\text{sign}(v_2) \, v_1f^M>^n )  + \frac{h}{2}  D_y <\Psi, \text{sign}(v_2) \, v_2 f^M>^n \right )\nonumber
\\ &+ \theta \left( C_x \mathbf{G}_x^{n+1}  + C_y \mathbf{G}_y^{n+1}  + \frac{h}{2}D_x  <\Psi,\text{sign}(v_1) \, v_1 f^M >^{n+1}  \right.  \nonumber
\\& \left. + \frac{h}{2} D_{xy}(<\Psi,\text{sign}(v_1) \, v_2f^M >^{n+1} + <\Psi,\text{sign}(v_2) \, v_1f^M>^{n+1}  )  + \frac{h}{2}  D_y <\Psi, \text{sign}(v_2) \, v_2 f^M>^{n+1}  \right ) =0\nonumber
\end{align}
For implicit method $\theta =1$, which gives
\begin{align}
&  M \frac{\mathbf{U}^{n+1}-\mathbf{U}^n }{\Delta t}  + \left( C_x \mathbf{G}_x^{n+1}  + C_y \mathbf{G}_y^{n+1}  + \frac{h}{2}D_x  <\Psi,\text{sign}(v_1) \, v_1 f^M >^{n+1}  \right.  \nonumber
\\& \left. + \frac{h}{2} D_{xy}(<\Psi,\text{sign}(v_1) \, v_2f^M >^{n+1} + <\Psi,\text{sign}(v_2) \, v_1f^M>^{n+1}  )  + \frac{h}{2}  D_y <\Psi, \text{sign}(v_2) \, v_2 f^M>^{n+1}  \right ) =0\nonumber
\end{align}
Decomposing vectors  $\mathbf{G}_1 = \mathbf{A}_1 \mathbf{U}$, \  $\mathbf{G}_2 = \mathbf{A}_2 \mathbf{U}$, \ $<\Psi,\text{sign}(v_1) \, v_1 f^M> = \mathbf{\tilde{D}_x} \mathbf{U}$, \newline $<\Psi,\text{sign}(v_2) \, v_2 f^M > = \mathbf{\tilde{D}_y} \mathbf{U}$, \ $<\Psi,\text{sign}(v_1) \, v_2 f^M > = \mathbf{\tilde{D}_{xy}} \mathbf{U}$, \   and $<\Psi,\text{sign}(v_2) \, v_1 f^M > = \mathbf{\tilde{D}_{yx}} \mathbf{U}$,  where $\mathbf{A}_1$, $\mathbf{A}_2$ are Jacobian matrixes defined previously and
\begin{equation}
\mathbf{\tilde{D}_x} =   \left [ 
\begin{array}{cccc}
u_1 \,\textrm{erf}(s_1) + \frac{e^{-s_1^2}}{\sqrt{\pi \beta}}  & 0& 0&0  \\
\frac{p\,\textrm{erf}(s_1) }{\rho} & u_1 \,\textrm{erf}(s_1) + \frac{e^{-s_1^2}}{\sqrt{\pi \beta}}  & 0  &0\\
0&0& u_1 \,\textrm{erf}(s_1) + \frac{e^{-s_1^2}}{\sqrt{\pi \beta}}  & 0  \\
0 &0 &0& \Theta_{11}
\end{array} \right ]  
\end{equation}
\begin{equation}
\mathbf{\tilde{D}_y} =   \left [ 
\begin{array}{cccc}
u_2 \,\textrm{erf}(s_2) + \frac{e^{-s_2^2}}{\sqrt{\pi \beta}}  & 0& 0&0  \\
0& u_2 \,\textrm{erf}(s_2) + \frac{e^{-s_2^2}}{\sqrt{\pi \beta}}  & 0  &0\\
\frac{p\,\textrm{erf}(s_2) }{\rho} &0& u_2 \,\textrm{erf}(s_2) + \frac{e^{-s_2^2}}{\sqrt{\pi \beta}}  & 0  \\
0 &0 &0& \Theta_{22}
\end{array} \right ]  
\end{equation}
\begin{equation}
\mathbf{\tilde{D}_{xy}} =   \left [ 
\begin{array}{cccc}
u_2 \,\textrm{erf}(s_1) & 0& 0&0  \\
0& 0& u_1 \,\textrm{erf}(s_1) + \frac{e^{-s_1^2}}{\sqrt{\pi \beta}}   &0\\
\textrm{erf}(s_1) \left( \frac{1}{2\beta} +u_2^2\right) &0& 0  & 0  \\
0 &0 &0& \Theta_{12}
\end{array} \right ]  
\end{equation}
\begin{equation}
\mathbf{\tilde{D}_{yx}} =   \left [ 
\begin{array}{cccc}
u_1 \,\textrm{erf}(s_2) & 0& 0&0  \\
\textrm{erf}(s_2) \left( \frac{1}{2\beta} +u_1^2\right) &0& 0  & 0  \\
0& u_2 \,\textrm{erf}(s_2) + \frac{e^{-s_2^2}}{\sqrt{\pi \beta}} &0  &0\\
0 &0 &0& \Theta_{21}
\end{array} \right ]  
\end{equation}
with 
\begin{align}
\Theta_{11} = & \frac{\left\{ \left( \frac{\gamma}{\gamma-1} p + \frac{1}{2} \rho (u_1^2+u_2^2)\right) u_1\text{erf}(s_1) + \left( \frac{\gamma+1}{2(\gamma-1)} p + \frac{1}{2} \rho (u_1^2+u_2^2)\right)\frac{e^{-s_1^2}}{\sqrt{\pi \beta}} \right\}}{\rho E}\nonumber
\end{align}
\begin{align}
\Theta_{22} = &  \frac{\left\{ \left( \frac{\gamma}{\gamma-1} p + \frac{1}{2} \rho (u_1^2+u_2^2)\right) u_2\text{erf}(s_2) + \left( \frac{\gamma+1}{2(\gamma-1)} p + \frac{1}{2} \rho (u_1^2+u_2^2)\right)\frac{e^{-s_2^2}}{\sqrt{\pi \beta}} \right\}}{\rho E} \nonumber
\end{align}
\begin{align}
 \Theta_{12} = & \left\{ \rho\,I_0  \, u_2 \text{erf}(s_1) +   +\frac{\rho}{2} \, \text{erf}(s_1)  \left(  \frac{3u_2}{2\beta} +u_2^3\right)   \right. \nonumber   
\\ & \left.+\frac{\rho \, u_2}{2} \sqrt{\frac{\beta}{\pi}} \left(\frac{2u_1}{\beta} + \frac{2}{\beta \sqrt{\beta}} \left[ \frac{-s_1e^{-s_1^2}}{2} + \frac{\sqrt{\pi}}{2} \text{erf}(s_1)\right]+\frac{2u_1}{\beta}  (e^{-s_1^2} - 1)+\frac{u_1^2}{\sqrt{\beta}} \text{erf}(s_1) \sqrt{\pi} \right) \right\} \nonumber
\\ & /(\rho E) \nonumber
\end{align}
\begin{align}
  \Theta_{21} = &  \left\{ \rho\,I_0  \, u_1 \text{erf}(s_2) 
+\frac{\rho}{2} \, \text{erf}(s_2) \left(  \frac{3u_1}{2\beta} +u_1^3\right)   \right.\nonumber
\\  & + \left. \frac{\rho \, u_1}{2}  \sqrt{\frac{\beta}{\pi}} \left(\frac{2u_2}{\beta} + \frac{2}{\beta \sqrt{\beta}} \left[ \frac{-s_2 e^{-s_2^2}}{2} + \frac{\sqrt{\pi}}{2} \text{erf}(s_2)\right]+\frac{2u_2}{\beta}  (e^{-s_2^2} - 1)+\frac{u_2^2}{\sqrt{\beta}} \text{erf}(s_2) \sqrt{\pi} \right)\right\} \nonumber
\\ & /(\rho E) \nonumber
\end{align} 
the final implicit equation is written as
\begin{align}
&  M \frac{\mathbf{U}^{n+1}-\mathbf{U}^n }{\Delta t}  + \left( C_x \mathbf{A}_1  + C_y \mathbf{A}_2  + \frac{h}{2} \left[ \mathbf{\tilde{D}_x} D_x  +   (\mathbf{\tilde{D}_{xy}} + \mathbf{\tilde{D}_{yx}} )D_{xy}+ \mathbf{\tilde{D}_y} D_y \right] \right ) \mathbf{U}^{n+1} = 0\nonumber
\end{align}

\section{Results and Discussion}

In this section various test cases are solved for inviscid Burgers as well as Euler equations to demonstrate the accuracy, efficiency and robustness of the proposed scheme. These codes are run on 3.10 GHz desktop machine.   

\subsection{1-D Inviscid Burgers equation test case}
Consider the inviscid Burgers equation in 1-D  
$$ \frac{\partial u}{\partial t} +  \frac{\partial }{\partial x} \left( \frac{u^2}{2}\right)  = 0$$ with initial conditions representing a square wave given by
\begin{equation}
u(x,0) = \begin{cases} 1 & \text{for} \,\, |x| < 1/3 \\   -1 & \text{for}  \,\, 1/3 < |x| \leq 1\end{cases}
\end{equation}
Note that $c = u/2$ for this equation.  

\begin{figure} [h] 
\centering
\includegraphics[scale=0.47]{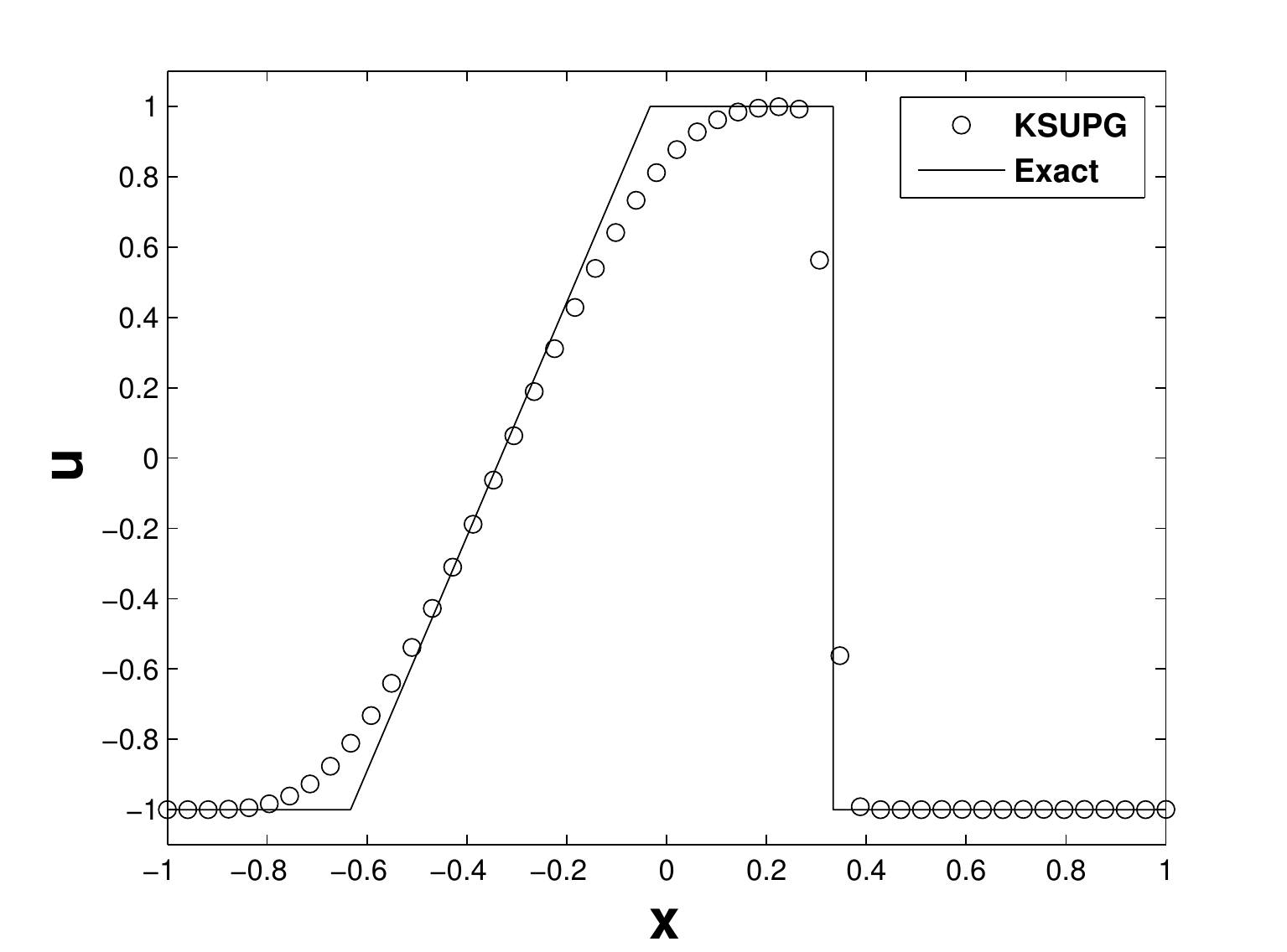}
\caption{Inviscid Burgers Equation Test Case}
\label{fig:ge43}
\end{figure}
In this test case, the jump from -1 to 1 at $x = -1/3$ creates an expansion fan whereas the jump from 1 to -1 creates a steady shock wave.  The solution at time 
$t=0.3$ contains a sonic point (where velocity $u = 0$) in the expansion fan.  Many schemes generate unphysical solutions at sonic points due to insufficient numerical diffusion.  Figure  ~\ref{fig:ge43} shows the numerical solution obtained with 50 grid points using CFL number 0.3. The proposed method does not encounter any sonic point problem.   This scheme captures the shock with just two  grid points.

\subsection{Test cases for 1D Euler equations} 
For 1D Euler equations we shall solve the following test cases.
\subsubsection{Sod's Shock Tube Problem: }
\begin{figure} [h!] 
\centering
\subfigure
{\includegraphics[scale=0.42]{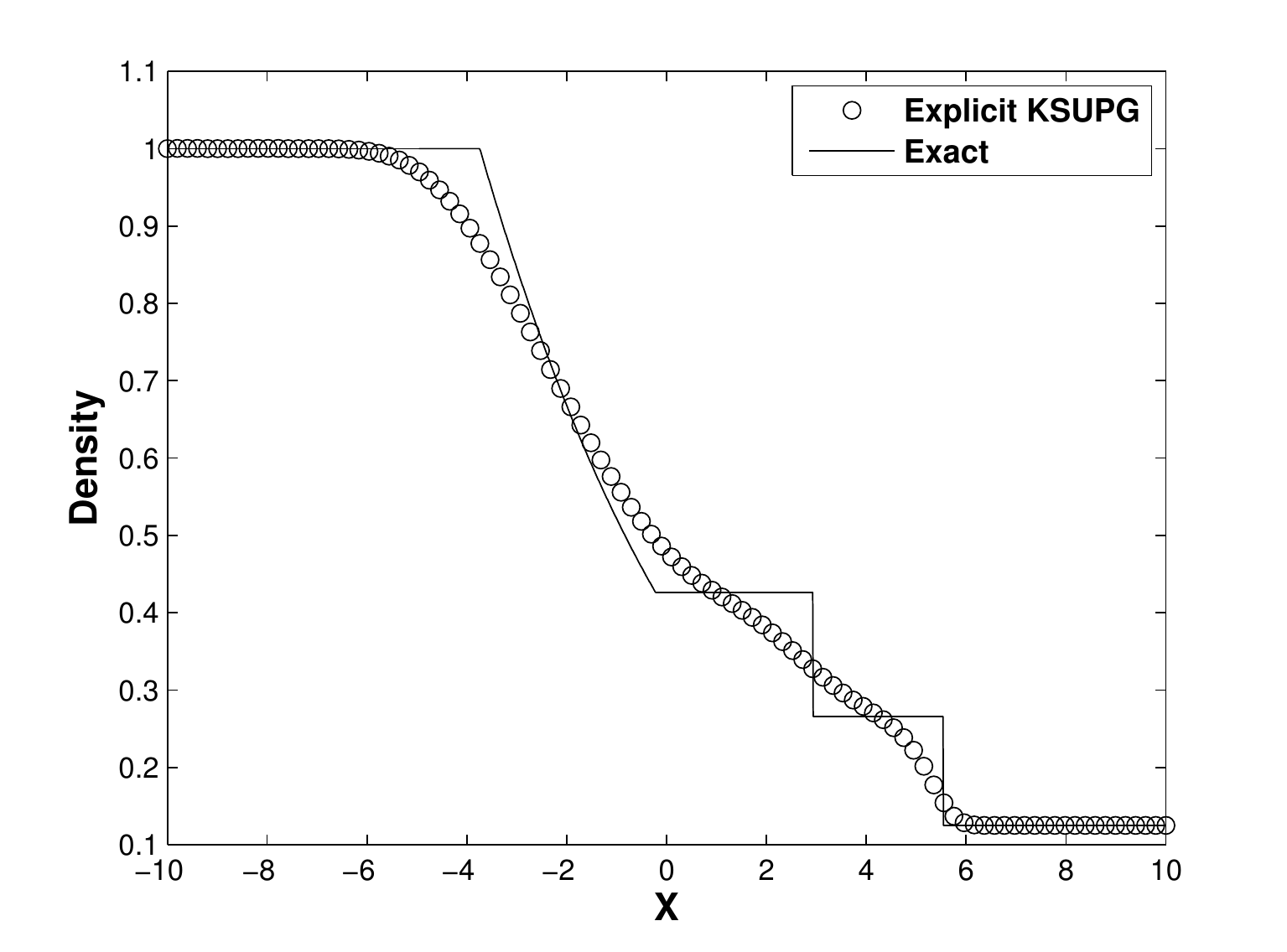}}
\subfigure
{\includegraphics[scale=0.42]{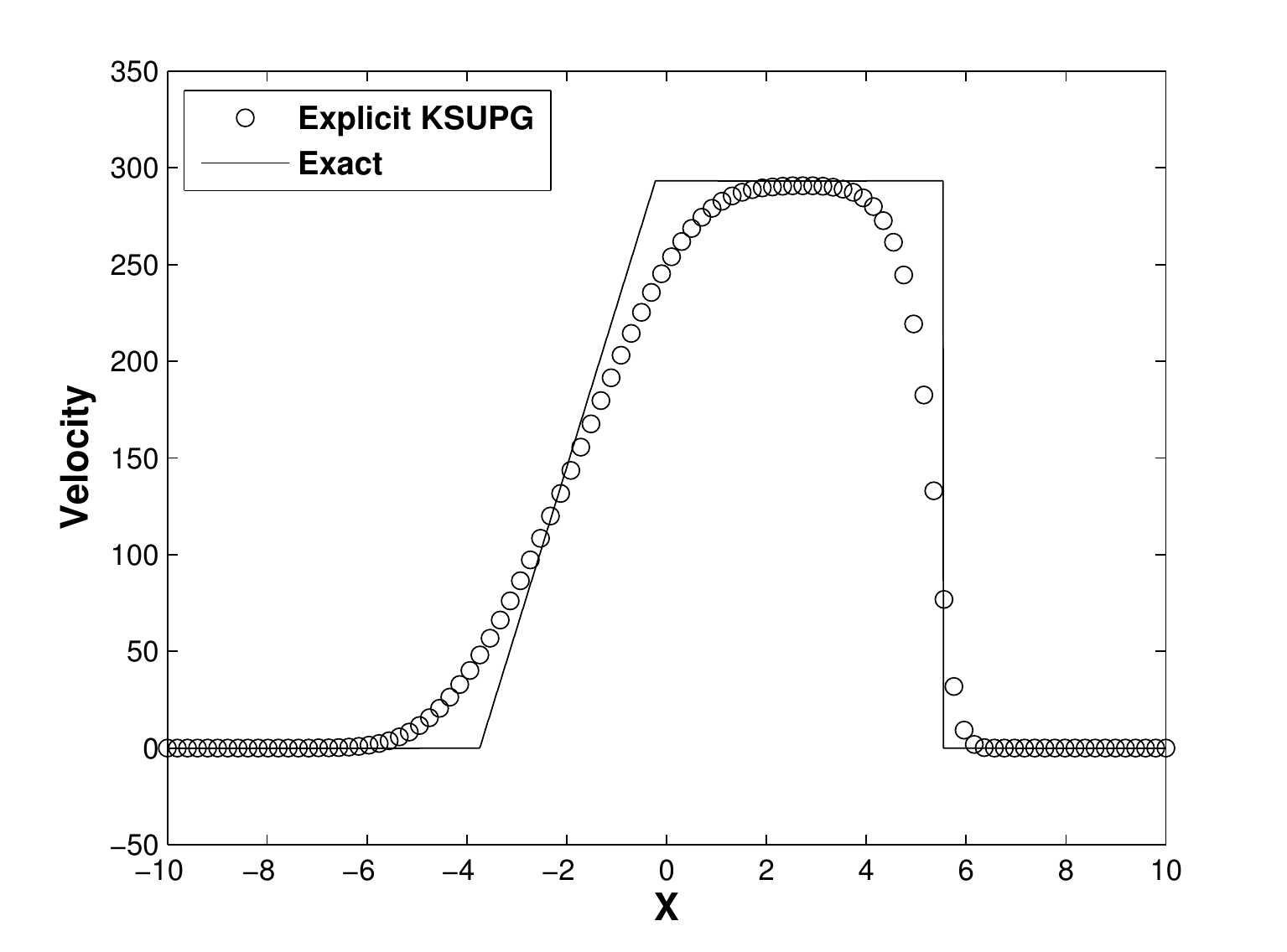}}
\subfigure
{\includegraphics[scale=0.42]{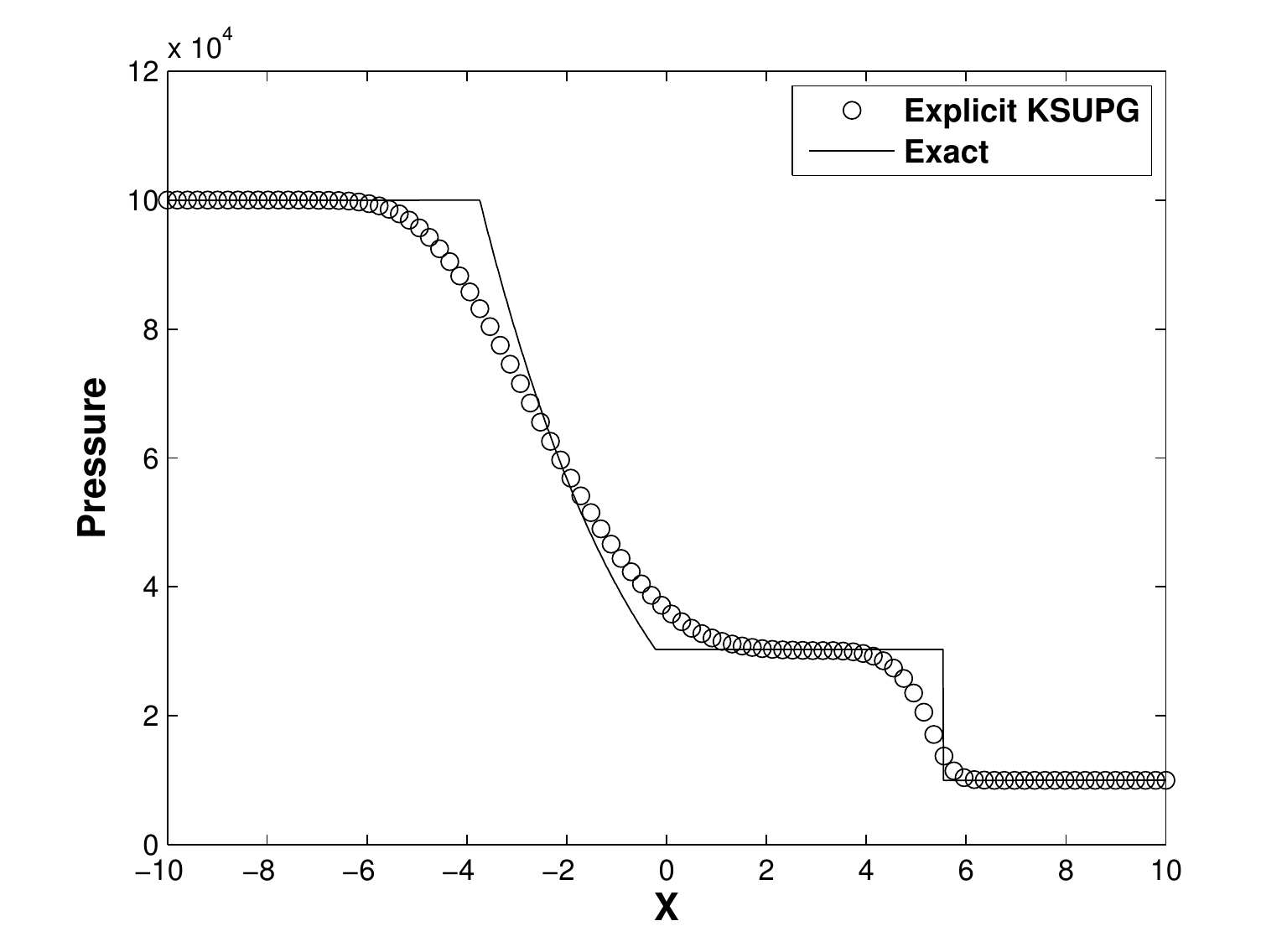}}
\subfigure
{\includegraphics[scale=0.42]{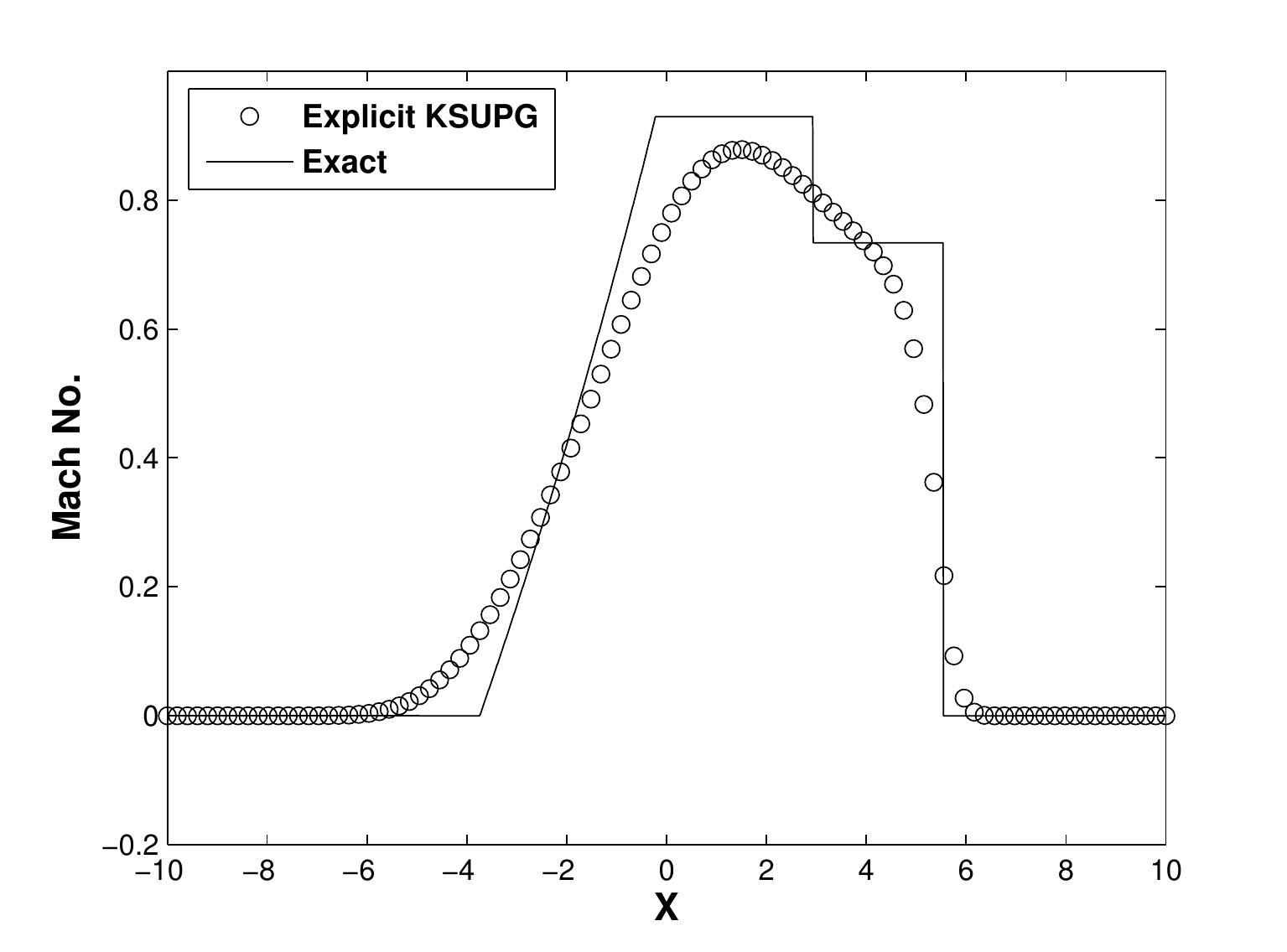}}
\caption{Sod's Shock Tube Problem}
\label{fig:kgE1}
\end{figure}

Sod's shock tube problem consists of a left rarefaction, a right shock wave and a contact discontinuity which separates the rarefaction and shock wave. The initial conditions are given by 
\begin{equation}
\rho (x,0), u(x,0), p(x,0) = \begin{cases} 1,0,100000 & \text{If}\,\,\, -10<x<0 \\ 0.125,0,10000 & \text{If}\,\,\, 0<x<10\end{cases}
\end{equation}
The number of node points are 100 and CFL number is 0.15. Final time is t = 0.01. Figure ~\ref{fig:kgE1} shows the density, velocity, pressure and Mach number plots.  Here, all the essential features like expansion wave, contact discontinuity and shock wave are captured reasonably well.  

\subsubsection{Shock Tube Problem of Lax: }
\begin{figure} [h!] 
\centering
\subfigure
{\includegraphics[scale=0.42]{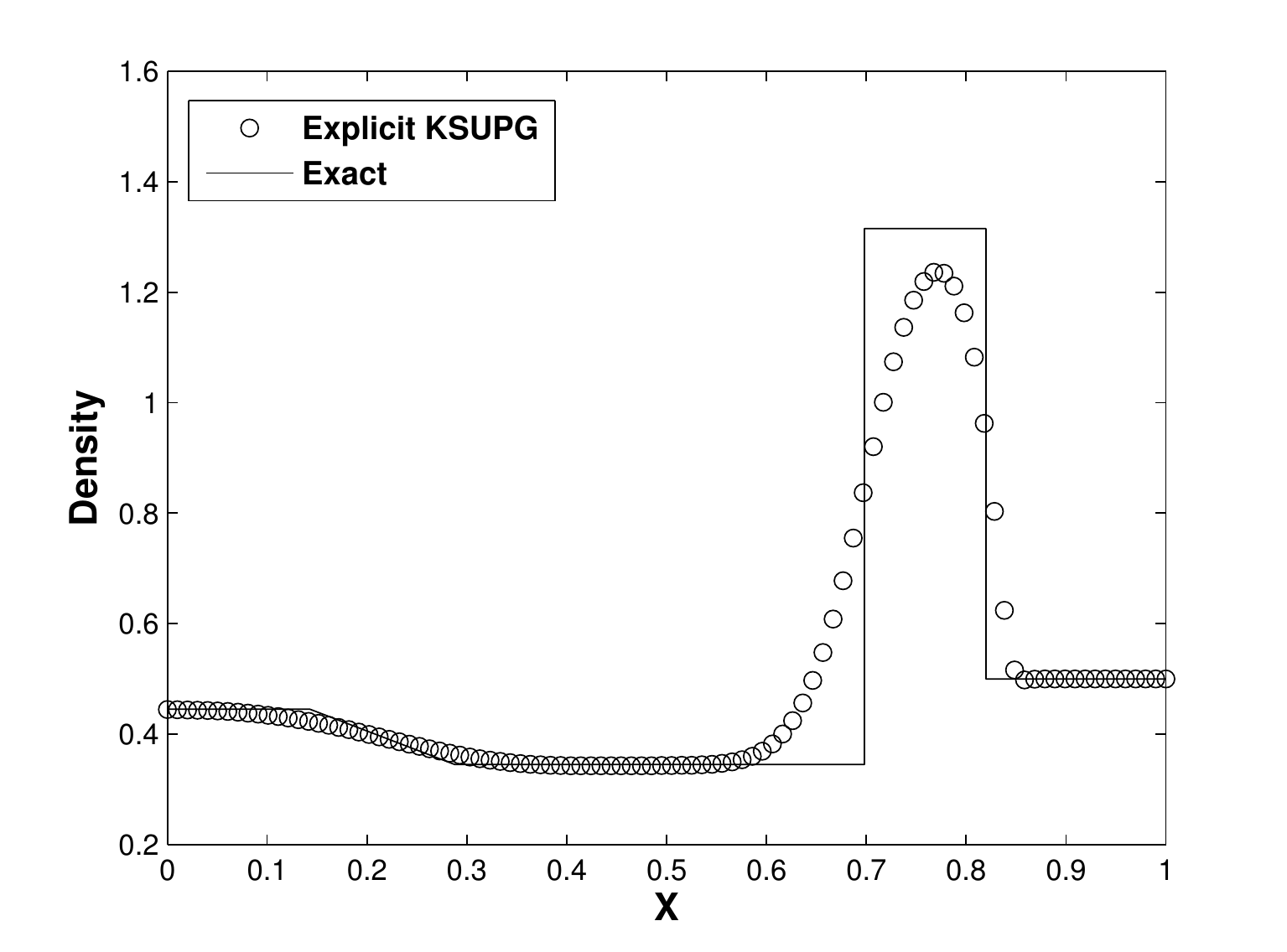}}
\subfigure
{\includegraphics[scale=0.42]{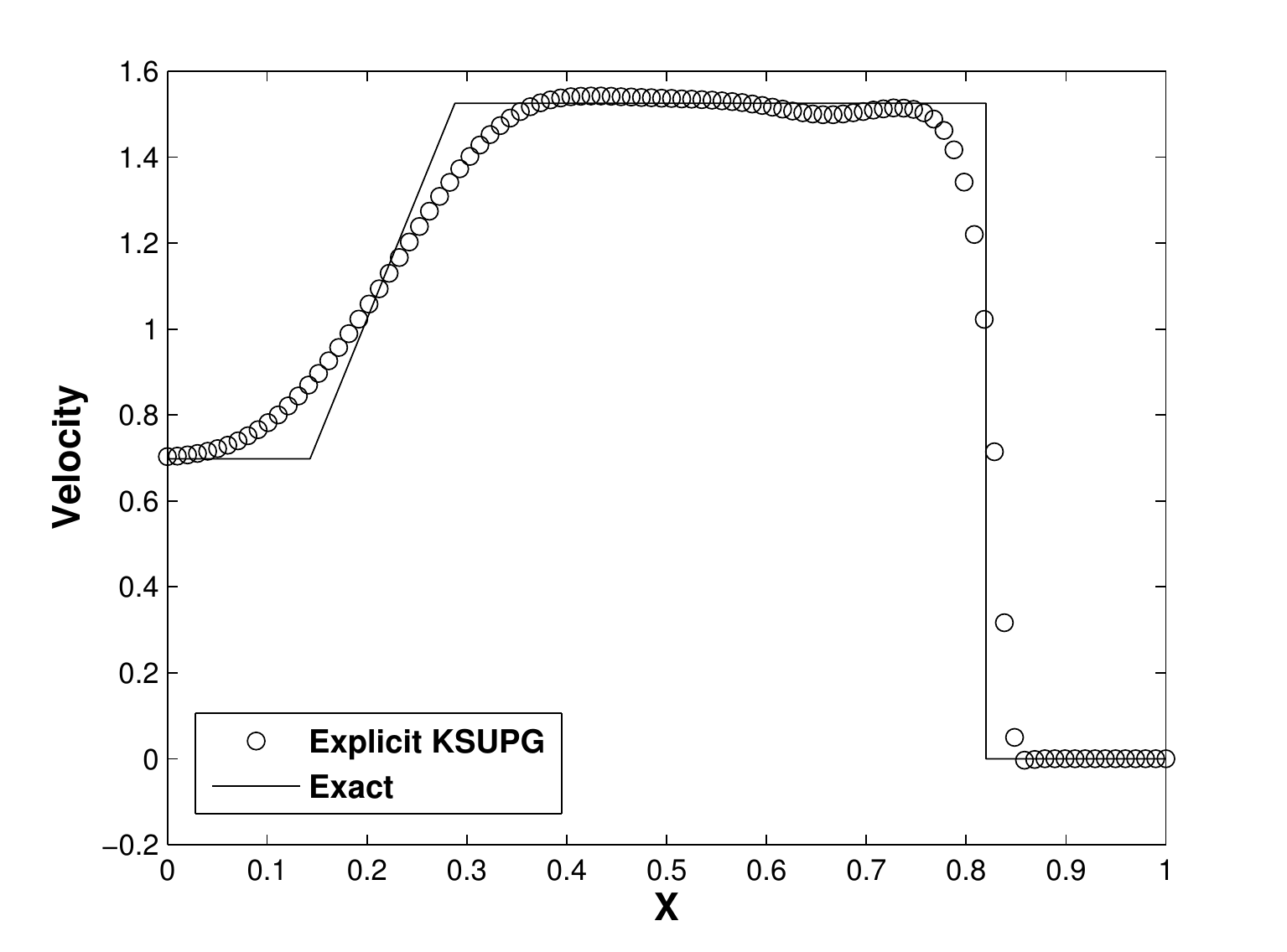}}
\subfigure
{\includegraphics[scale=0.42]{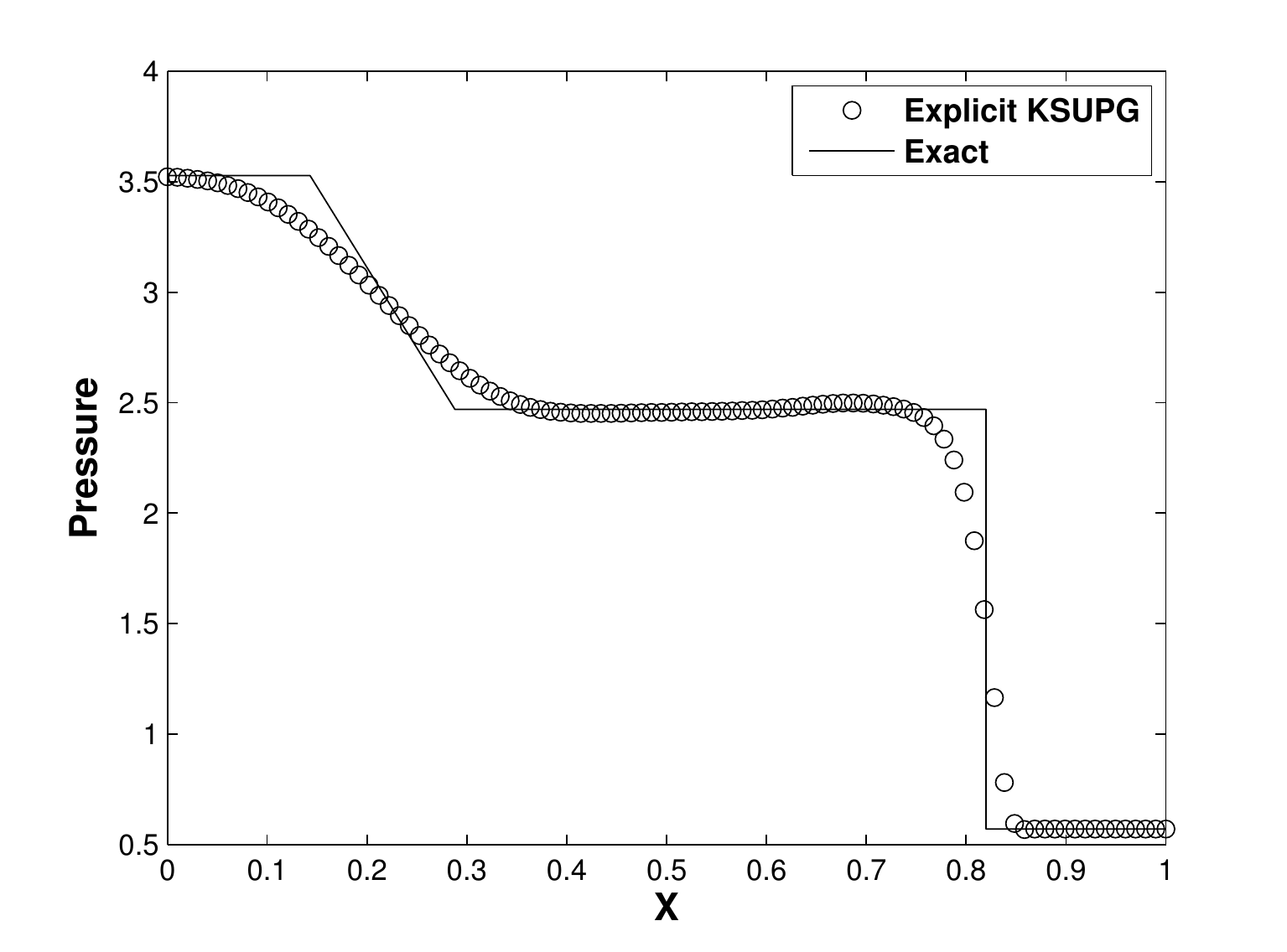}}
\subfigure
{\includegraphics[scale=0.42]{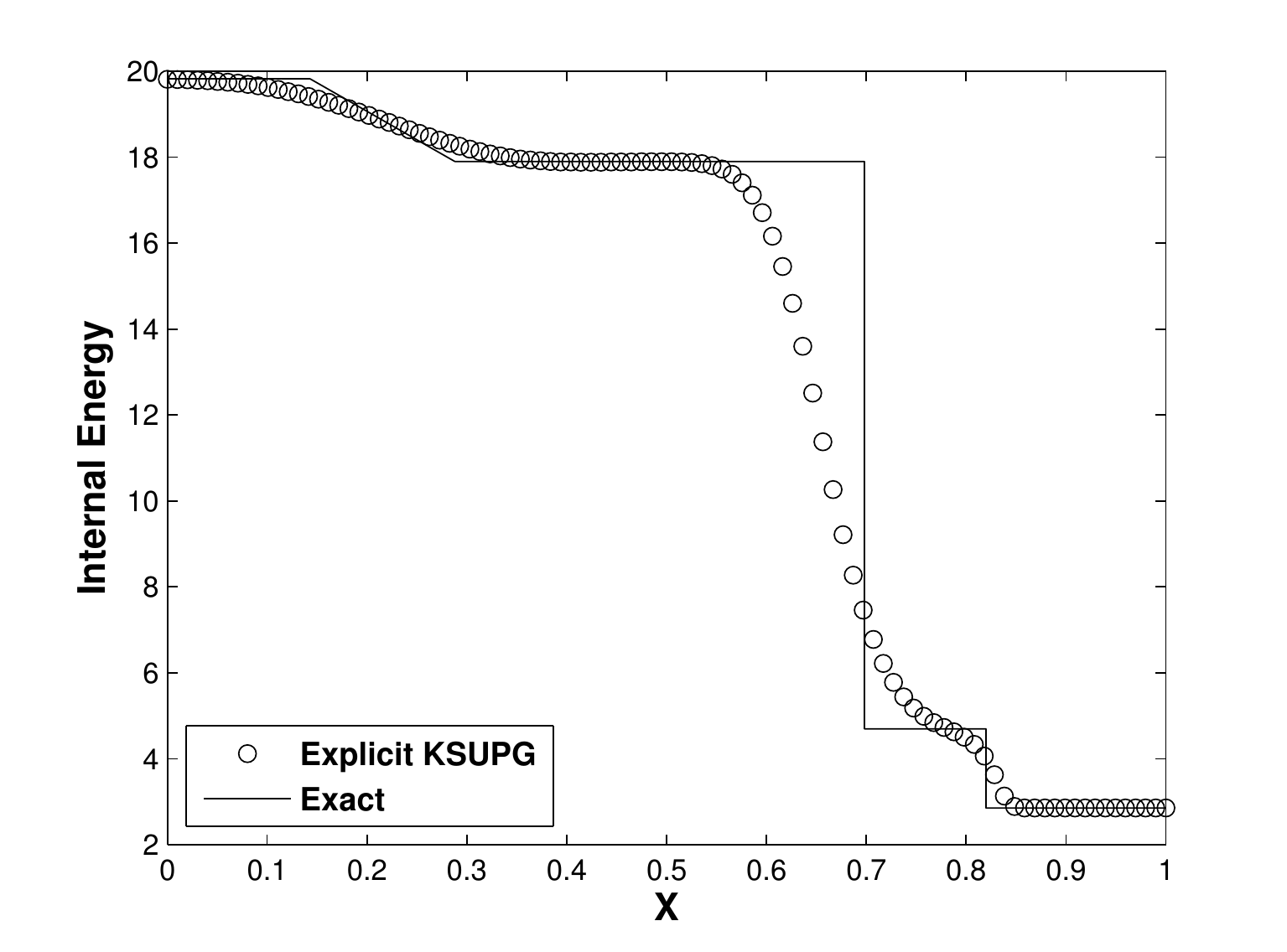}}
\caption{Shock Tube Problem of Lax}
\label{fig:ksla}
\end{figure}

The initial conditions for the shock tube problem of Lax are
\begin{equation}
\rho (x,0), u(x,0), p(x,0) = \begin{cases} 0.445,0.698,3.528 & \text{If}\,\,\, 0<x<0.5 \\ 0.5,0,0.571 & \text{If}\,\,\, 0.5<x<1\end{cases}
\end{equation}
The number of node points are 100 and CFL number is 0.1.  Final time is t=0.13.    Figure ~\ref{fig:ksla} shows the density, velocity, pressure and internal energy plots.

\subsubsection{Strong Rarefactions Riemann Problem}
\begin{figure} [h!] 
\centering
\subfigure
{\includegraphics[scale=0.42]{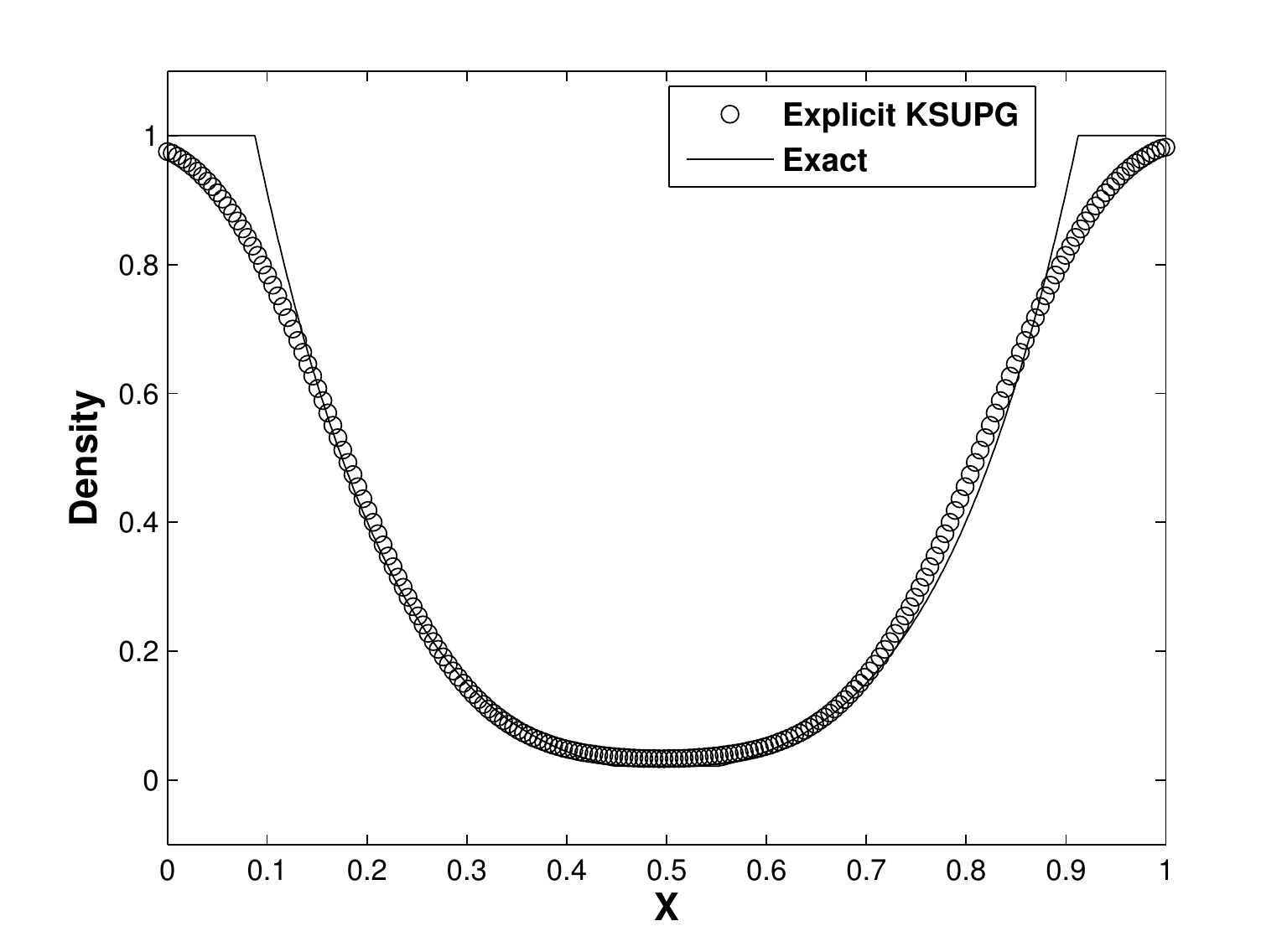}}
\subfigure
{\includegraphics[scale=0.42]{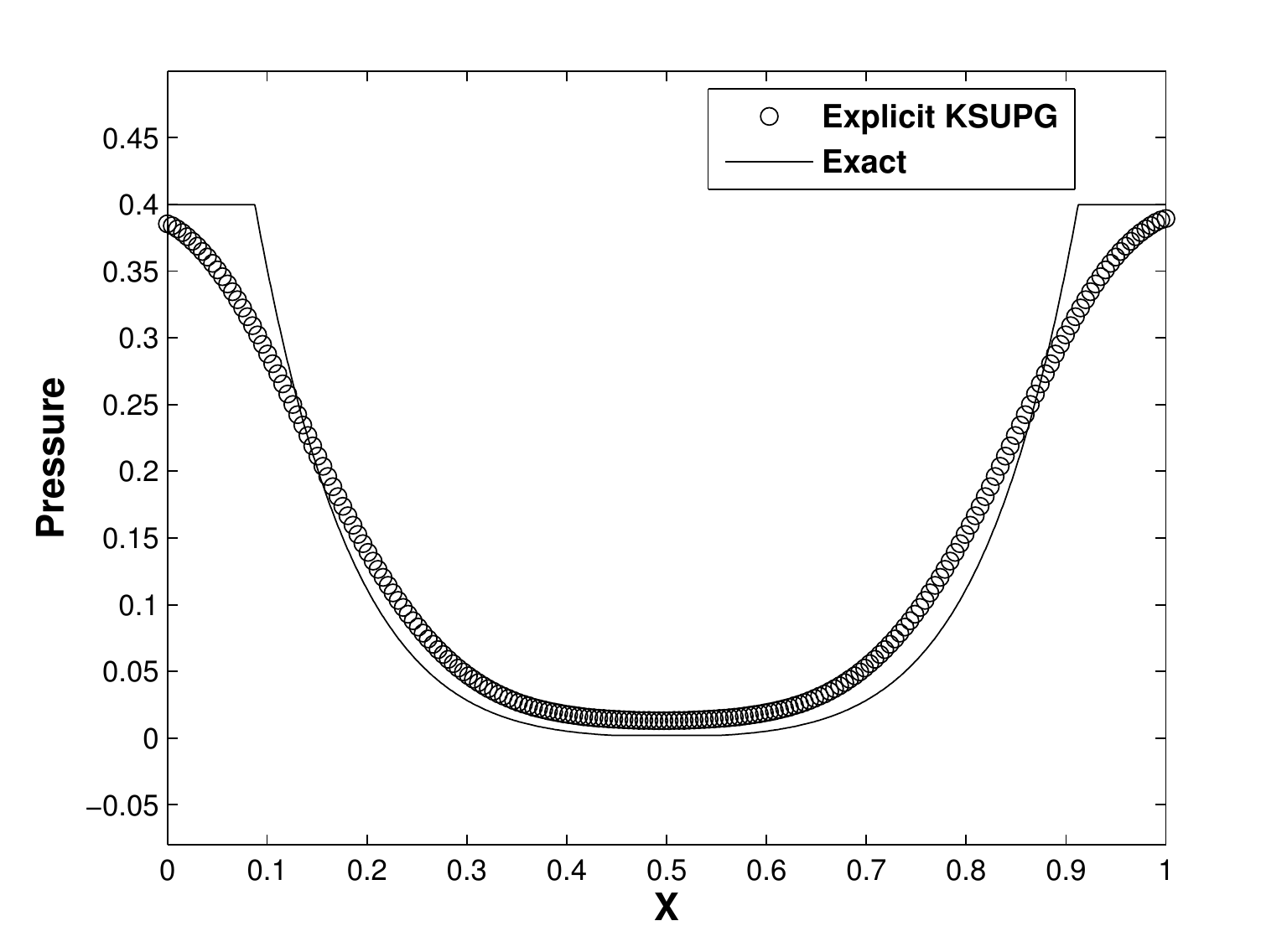}}
\subfigure
{\includegraphics[scale=0.42]{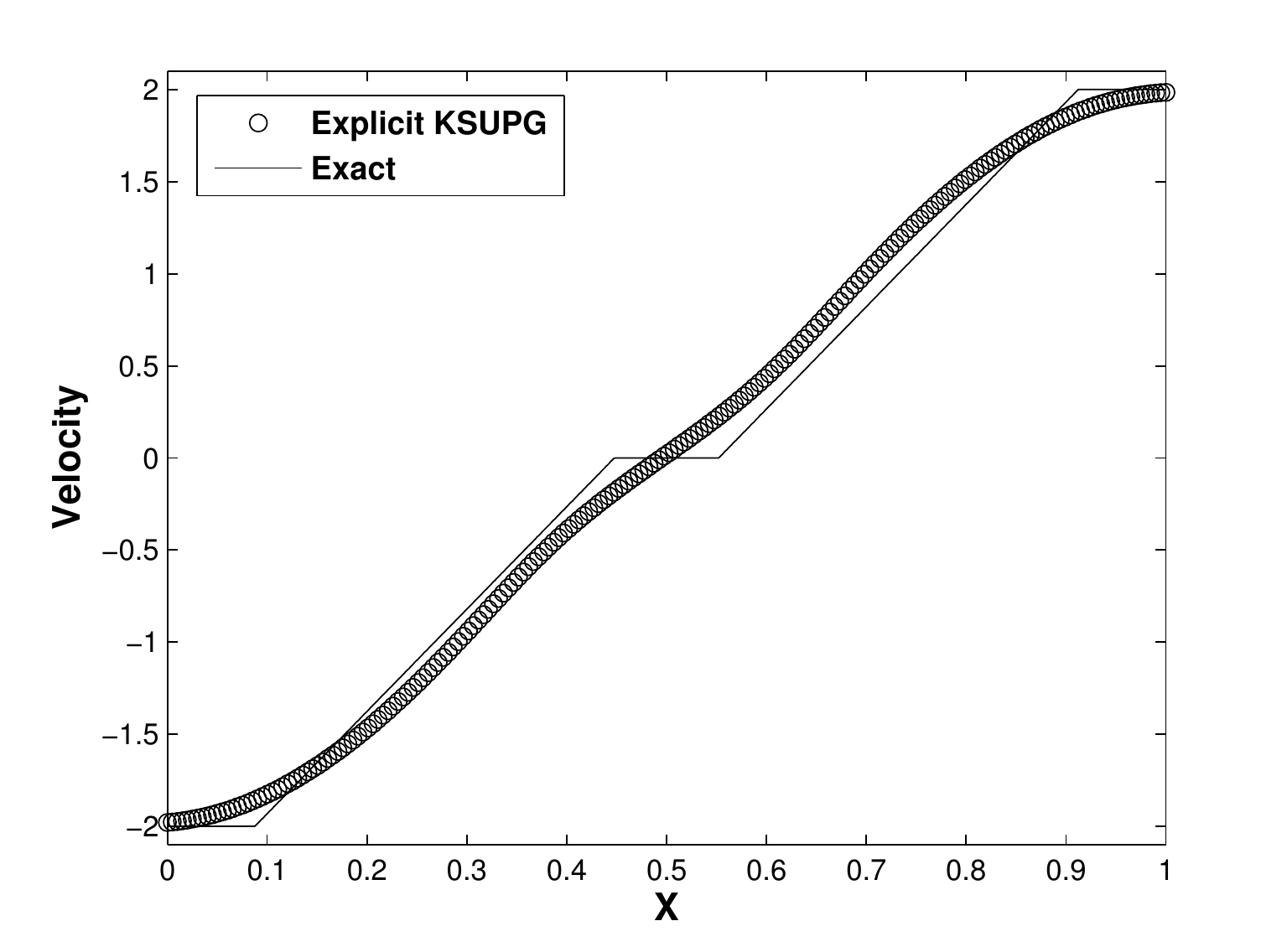}}
\caption{Strong Rarefactions Riemann Problem}
\label{fig:ksrr}
\end{figure}

The initial conditions for the Riemann problem are
\begin{equation}
\rho (x,0), u(x,0), p(x,0) = \begin{cases} 1,-0.2,0.4 & \text{If}\,\,\, 0<x<0.5 \\ 1,2,0.4 & \text{If}\,\,\, 0.5<x<1\end{cases}
\end{equation}
In this test case a near vacuum state is reached. Many popular schemes like linearized Riemann solver fails to predict the correct pressure and density and typically give negative pressure and density.  The number of node points are 200 and CFL number is 0.1. Final time is t=0.15.  Figure ~\ref{fig:ksrr} shows the density, pressure and velocity plots. 

\subsection{Two dimensional Burgers equation}
Two dimensional Burgers equations is given by 
\begin{equation}
\frac{\partial u}{\partial t} + \frac{\partial \left(\frac{u^2}{2}\right)}{\partial x} + \frac{\partial u}{\partial y} = 0
\end{equation}
The boundary conditions are:
$$u(0,y) = 1 \,\,\text{and}\,\, u(1,y) = -1, \ \ 0<y<1$$  and $$u(x,0) = 1-2x, \ \ 0<x<1$$
Exact solution is given in \cite{Spekreijse}.  

\begin{figure} [h!] 
\centering
\includegraphics[scale=0.43]{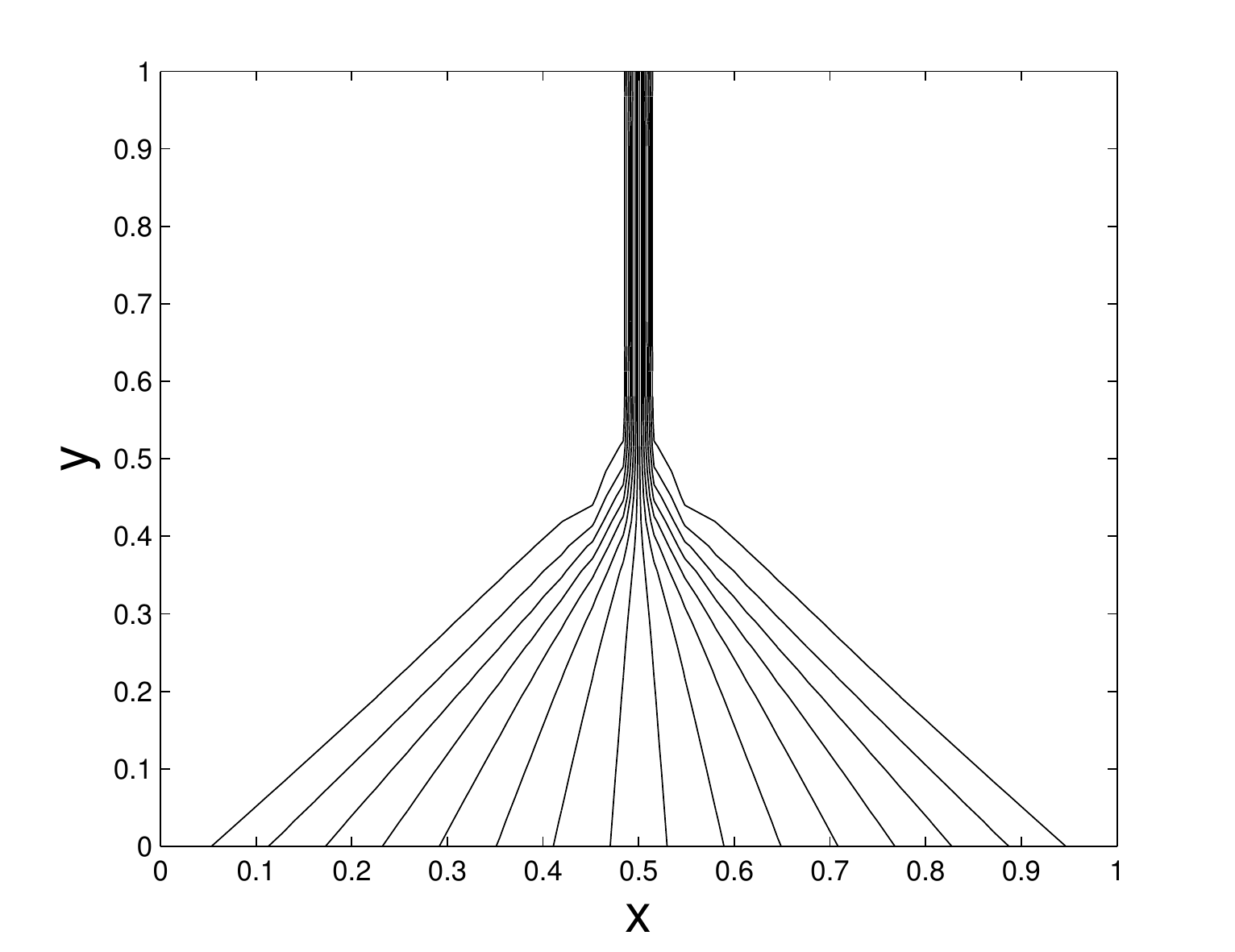}
\includegraphics[scale=0.43]{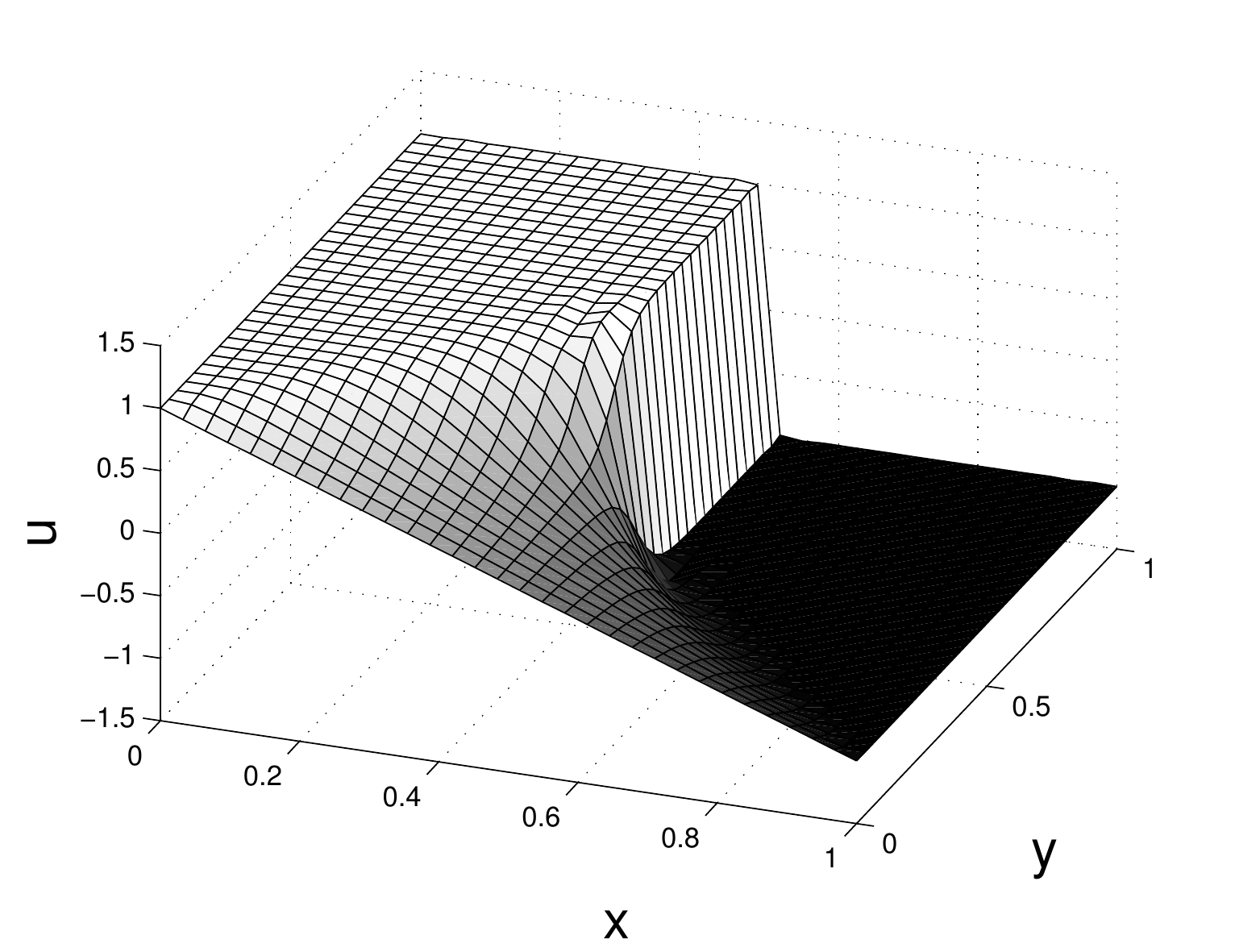}
\caption{Solution of 2D Burgers equation on Q4 mesh}
\label{fig:k2p3aQ}
\end{figure}

\begin{figure} [h!] 
\centering
\includegraphics[scale=0.43]{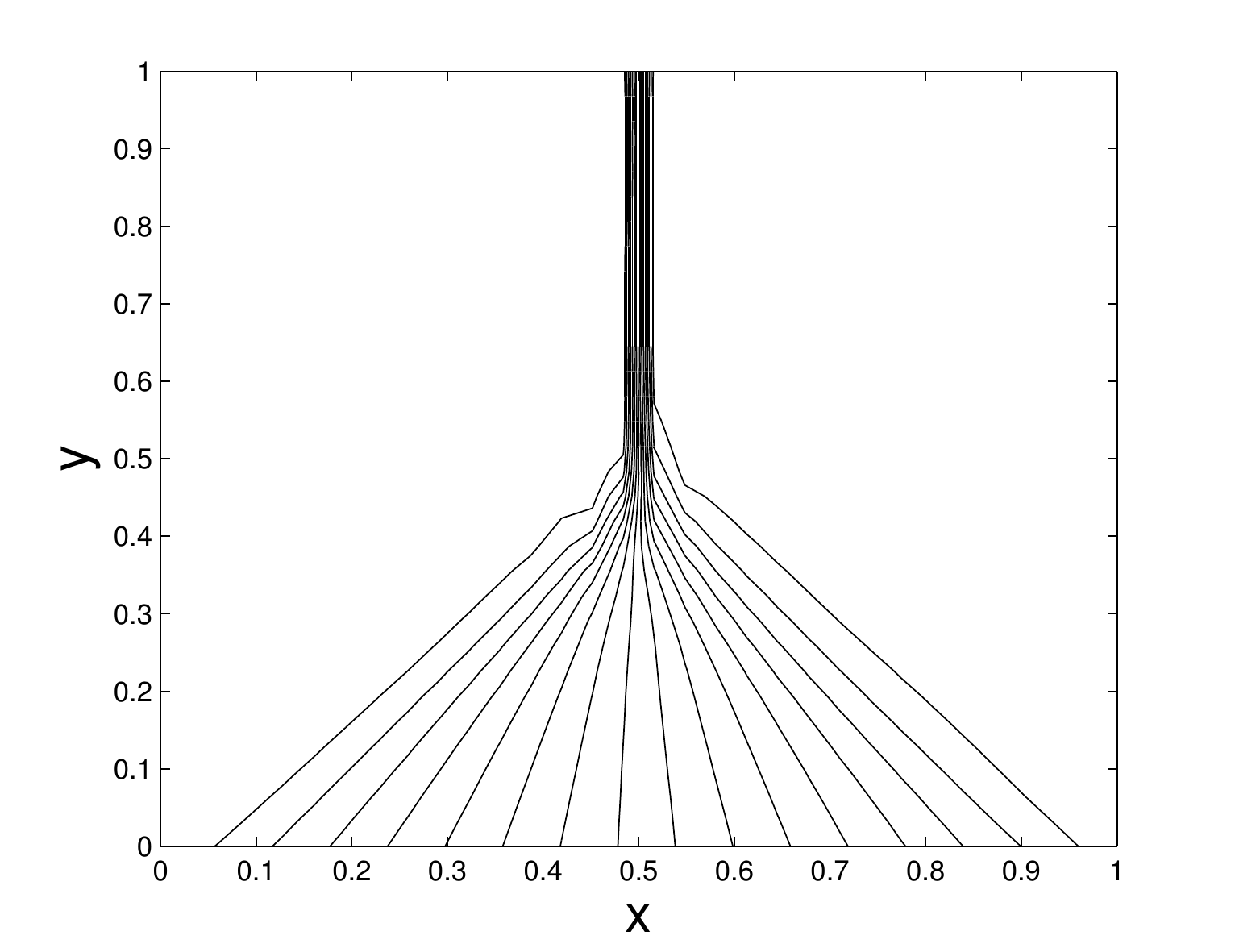}
\includegraphics[scale=0.43]{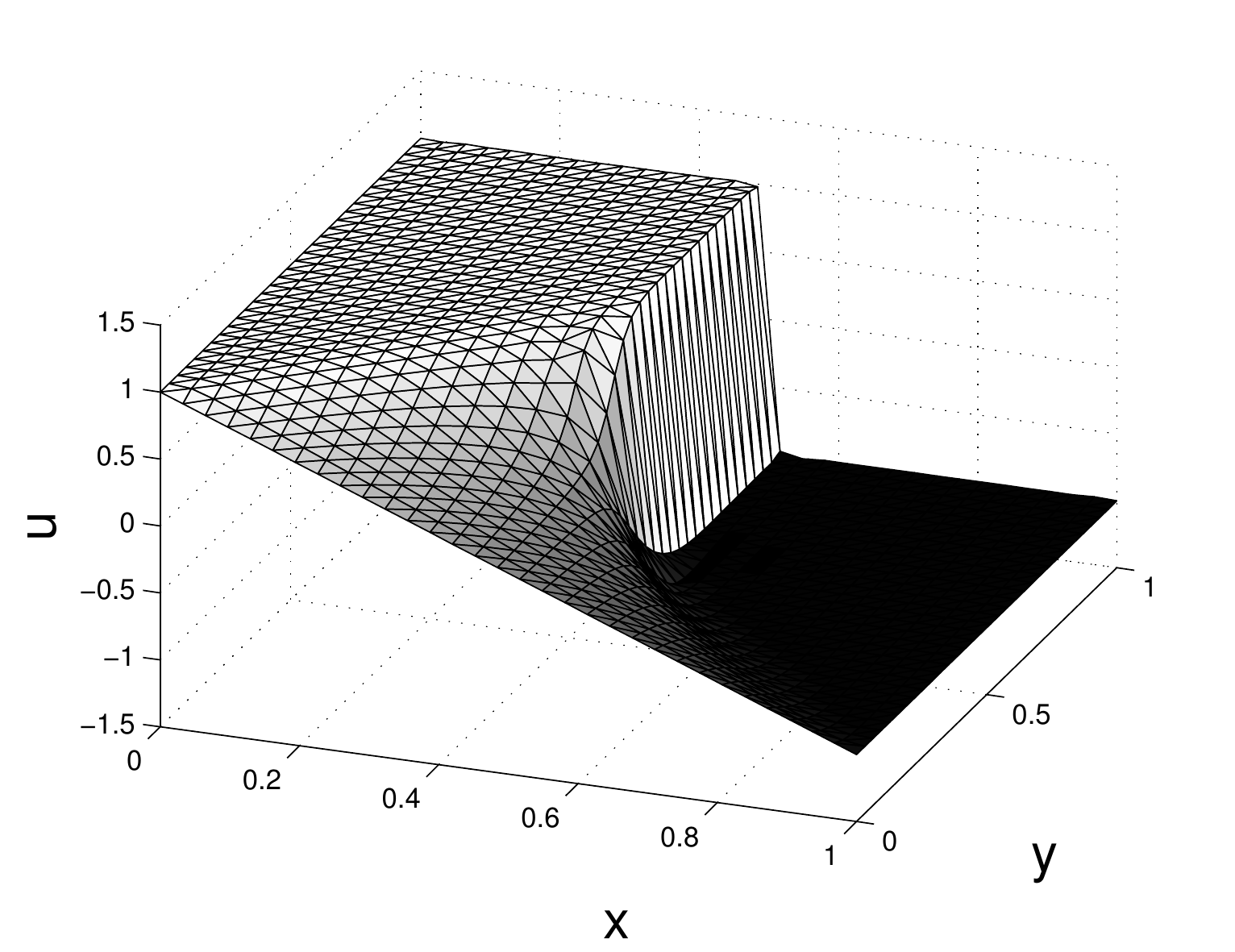}
\caption{Solution of 2D Burgers equation on T3 mesh}
\label{fig:k2p3aT}
\end{figure}

\begin{figure} [h!] 
\centering
\includegraphics[scale=0.5]{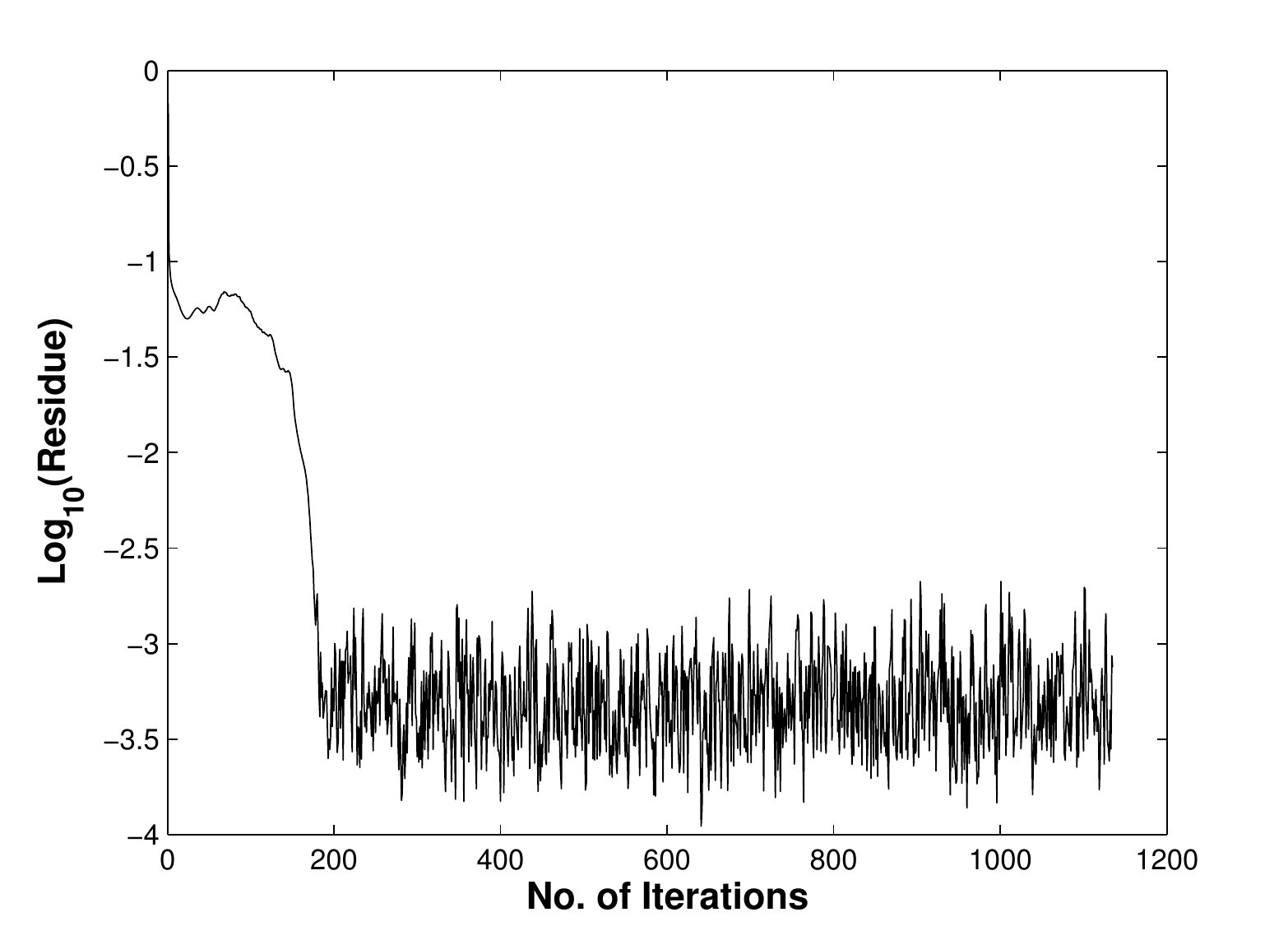}
\caption{Residue plot for 2D Burgers equation on Q4 mesh}
\label{fig:resolp}
\end{figure}

Here, $c_1 = u/2$ and $c_2 = 1$.  Figures ~\ref{fig:k2p3aQ} and ~\ref{fig:k2p3aT}  show the contour and surface plots of steady state solution on $ 32 \times 32$ Q4 and T3 meshes respectively. The normal shock wave is captured quite accurately  using such a coarse grid. Figure ~\ref{fig:resolp} shows the residue plot for Q4 mesh where residue is given by
$$ \text{Residue} = \frac{||u^{n+1}- u^n||_{\mathbb{L}_2}}{||u^{n+1}||_{\mathbb{L}_2}}$$

\begin{figure} [h!] 
\centering
\includegraphics[scale=0.45]{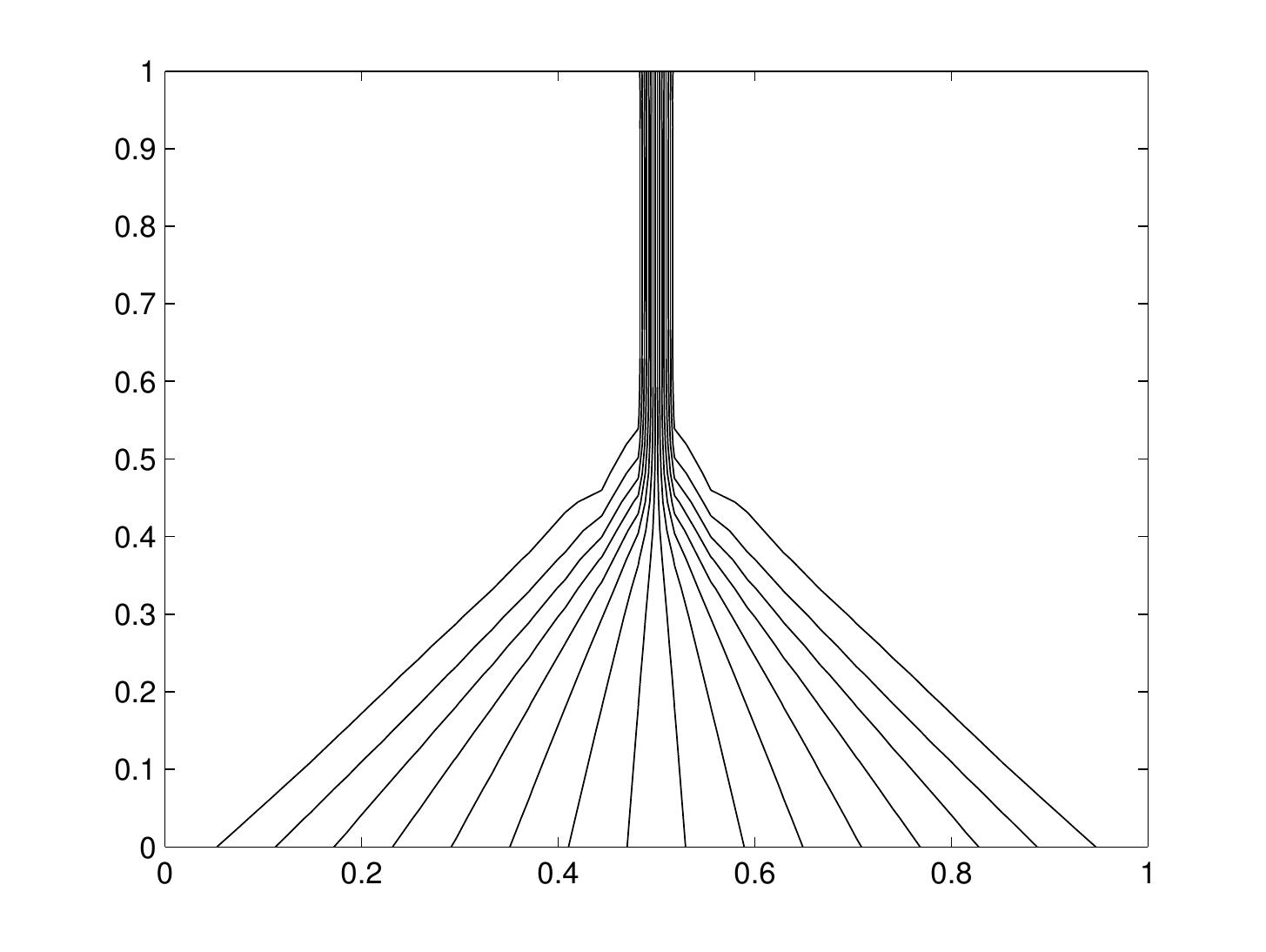}
\includegraphics[scale=0.42]{Spe3aQc.pdf}
\caption{Solution of Burgers equation with original SUPG method using Delta parameter (Left) and the proposed method without any Delta parameter (Right) on $ 32 \times 32 $ Q4 Mesh.}
\label{fig:kqr34}
\end{figure}

\subsubsection{Comparison of Standard SUPG and the Proposed KSUPG Scheme}
The comparison has been done with the standard SUPG method with shock capturing parameter used in \cite{TS} with the proposed KSUPG scheme.  Figure ~\ref{fig:kqr34} shows the result of 2D Burgers equation over $32 \times 32$ Q4 element grid. Accuracy of the solution is more in the proposed scheme even without using shock capturing parameter.

\subsection{2D Euler Test Cases}
Unlike in the case of 2D Burgers equation, additional diffusion is required for simulating 2D Euler equations and thus the shock capturing parameter given by equation \eqref{SCPn} is added in both explicit and implicit KSUPG formulation for 2D Euler equations.  

\subsubsection{Oblique Shock \cite{Hendriana_Bathe}:}
The domain is square $[0,\,1] \times [0,\,1]$, the left boundary is the inlet with Mach 2 at an angle of $-10^o$ to the bottom boundary. Bottom boundary is the wall from where oblique shock wave is generated which makes an angle of $29.3^o$ with the wall. The Dirichlet boundary conditions on left and top boundaries are $ \rho = 1, u_1 = \cos 10^o, u_2 = -\sin 10^o, p = 0.179$. At the wall, no-slip condition is applied, \textit{i.e.}, $\mathbf{v} .n = 0$ where $\mathbf{v}$ is a velocity vector in two dimensions and at right boundary where the flow is supersonic all primitive variables $\rho$, $u $, $v$ and $p$ are extrapolated with first order approximation.
\begin{figure} [h] 
\centering
\includegraphics[scale=0.3]{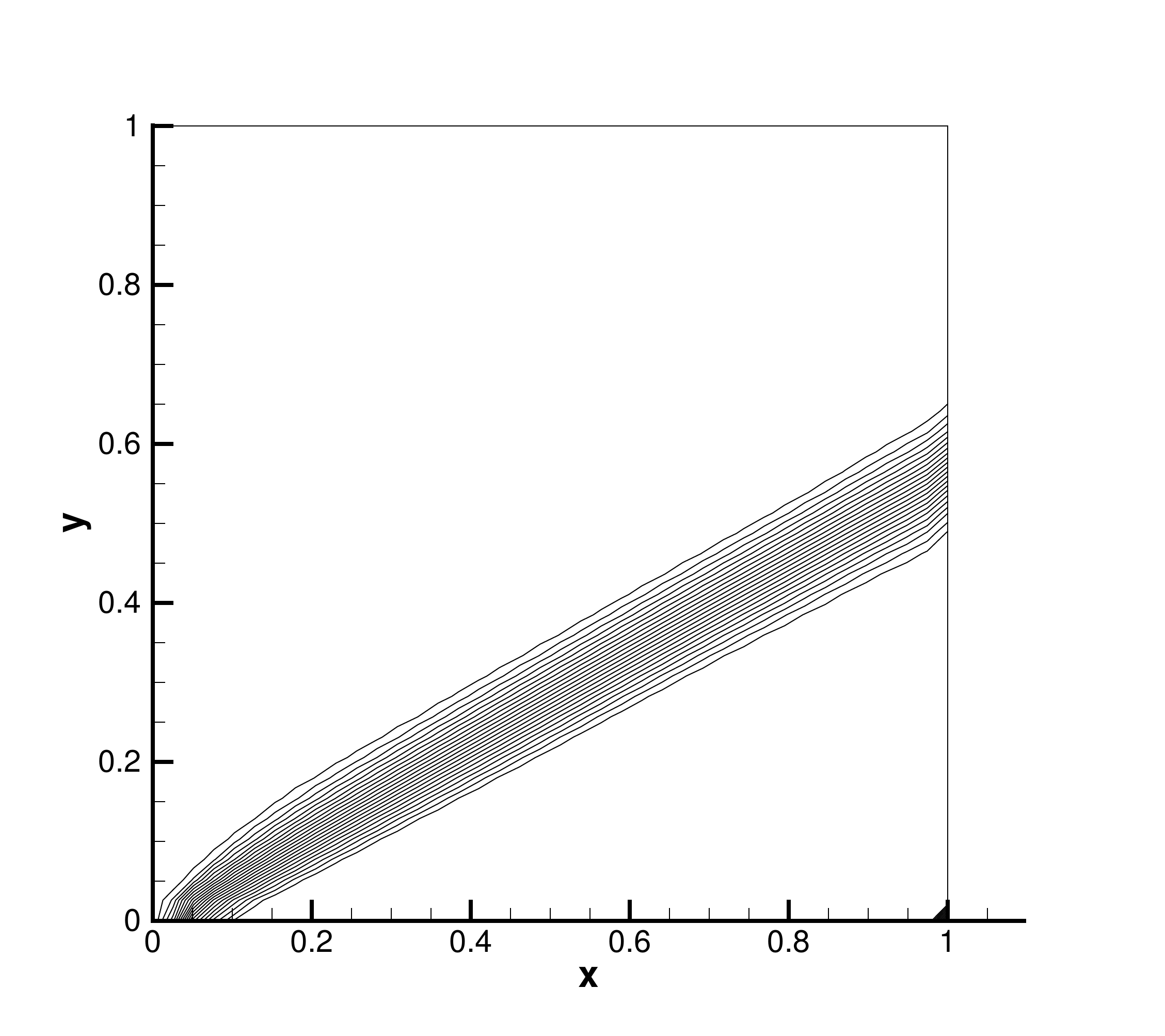}
\caption{Pressure contours (0.18:0.005:0.29) using $40 \times 40$ Q4 elements. }
\label{fig:oSI312}
\end{figure}
 Figure  ~\ref{fig:oSI312} shows the pressure contours using $40 \times 40$ Q4 mesh.

\subsubsection{Oblique Shock Reflection from a flat plate: }
\begin{figure} [h!] 
\centering
\includegraphics[trim=1cm 0.9cm 1cm 12cm, clip=true, scale=0.6, angle = 0]{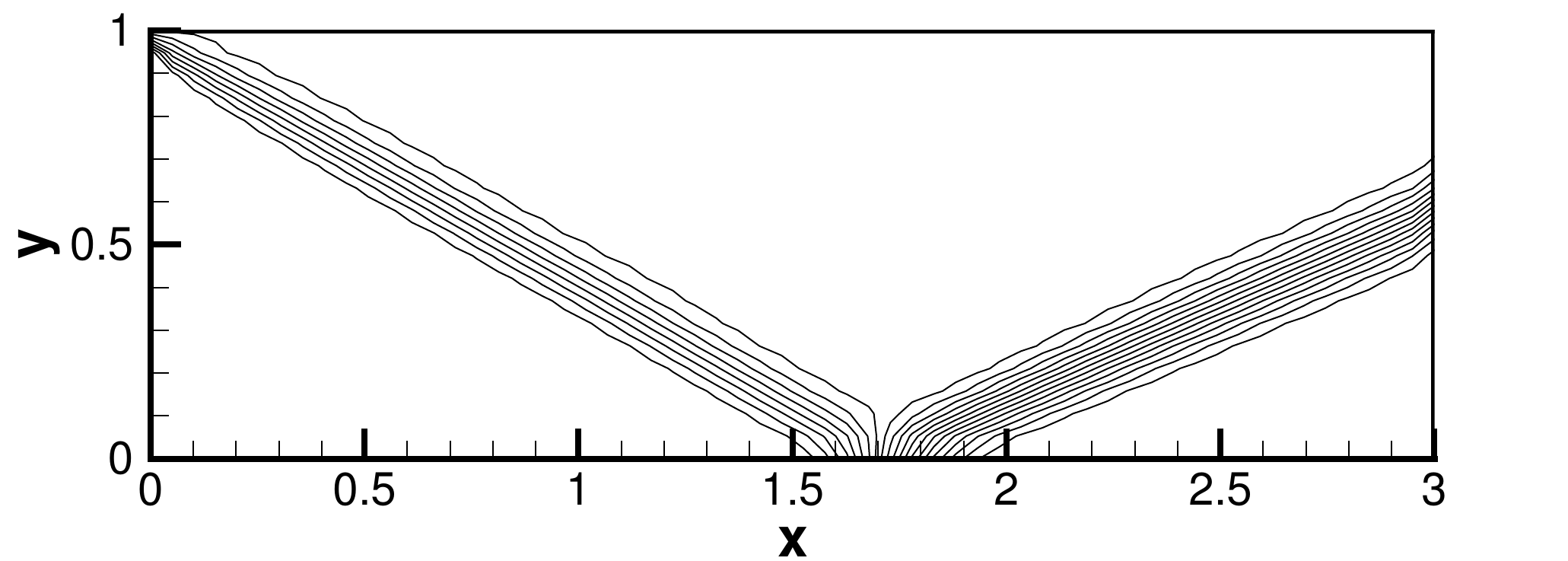}
\includegraphics[trim=1cm 0.9cm 1cm 12cm, clip=true, scale=0.6, angle = 0]{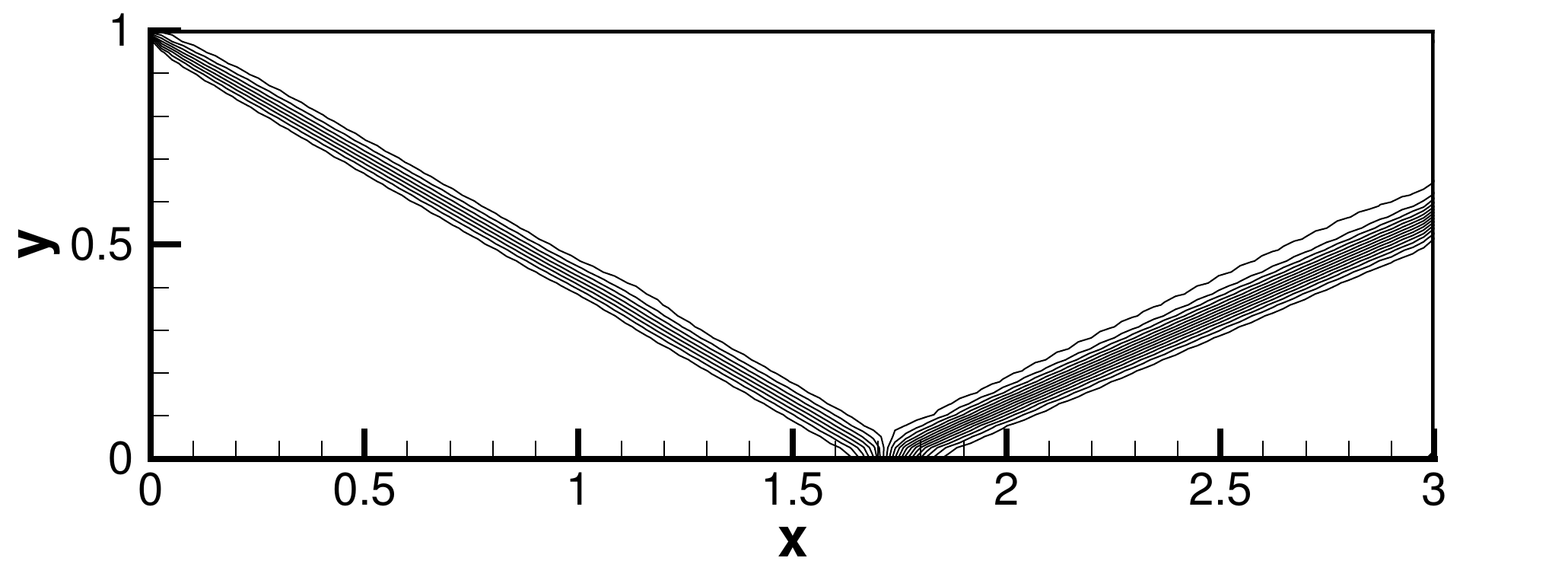}
\includegraphics[trim=1cm 0.9cm 1cm 12cm, clip=true, scale=0.6, angle = 0]{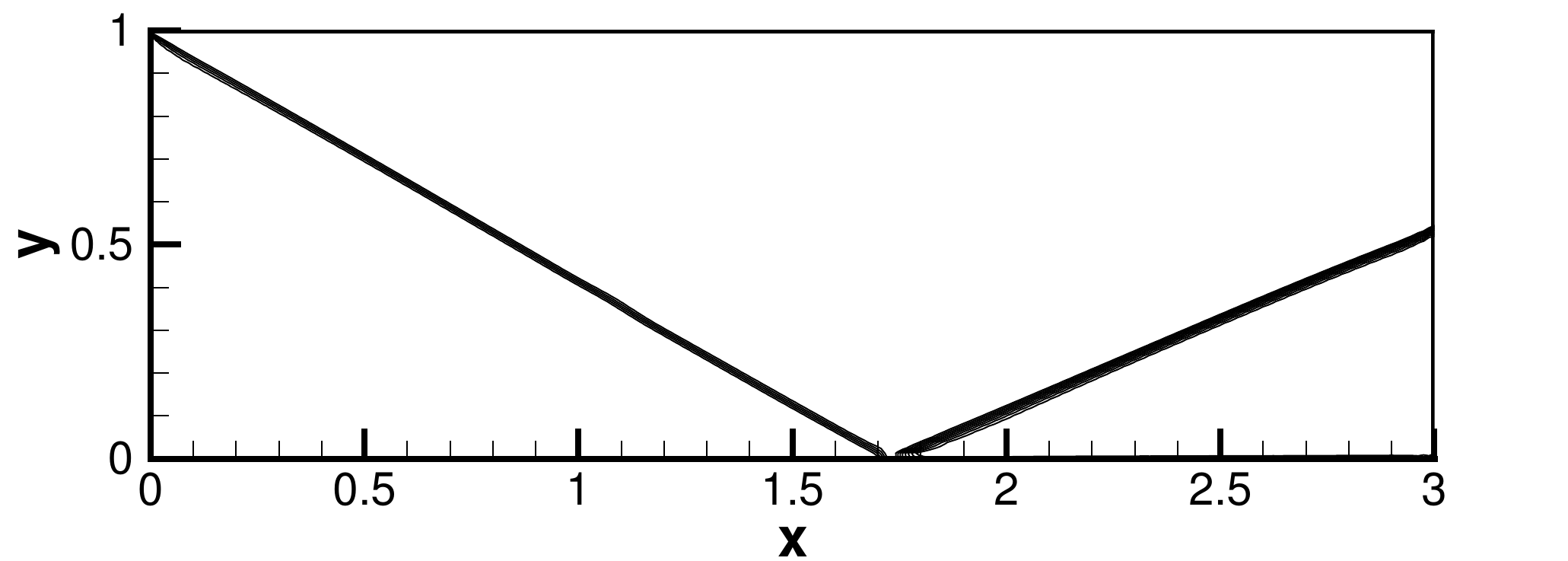}
\caption{Pressure contours (0.8:0.1:2.8) for $60 \times 20$, $120 \times 40$ and $240 \times 80$ quadrilateral mesh using Q4 element.}
\label{fig:SRq4}
\end{figure}

\begin{figure} [h!] 
\centering
\includegraphics[scale=0.6]{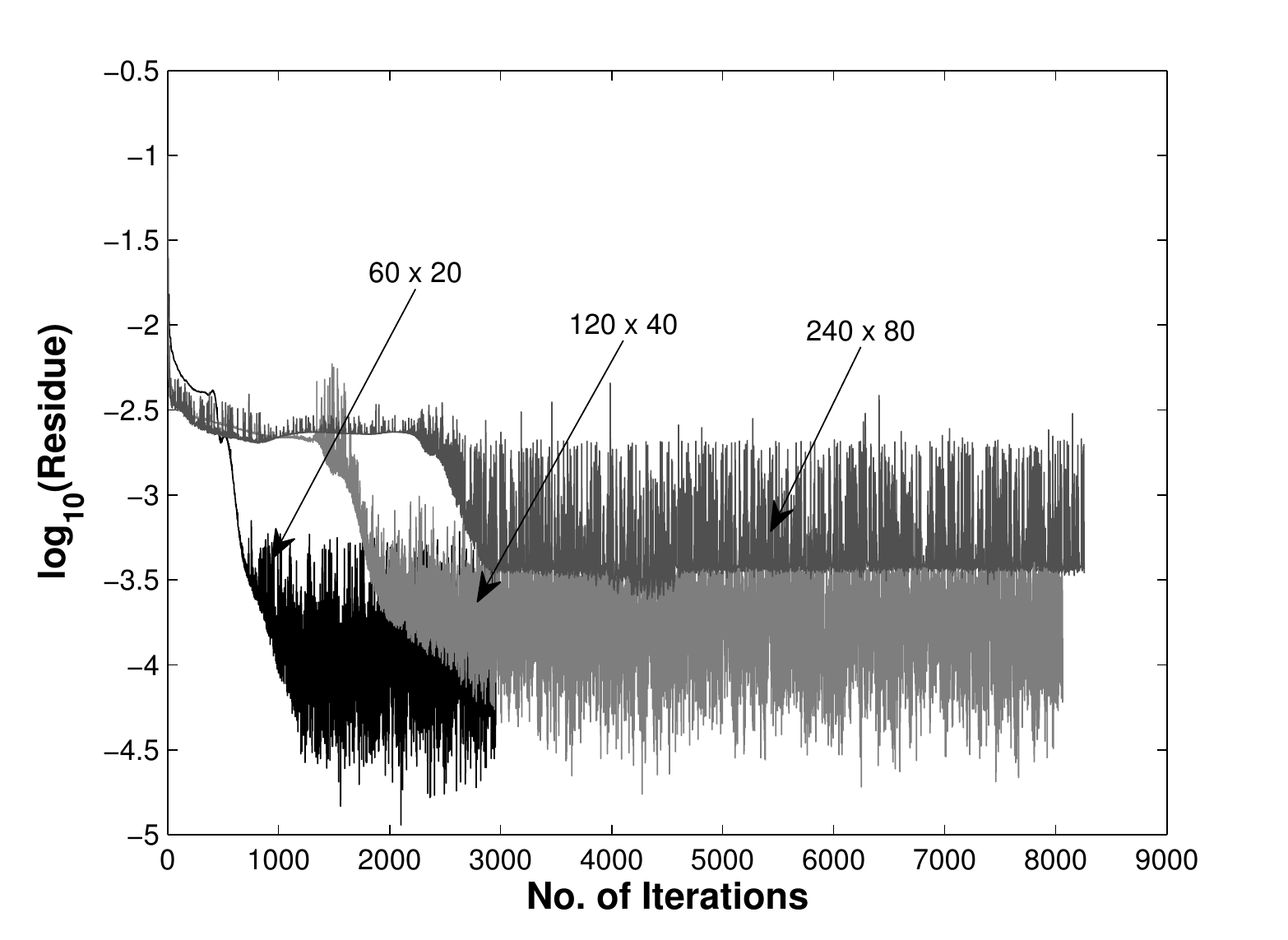}
\caption{Residue plots for $60 \times 20$, $120 \times 40$ and  $240 \times 80$ }
\label{fig:RS123Q}
\end{figure}

\begin{figure} [h!] 
\centering
\includegraphics[trim=1cm 0.9cm 1cm 12cm, clip=true, scale=0.6, angle = 0]{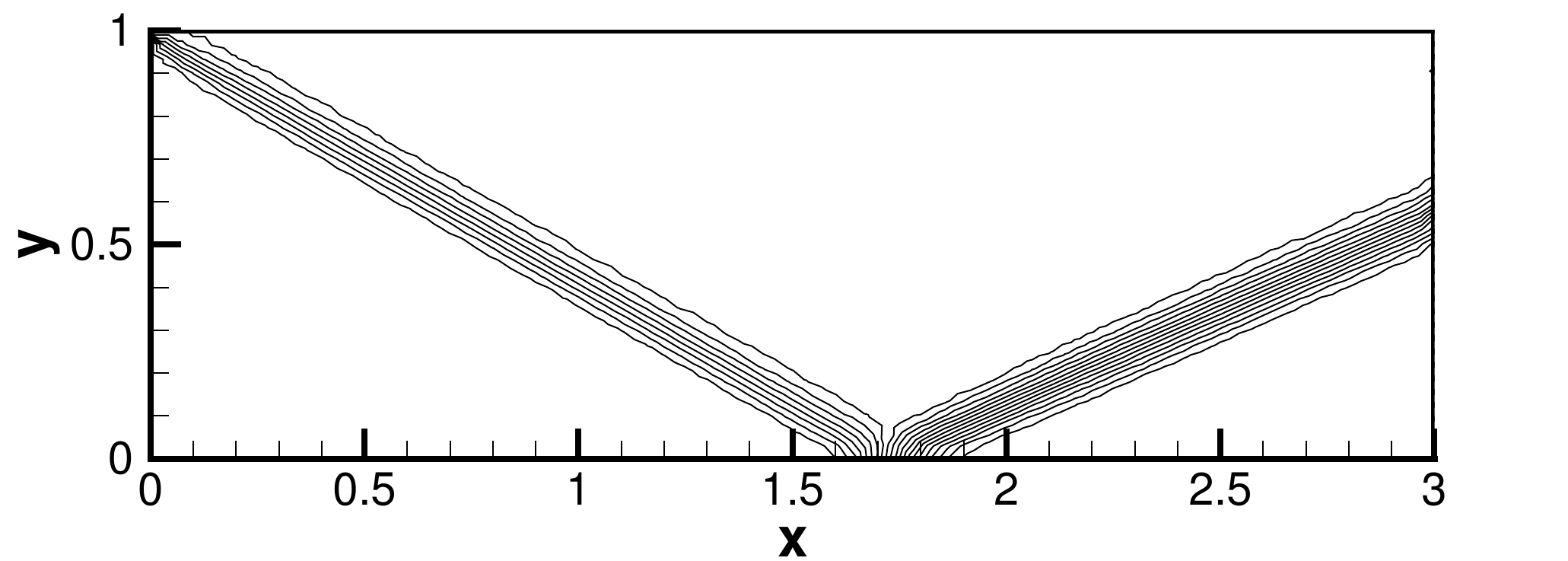}
\includegraphics[trim=1cm 0.9cm 1cm 12cm, clip=true, scale=0.6, angle = 0]{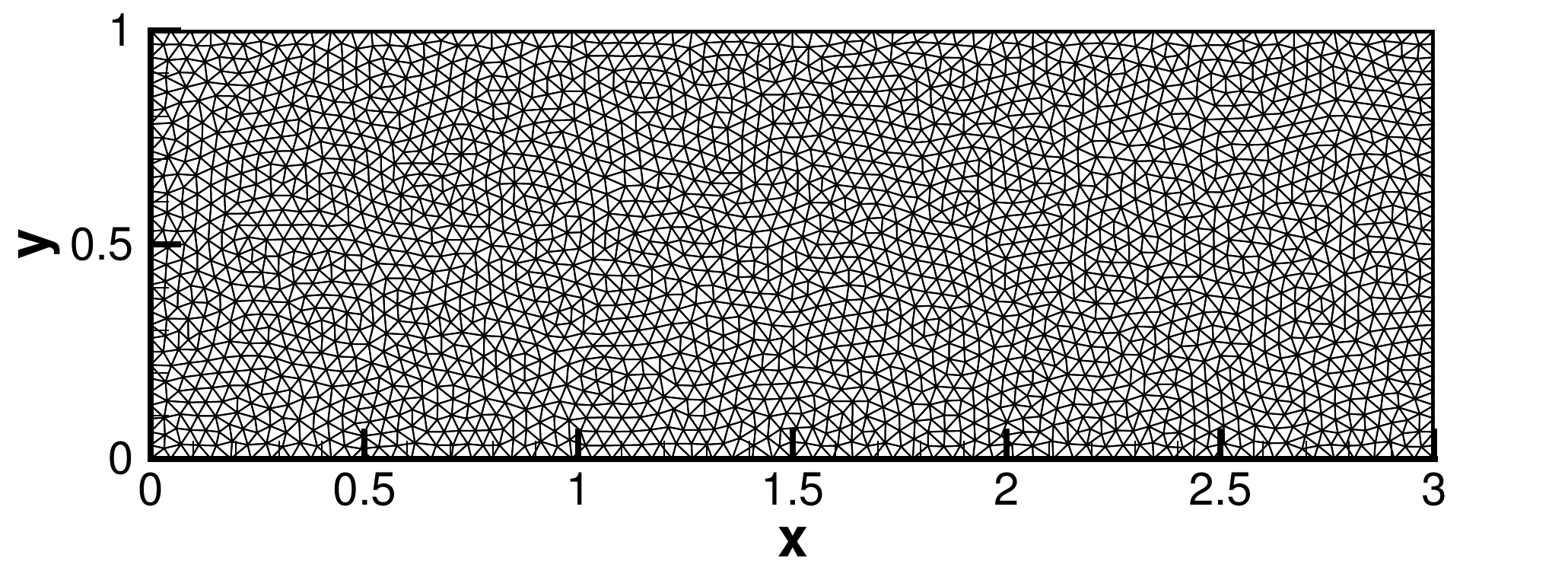}
\caption{Pressure contours (0.8:0.1:2.8) using T3 element.}
\label{fig:unsq34}
\end{figure}

\begin{figure} [h!] 
\centering
\includegraphics[scale=0.6]{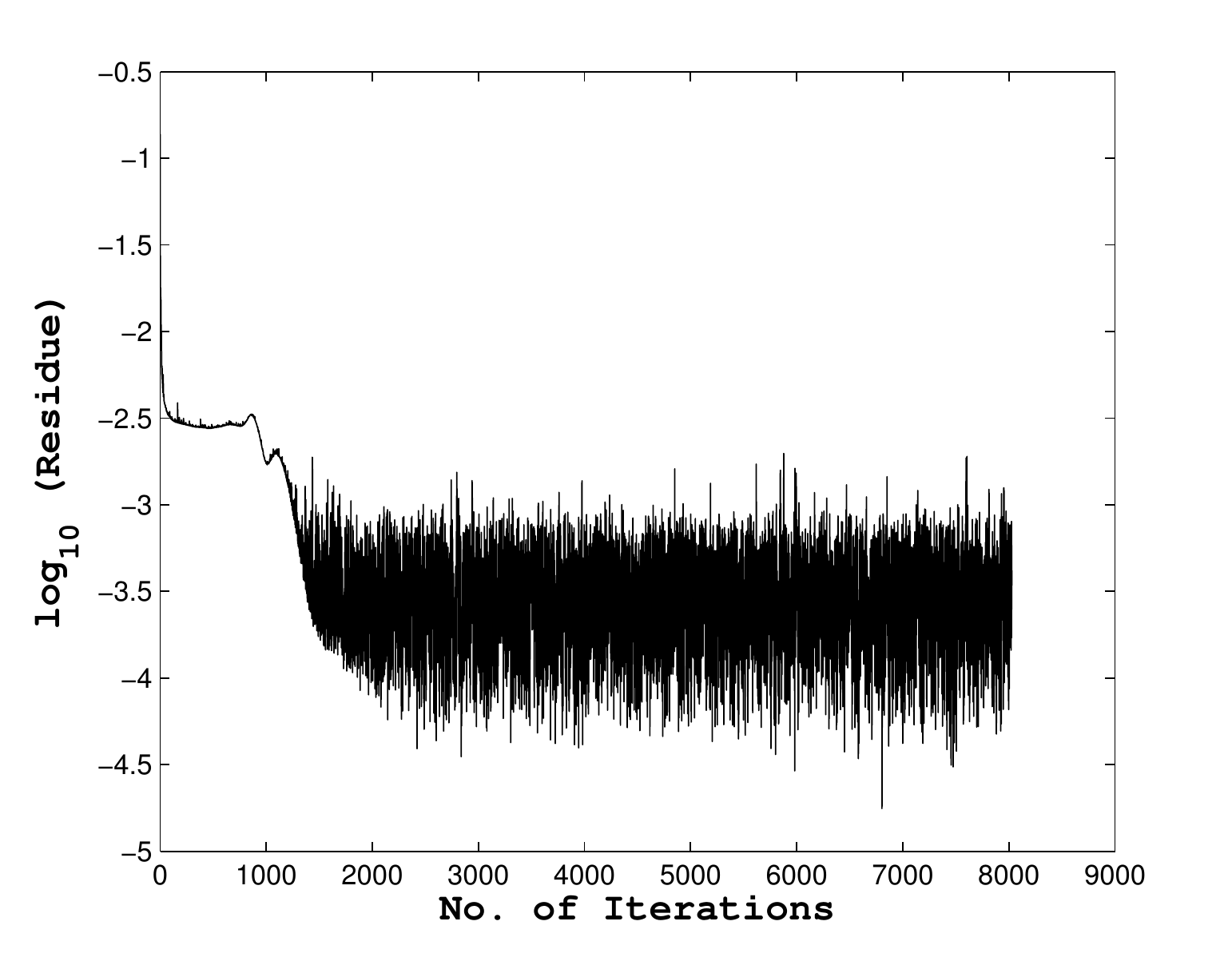}
\caption{Residue plots for a triangular elements.}
\label{fig:unsqRES}
\end{figure}

In this test case \cite{Yee_Warming_Harten} the domain is rectangular $[0 ,\,\,3] \times [0,\,\,1]$. The boundary conditions are
\begin{enumerate}
\item Inflow (left boundary) : $\rho =1, u=2.9, v=0, p=1/1.4$
\item Post shock condition (top boundary) : $\rho =1.69997, u=2.61934, v=-0.50633, p=1.52819$
\item Bottom boundary is a solid wall where slip boundary condition is applied, \textit{i.e.}, $\mathbf{v} .n = 0$.
\item At right boundary where the flow is supersonic all primitive variables $\rho$, $u $, $v$ and $p$ are extrapolated with first order approximation.
\end{enumerate}
Pressure plots for  $60 \times 20$, $120 \times 40$ and $240 \times 80$ quadrilateral mesh are given in figure ~\ref{fig:SRq4}.
The comparison of residue plots are given in figure ~\ref{fig:RS123Q}.  
For triangular unstructured mesh  (number of nodes: 2437 and number of triangles: 4680) the pressure contours are given in  figure ~\ref{fig:unsq34} and residue plot is shown in figure ~\ref{fig:unsqRES}.  
The incident and reflected shocks are captured quite accurately at correct positions.  

\subsubsection{Supersonic flow over a half cylinder: }
Two supersonic test cases with inflow Mach numbers 2 and 3 are tested on a half cylinder \cite{HVi}. The domain is half circular, the left outer circle is inflow boundary. Small circle inside the domain is a cylinder wall and the straight edges on right sides are supersonic outflow boundaries.    
\begin{figure} [h!] 
\centering
\includegraphics[trim=1cm 0.9cm 11cm 1cm, clip=true, scale=0.41, angle = 0]{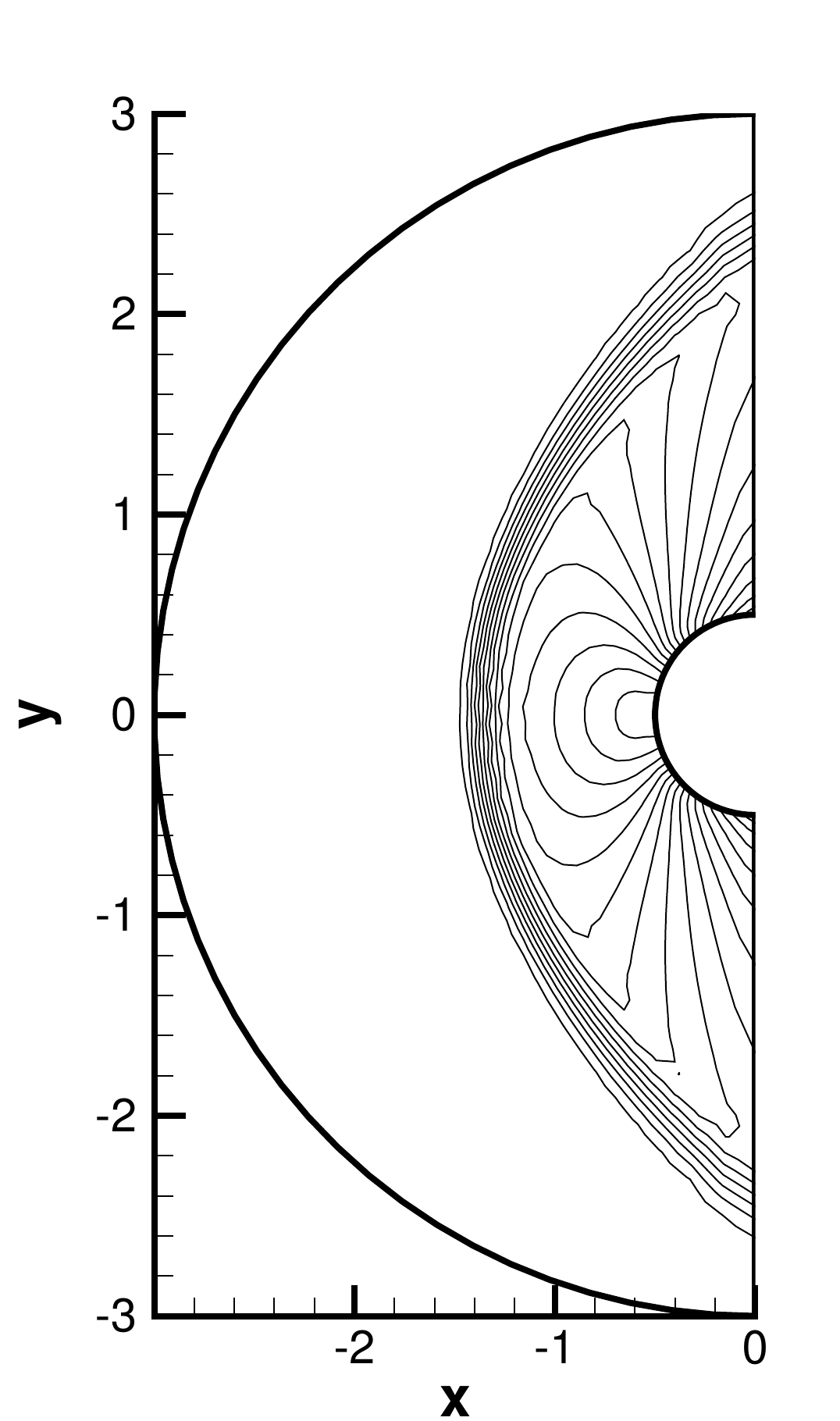}
\includegraphics[trim=1cm 0.9cm 11cm 1cm, clip=true, scale=0.41, angle = 0]{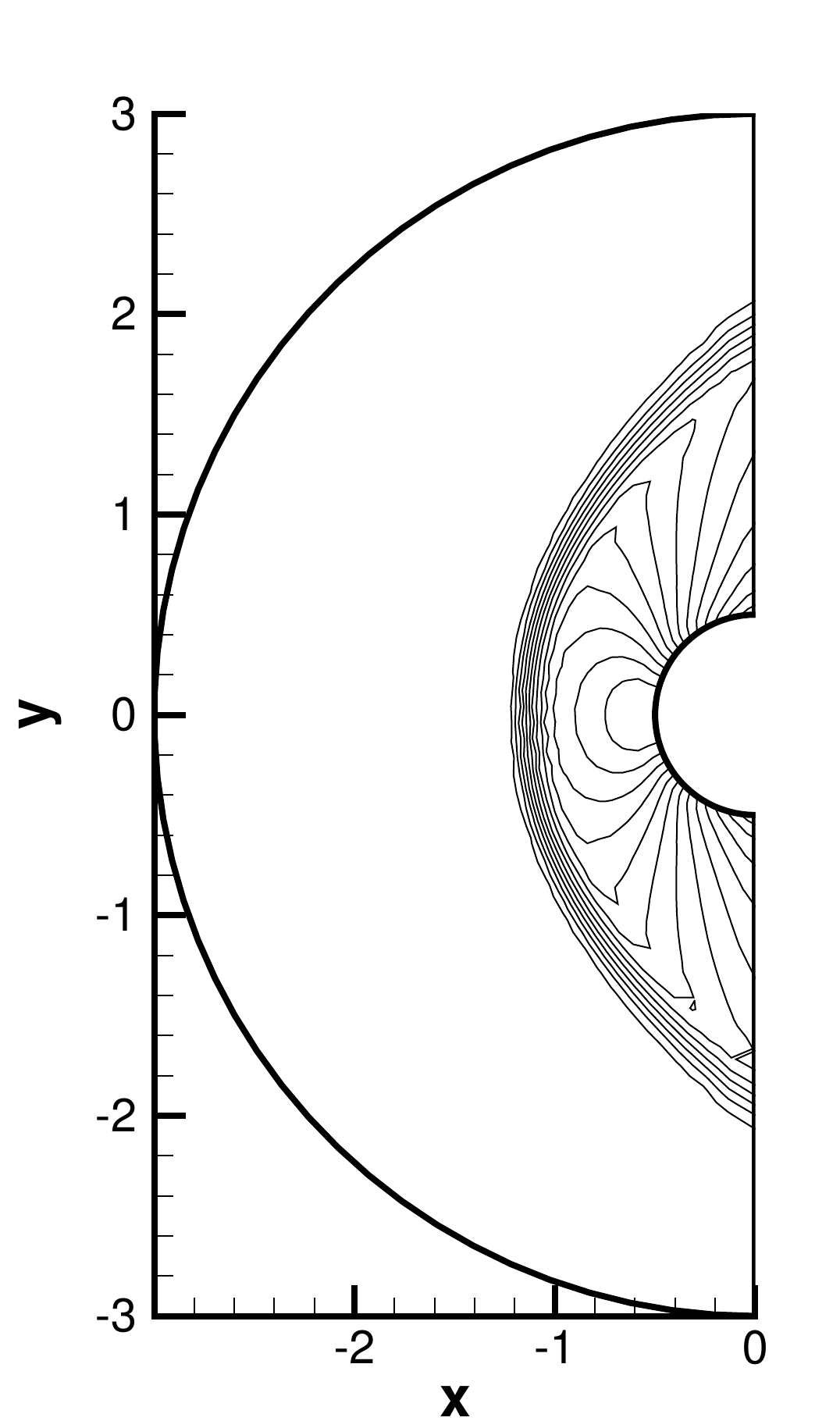}
\caption{Pressure contours (0.6:0.2:3.6) for Mach 2 (left) and (1:0.5:7) for Mach 3 (right) using $46 \times 46$ quadrilateral mesh.}
\label{fig:HalfCyl}
\end{figure}
Pressue plots (see figure ~\ref{fig:HalfCyl}) show that the bow shock in front of the half-cylinder is captured accurately at the right position in each case which are compared with existing results \cite{FSBillig}.

\subsection{Numerical Experiments for Implicit KSUPG}
All previously solved 1D test cases are again solved for implicit KSUPG method.

\subsubsection{Sod's Shock Tube Problem:}
\begin{figure} [h!] 
\centering
\subfigure
{\includegraphics[scale=0.42]{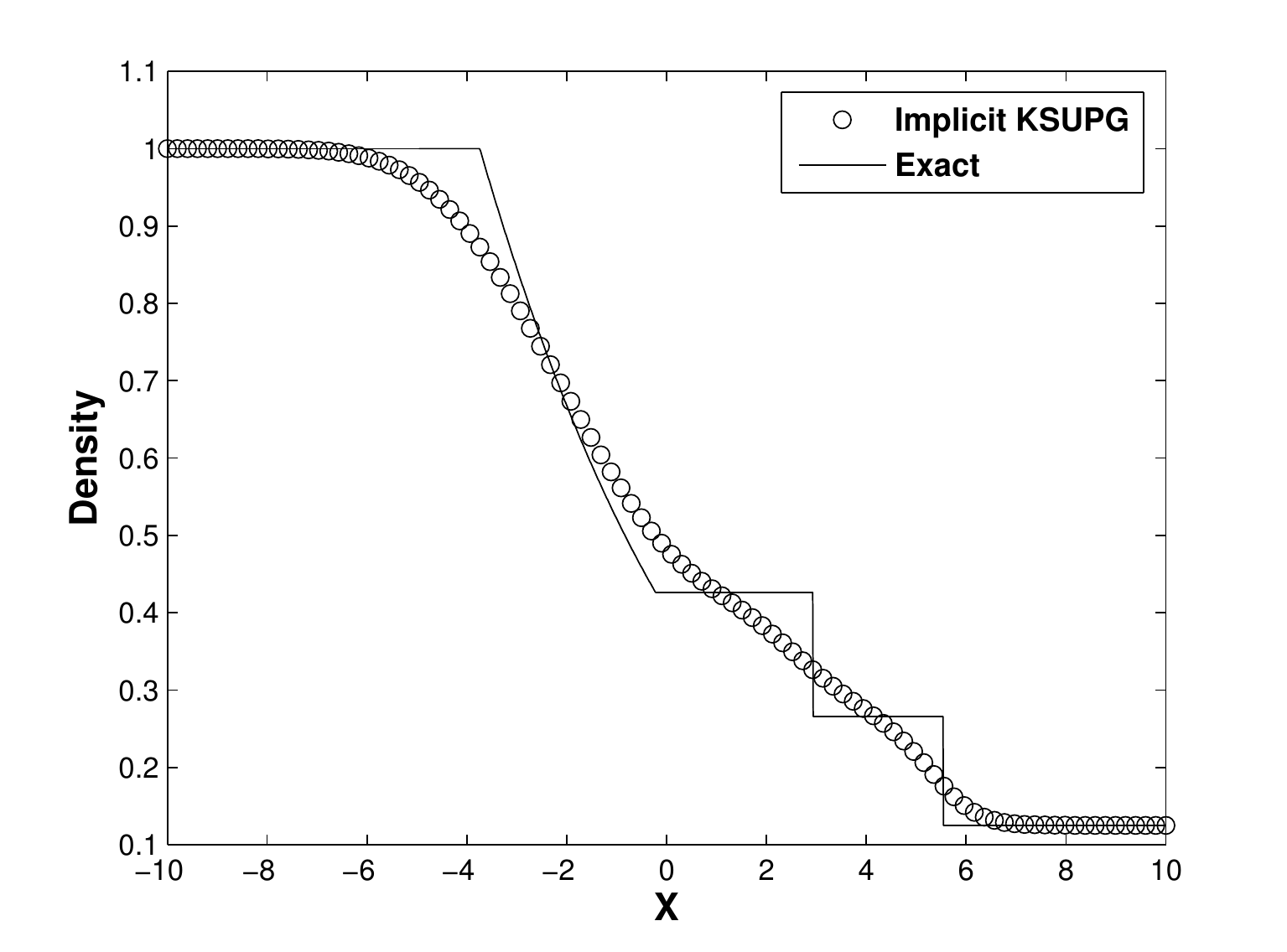}}
\subfigure
{\includegraphics[scale=0.42]{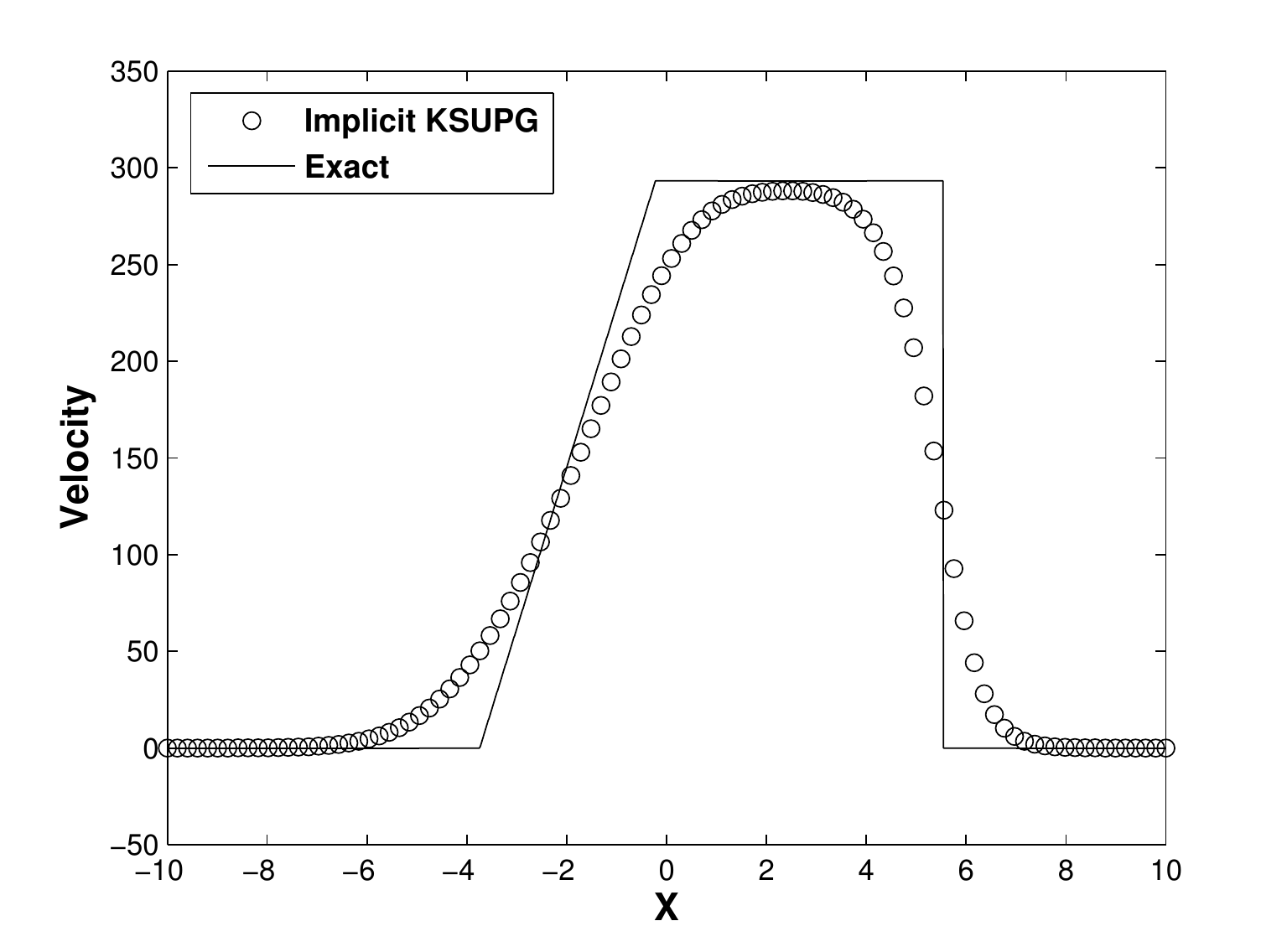}}
\subfigure
{\includegraphics[scale=0.42]{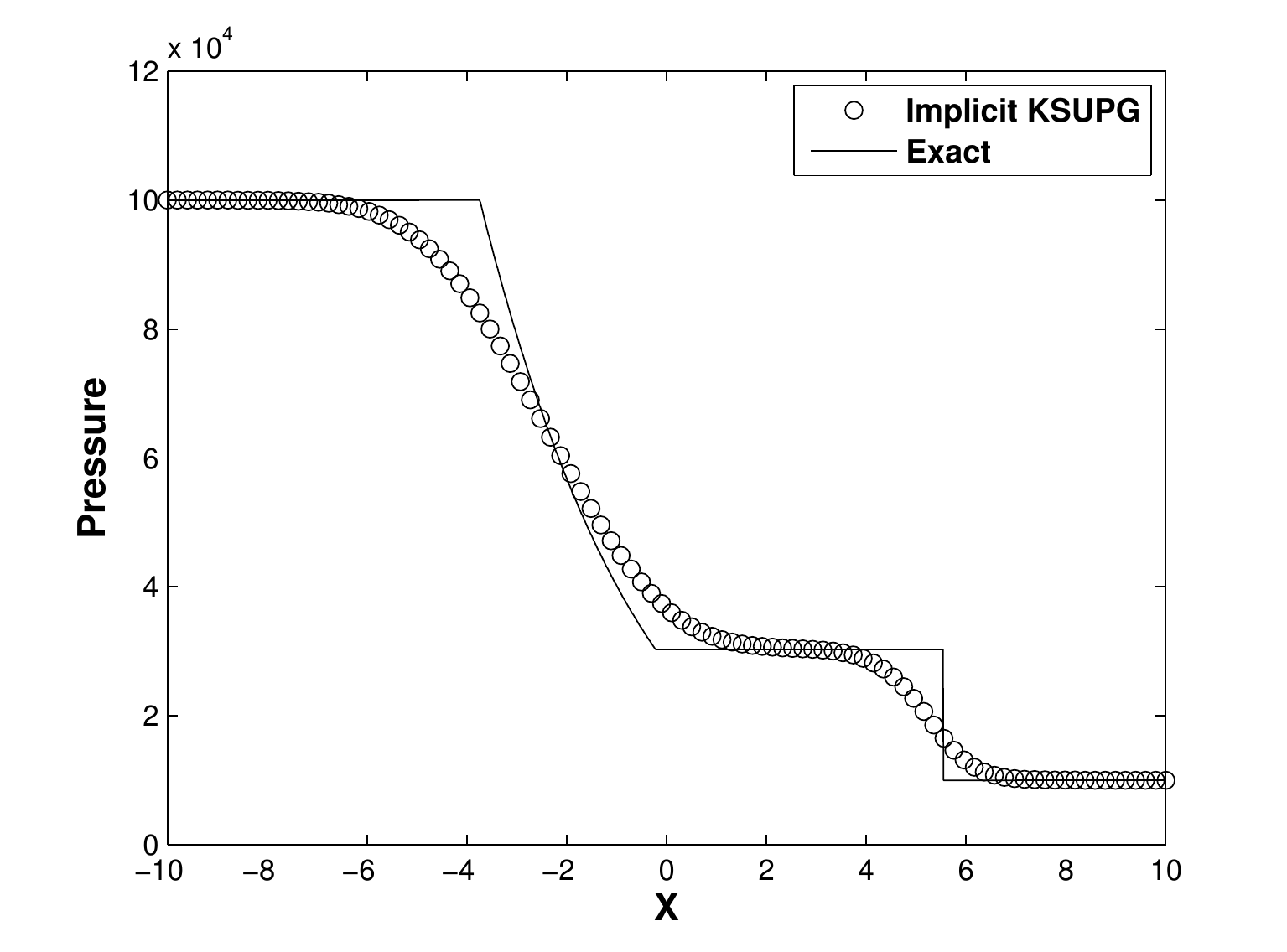}}
\subfigure
{\includegraphics[scale=0.42]{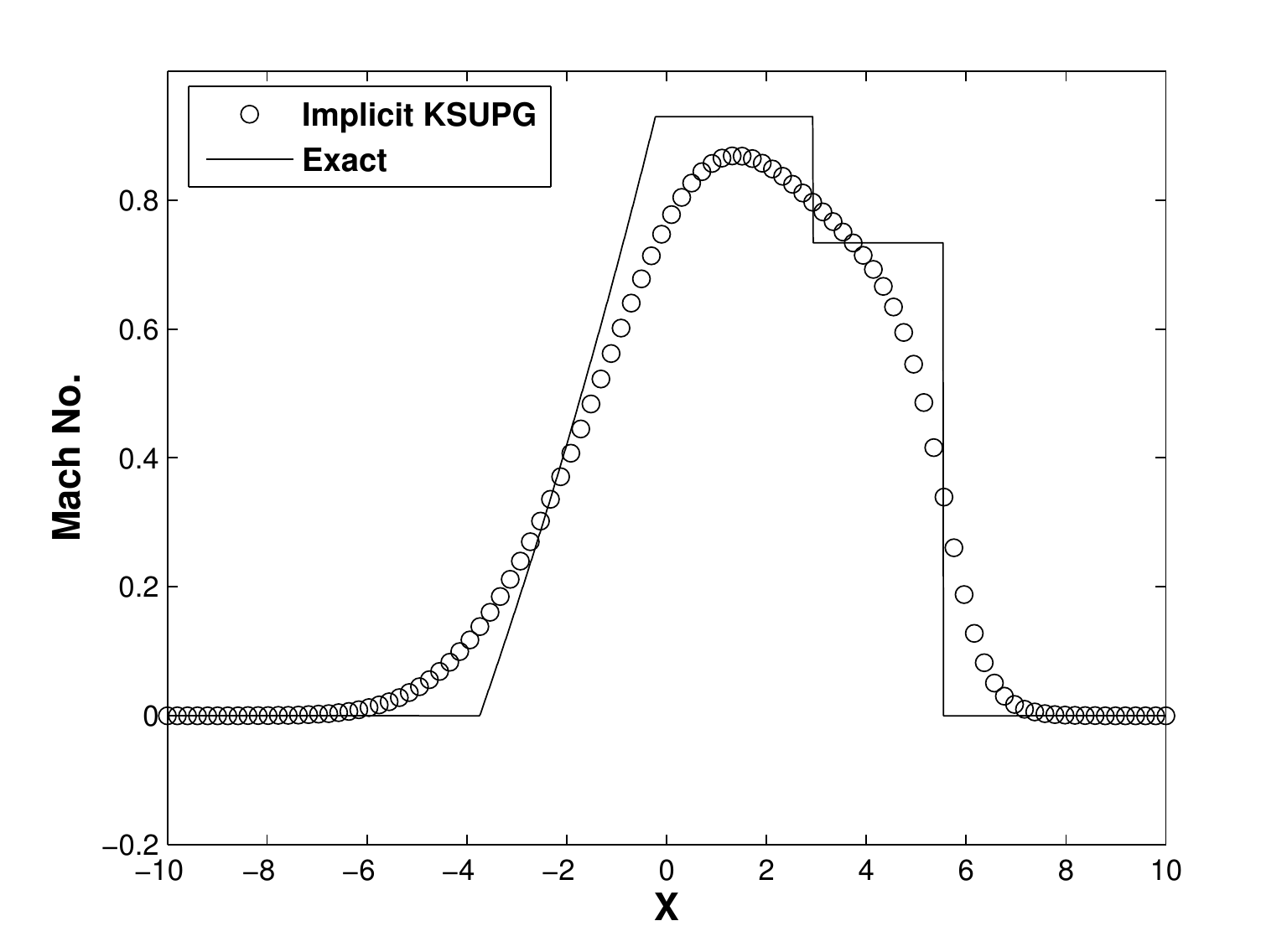}}
\caption{Sod's Shock Tube Problem}
\label{fig:kgE1II}
\end{figure}

The number of node points are 100 and CFL number is 0.6. Final time is t = 0.01. Figure ~\ref{fig:kgE1II} shows the density, velocity, pressure and mach number plots.

\subsubsection{Shock Tube Problem of Lax:}
\begin{figure} [h!] 
\centering
\subfigure
{\includegraphics[scale=0.42]{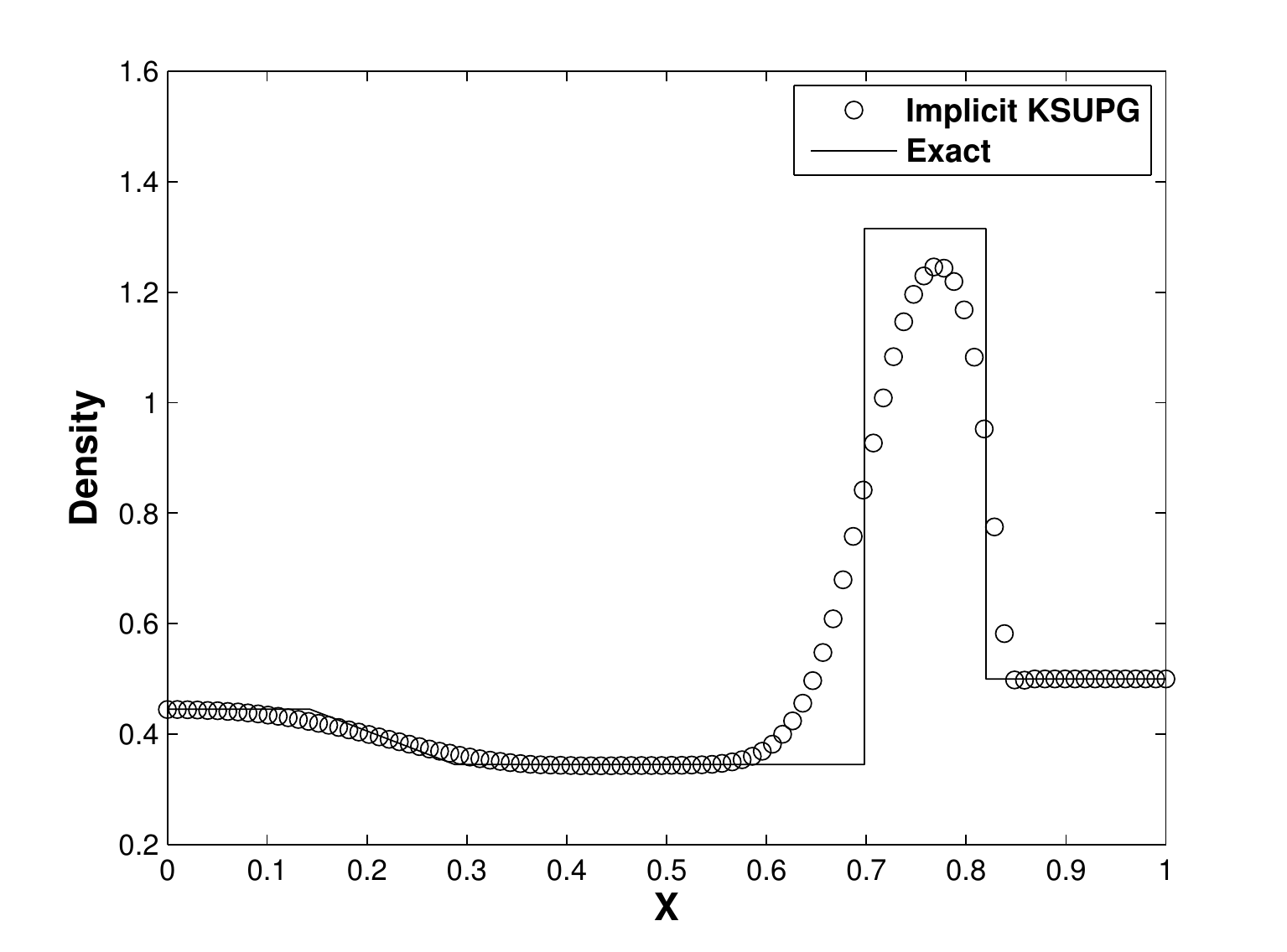}}
\subfigure
{\includegraphics[scale=0.42]{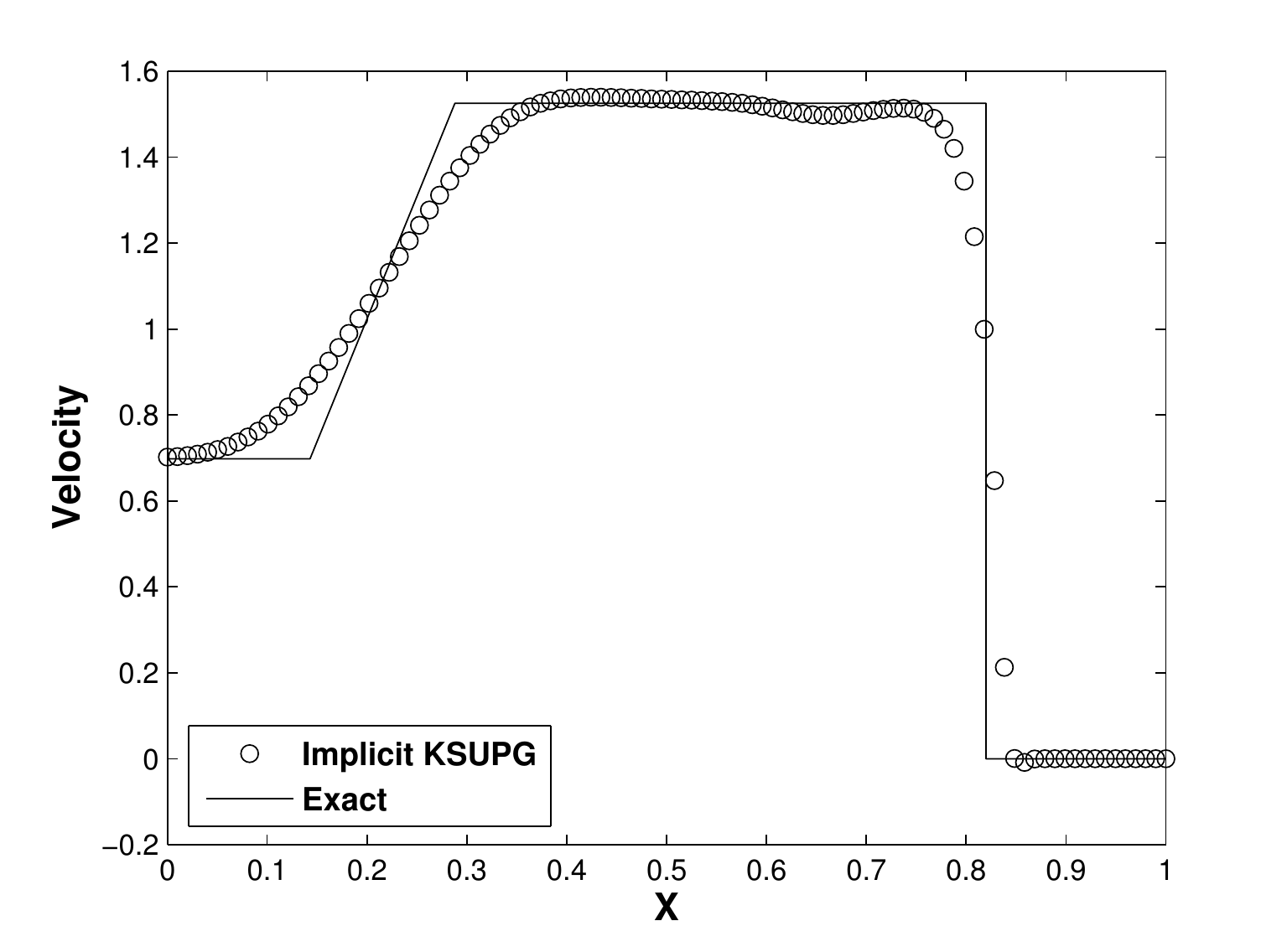}}
\subfigure
{\includegraphics[scale=0.42]{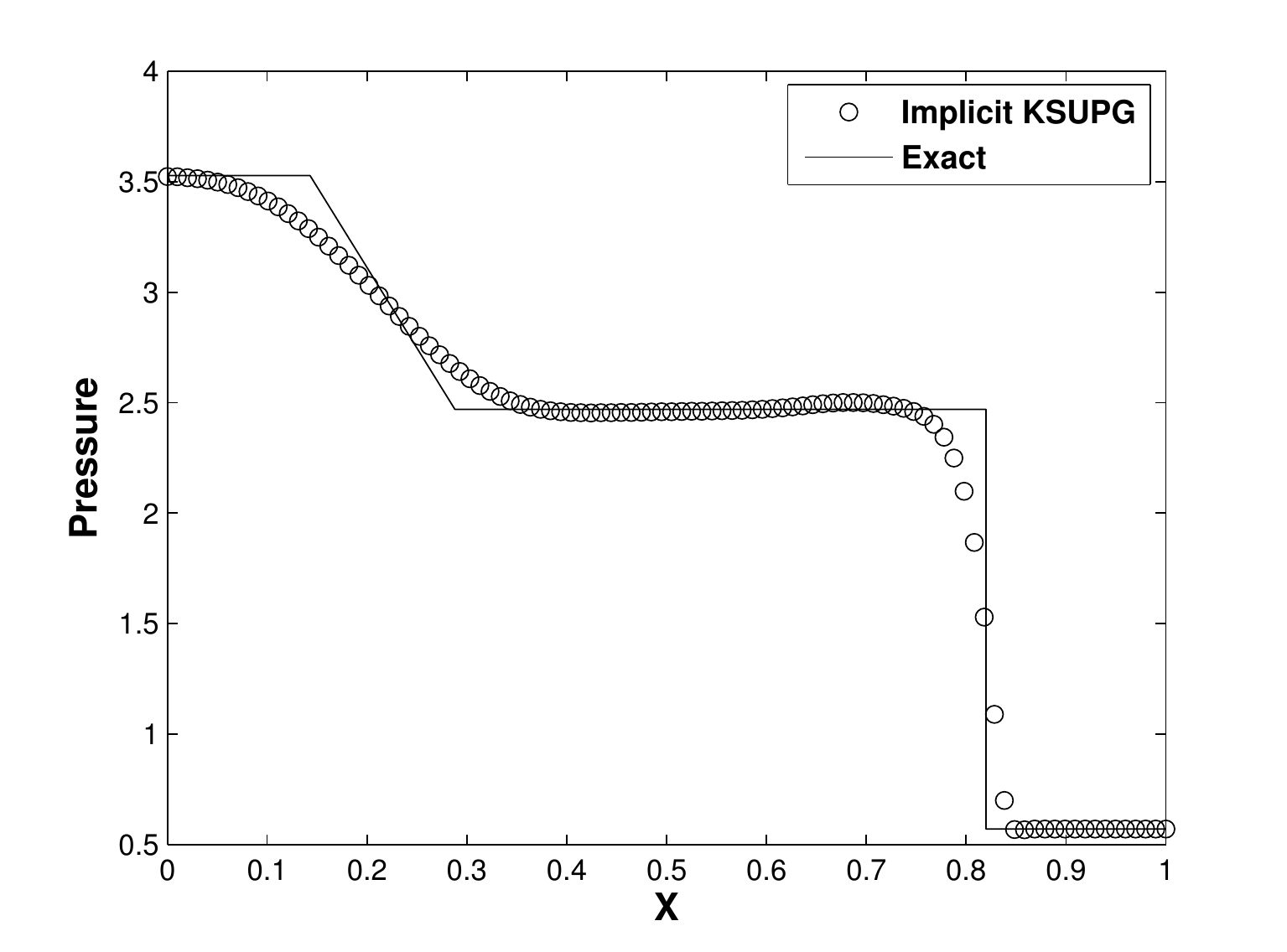}}
\subfigure
{\includegraphics[scale=0.42]{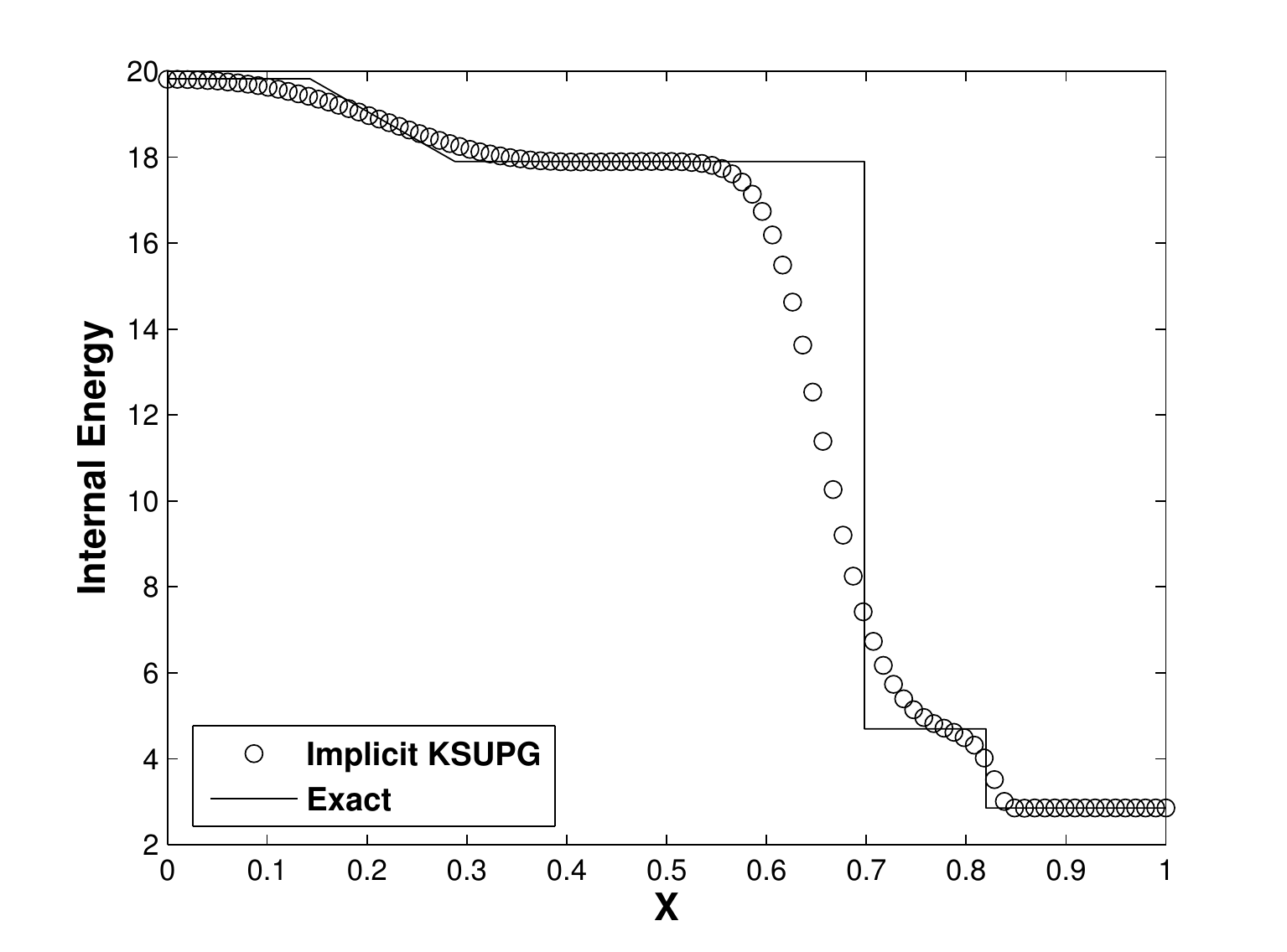}}
\caption{Shock Tube Problem of Lax}
\label{fig:kslaII}
\end{figure}

The number of node points are 100 and CFL number is 0.6.  Final time is t=0.13.  Figure ~\ref{fig:kslaII} shows the density, velocity, pressure and internal energy plots.

\subsubsection{Strong Rarefactions Riemann Problem:}
\begin{figure} [h!] 
\centering
\subfigure
{\includegraphics[scale=0.42]{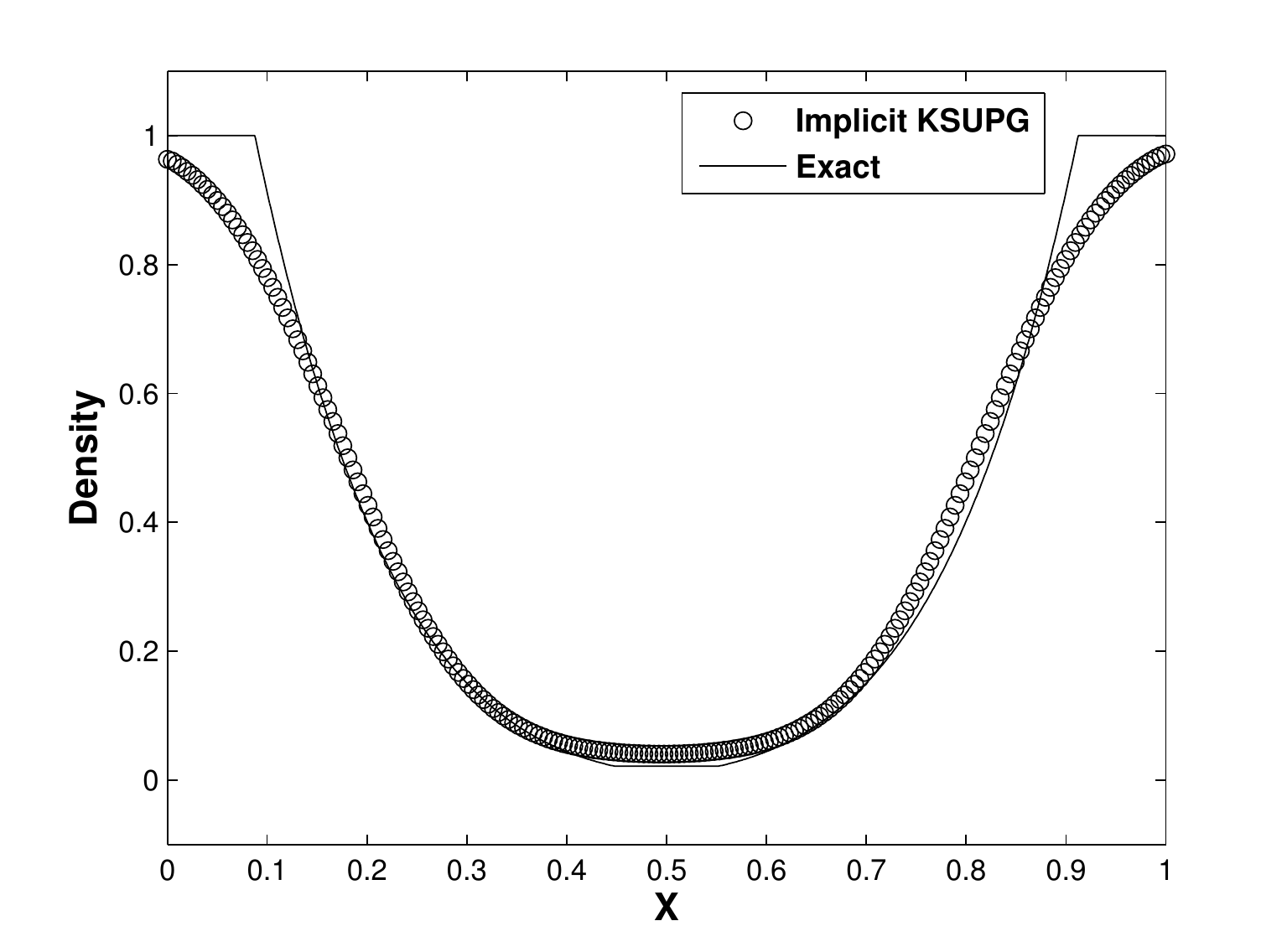}}
\subfigure
{\includegraphics[scale=0.42]{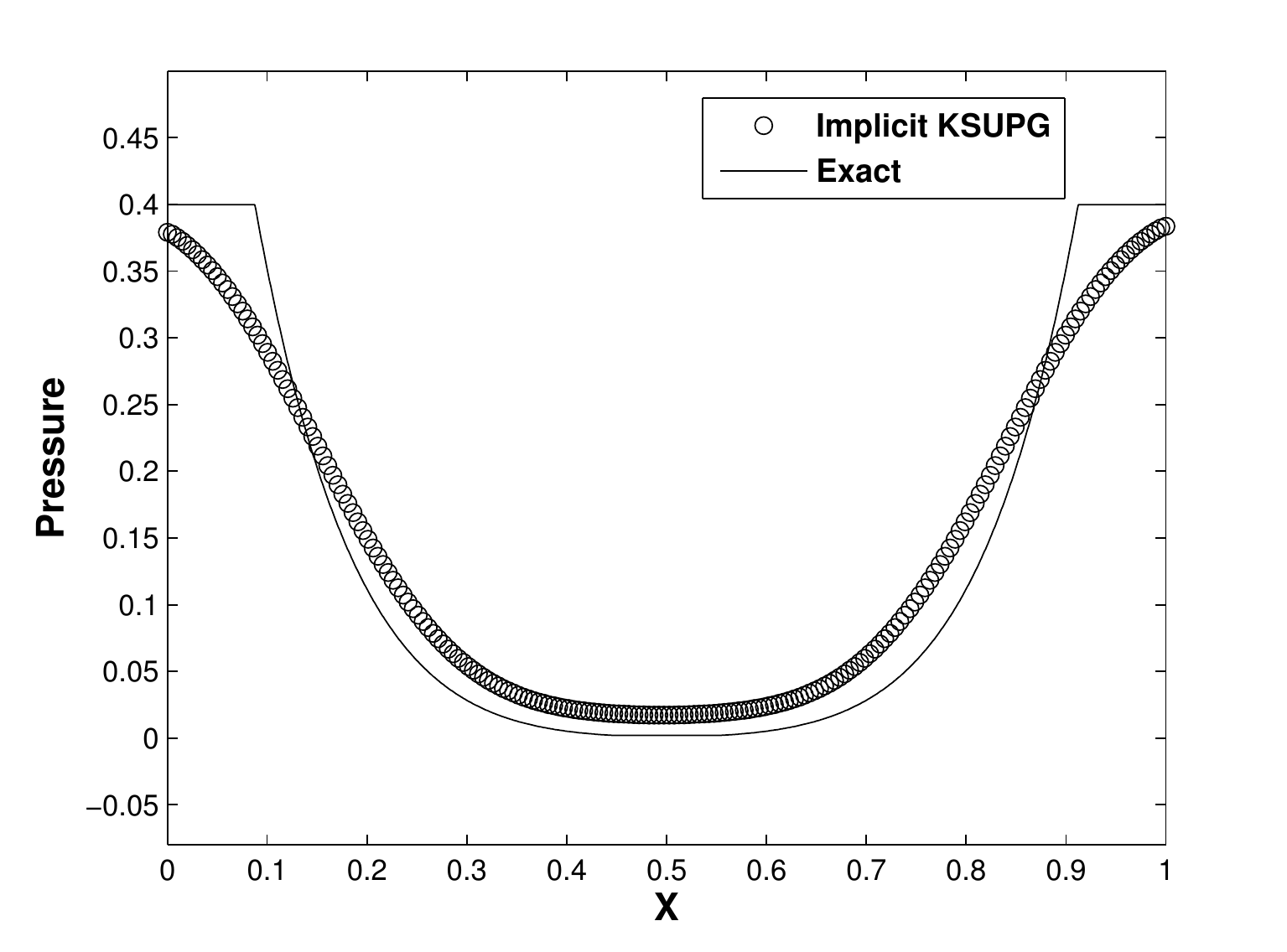}}
\subfigure
{\includegraphics[scale=0.42]{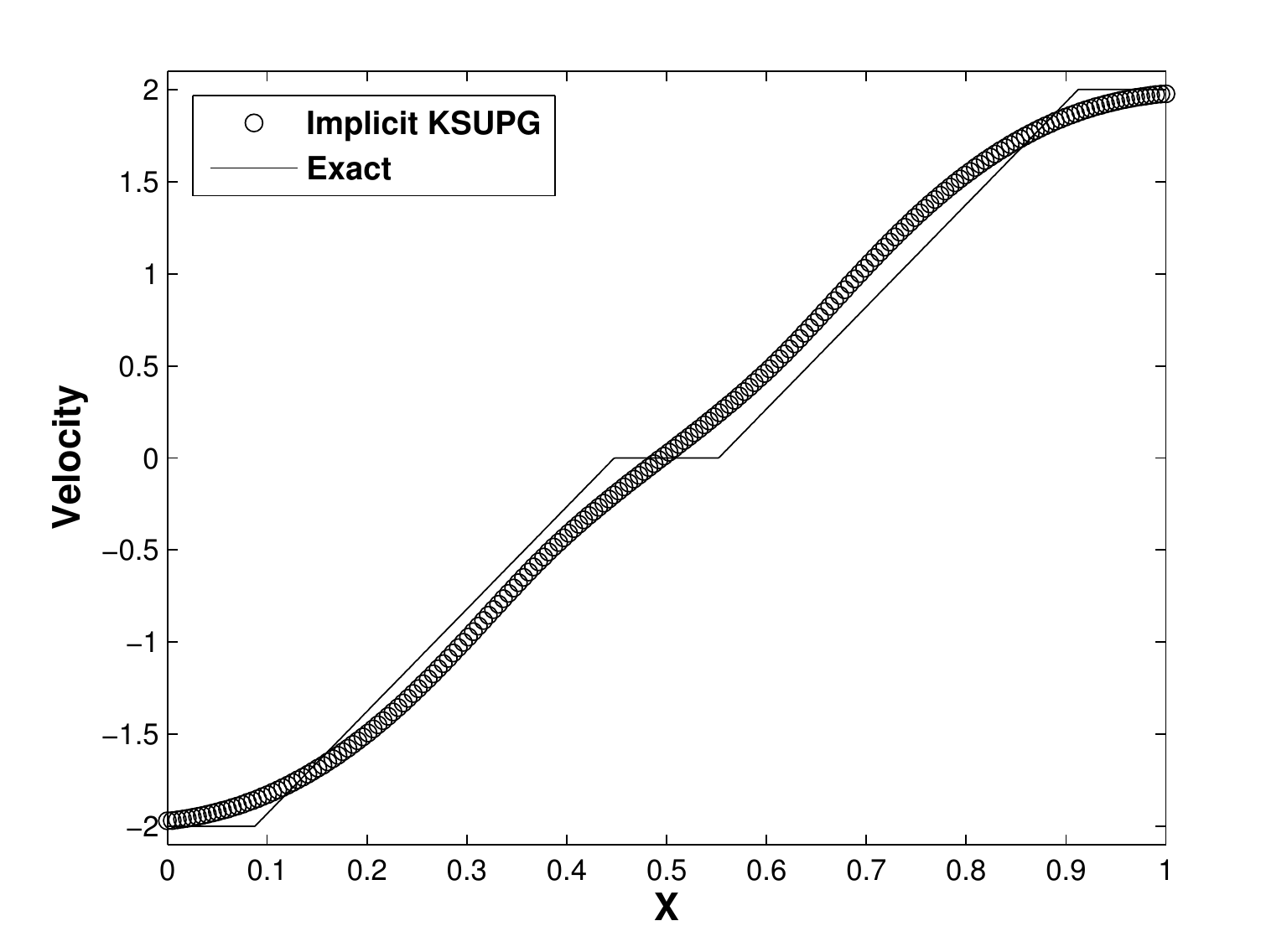}}
\caption{Strong Rarefactions Riemann Problem}
\label{fig:ksrrII}
\end{figure}

 The number of node points are 200 and CFL number is 0.6. Final time is t=0.15.  Figure ~\ref{fig:ksrrII} shows the density,  pressure and velocity plots.    

\subsection{Comparison of Explicit and Implicit KSUPG scheme for 2D Euler test case}
\begin{figure} [h!] 
\centering
\includegraphics[trim=1cm 0.9cm 1cm 12cm, clip=true, scale=0.6, angle = 0]{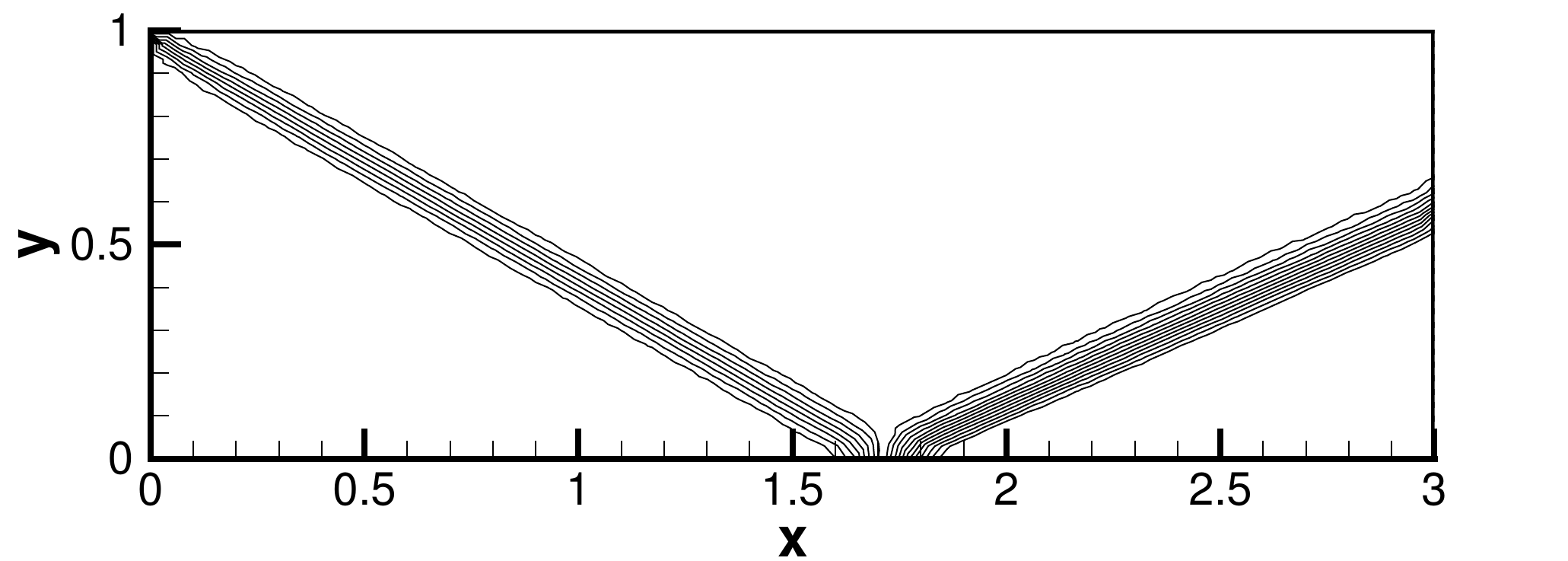}
\includegraphics[scale=0.55]{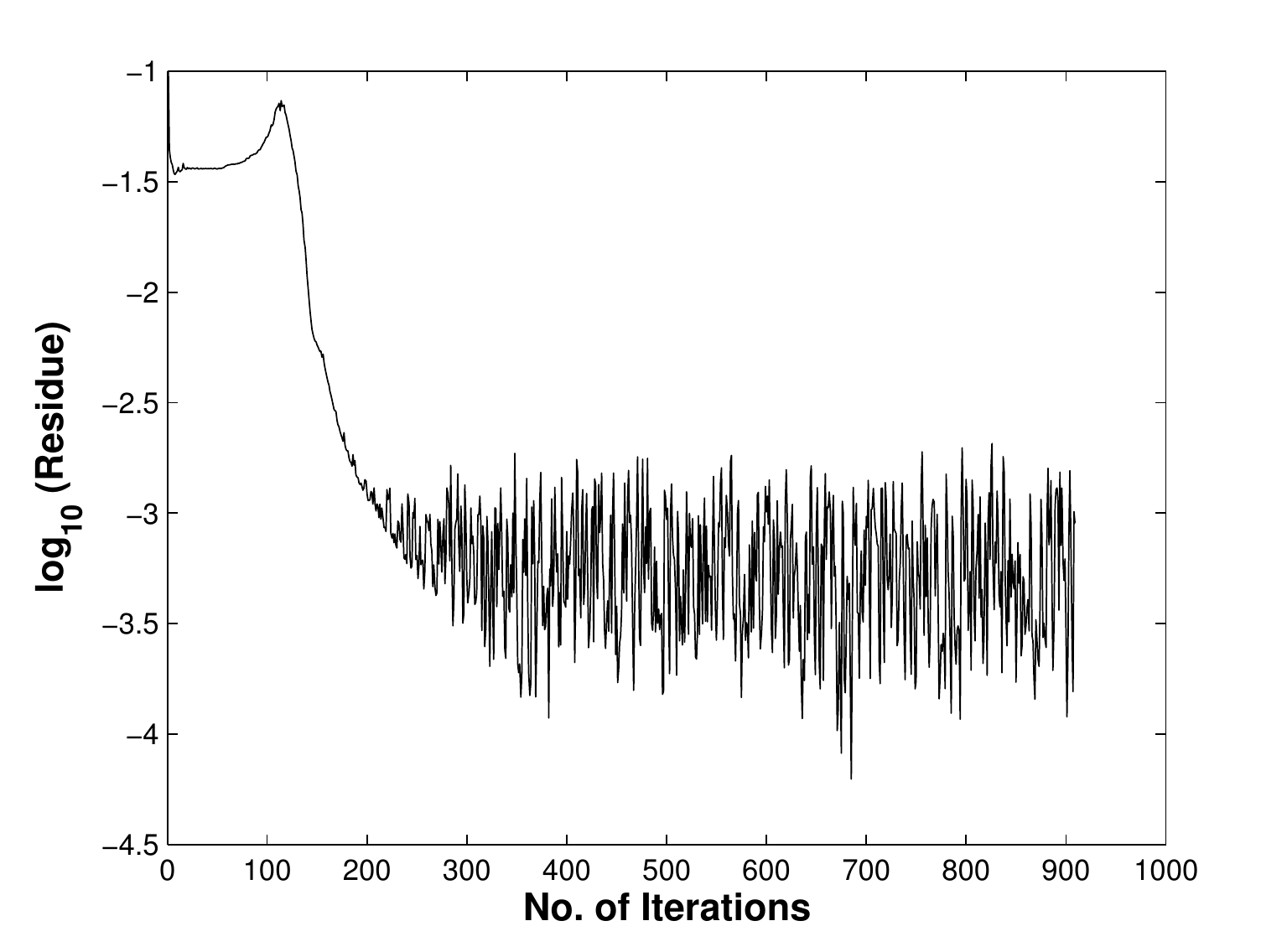}
\caption{Pressure plot (0.8:0.1:2.8) for shock reflection test case using Implicit KSUGP on $120 \times 40$ mesh and the residue plot.}
\label{fig:SRimK}
\end{figure}

In case of implicit KSUPG scheme, 2D Euler test cases are solved using CFL =1.  
For the comparison of explicit and implicit KSUPG schemes shock reflection test case is solved using different grids. Figure ~\ref{fig:SRimK} shows the pressure contour plot on $120 \times 40$ Q4 grid and the residue vs number of iterations.

Iteration speed-up ratio is defined as the ratio of number of iterations required for the residue to drop below a predefined tolerance value for an explicit method to that of an implicit method.  Similarly, one can define the computational speed-up ratio which is the ratio of total computational time required for explicit method to that of implicit method.  Following tables shows comparison of computational cost and number of iterations taken for explicit and implicit KSUPG schemes for oblique shock reflection test case. 

\begin{center} 
\small \begin{tabular}{|c|c|c|c|} \hline
& \multicolumn{3}{c|} {Total  Computational   Cost} \\ \hline 
&&& \\ 
 Grid Size & $60 \times 20$&  $120 \times 40$ &$240 \times 80$  \\ 
&&& \\ \hline
&&& \\ 
Explicit  &8388 sec &69948 sec &  189000 sec\\ 
KSUPG &&& \\ \hline
&&& \\ 
Implicit   &1338 sec &8784 sec &  37288 sec\\ 
KSUPG &&&\\ \hline 
&&& \\ 
Computational  &6.26 &7.96 &  5.06\\ 
Speed-up Ratio &&&\\ \hline  
 \end{tabular}
\end{center}

\begin{center}
\small \begin{tabular}{|c|c|c|c|} \hline
&&&\\ 
 Grid Size & $60 \times 20$&  $120 \times 40$ &$240 \times 80$  \\ 
&&& \\ \hline
&&& \\ 
Tolerance & $10^{-3.5}$&  $10^{-3.5}$ &$10^{-3.5}$   \\ 
Value &&& \\ \hline
&&& \\ 
No. of Iterations taken  & 749 &1720 &3862 \\ 
for Explicit KSUPG &&& \\ \hline
&&& \\ 
No. of Iterations taken  & 185 &310&504 \\ 
for Implicit KSUPG &&& \\ \hline 
&&& \\ 
 Iteration & 4.07  &5.54 & 7.66   \\ 
Speed-Up Ratio &&& \\ \hline 
 \end{tabular}
\end{center}

\subsubsection{Oblique shock}
Oblique shock test case is solved using implicit KSUPG method.
\begin{figure}[h!]
\centering
\includegraphics[scale=0.3]{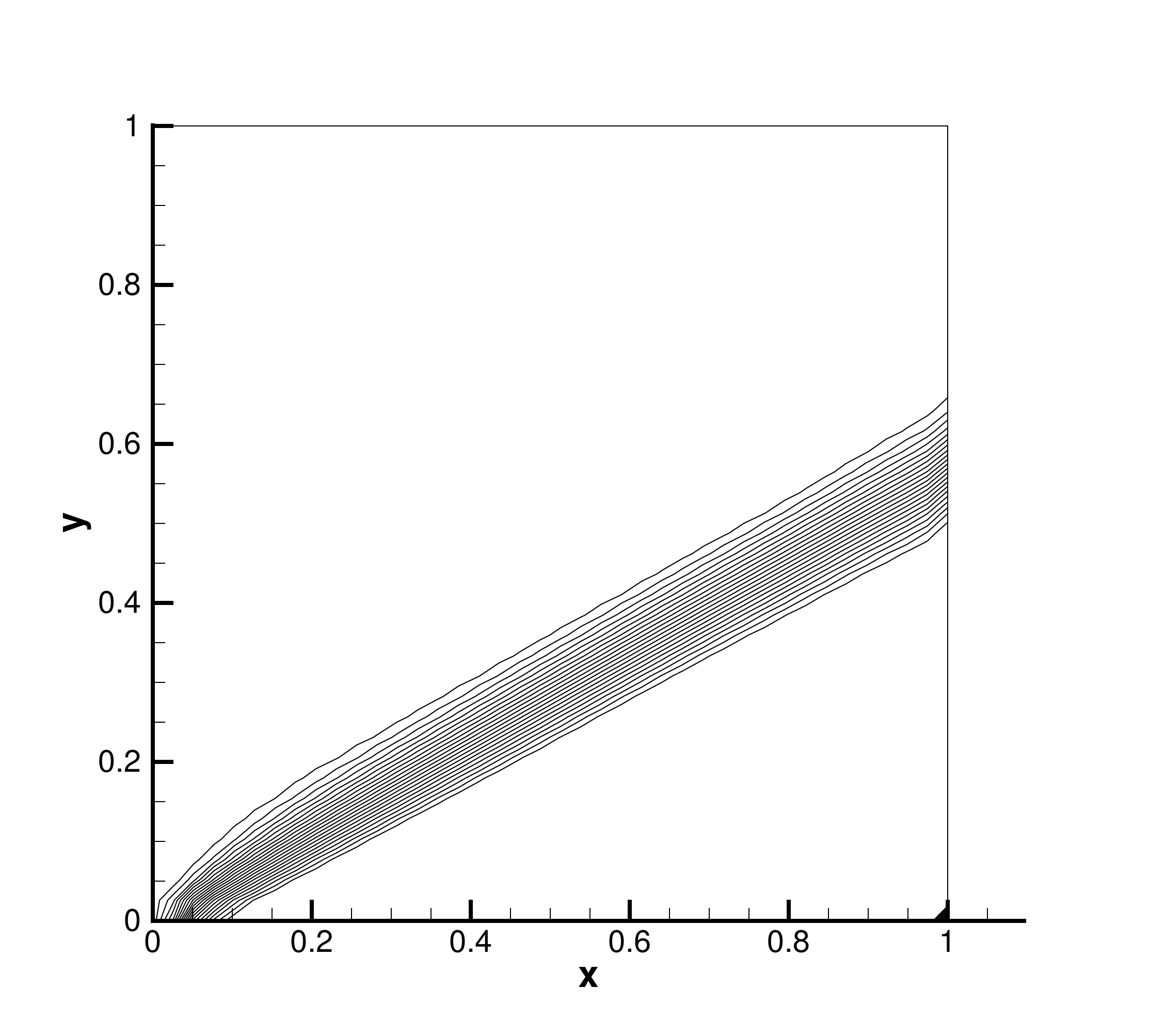}
\caption{Pressure contours (0.18:0.05:0.29) using $40 \times 40$ Q4 elements}
\label{fig:oSI312im}
\end{figure}
 Figure  ~\ref{fig:oSI312im} shows the pressure contours using $40 \times 40$ Q4 mesh.

\subsection{Sparsity Pattern in Explicit and Implicit KSUPG Method}
\begin{figure} [h!] 
\centering
\includegraphics[scale=0.35]{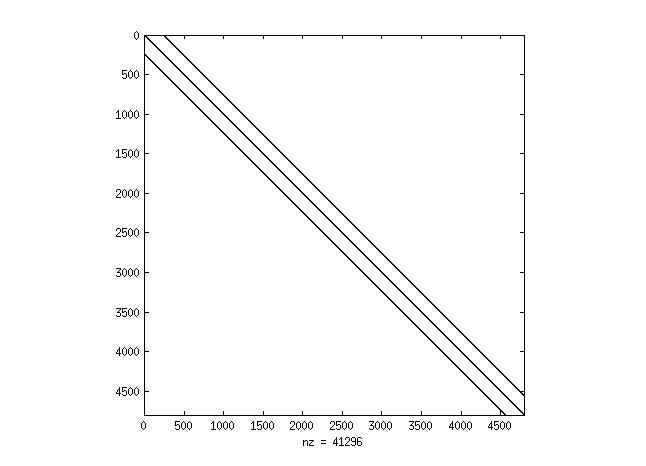}
\includegraphics[scale=0.35]{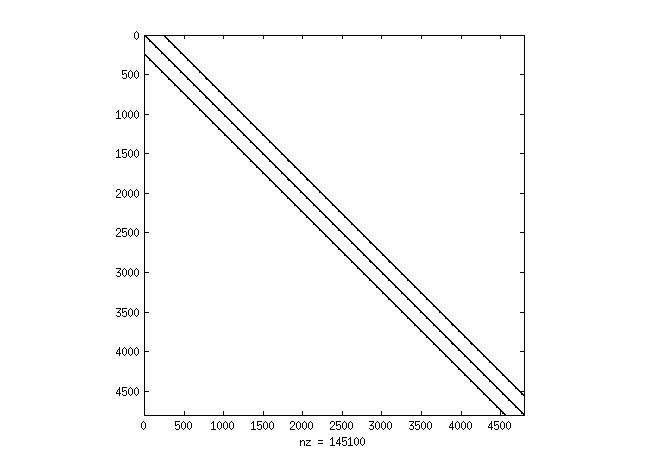}
\caption{Sparsity pattern in explicit and implicit KSUPG}
\label{fig:SPPP}
\end{figure}

Figure ~\ref{fig:SPPP} shows the sparsity patterns for explicit and implicit KSUPG method for shock reflection test case over $60 \times 20$ mesh size using Q4 elements. Both matrices are unsymmetric. The number of nonzero entries in explicit method is 41296 and in case of implicit method it is 145100 which are far more than the previous case.
The half band-width of explicit KSUPG matrix is 245 whereas for implicit KSUPG it is 248. In case of implicit KUPG method the condition number of coefficient matrix of assembled system is high compared to explicit KSUPG scheme. Following table shows the condition number calculated using $L_2$ norm for different grid size for shock reflection test case.

\begin{center}
\small \begin{tabular}{|c|c|c|} \hline
& \multicolumn{2}{c|} {Condition Number} \\ \hline 
&& \\ 
 Grid Size & Implicit KSUPG (CFL =1)&  Explicit KSUPG  \\ 
&& \\ \hline
&& \\ 
$60 \times 20$  &2.1670e04 &8.3722e03\\  
&& \\ \hline
&& \\ 
$120 \times 40$  &9.4669e04 &3.4646e04\\  
&& \\ \hline
&& \\
$240 \times 80$   &4.18822e05&1.40931e05\\ 
 &&\\ \hline 
 \end{tabular}
\end{center}

\section{Conclusions}

In this paper we presented a novel explicit as well as implicit kinetic theory based streamline upwind Petrov Galerkin scheme (KSUPG scheme) in finite element framework for both scalar case (inviscid Burgers equation in 1D and 2D) and vector case (1D and 2D Euler equations of gas dynamics). The proposed numerical scheme is simple and easy to implement. The important advantage in using a kinetic scheme in finite element framework is that, instead of dealing with nonlinear hyperbolic conservation laws, one needs to deal with a simple linear convection equation.  In comparison with the standard SUPG scheme, the advantage of the proposed scheme is that, it does not contain any complicated expression for $\tau$ (which is the intrinsic time scale) especially for vector equations. Also, for the multidimensional Burgers equation, standard SUPG scheme requires additional diffusion term (shock capturing parameter) which is not required in the proposed scheme. The accuracy and robustness of the scheme is demonstrated by solving various test cases for Burgers equation and Euler equations. Spectral stability analysis is done for 2D linear equation which gives an implicit expression of stable time step.  Finally, comparison between explicit and implicit versions of KSUPG scheme is done with respect to the number of iterations, computational cost, sparsity pattern and the condition number of a global system of equations.


\end{document}